\documentclass{article}
\usepackage{amsmath,amssymb}
\usepackage[arc,ps,dvips,curve]{xy}
\usepackage{epsf}
\UseDVIPSspecials\relax
\def\connsumtransversevertex#1#2#3{\setbox0=\hbox{$\scriptstyle\scriptstyle\rightarrow$}#1\#_{\vtop{\hbox to \wd0{\hfil$\scriptstyle #3$\hfil}\kern-7pt\hbox{$\scriptstyle\rightarrow$}}} #2}
\def\connsumparallelvertex#1#2#3{#1\#_{#3\uparrow} #2}
\def\set#1\endset{\{\,#1\,\}}
\def\proof{\ifdim\lastskip<\smallskipamount\relax\removelastskip
  \vskip\smallskipamount\fi\leavevmode\hbox to 10pt{\hfil}{\it Proof. }}
\def\strutdepth{\dp\strutbox}
\def\epmarker{\vbox to \strutdepth{\baselineskip\strutdepth\vss\hfill{%
\hbox to 0pt{\hss\vrule height 4pt width 4pt depth 0pt}\null}}}
\def\edproofmarker{\strut\vadjust{\kern-2\strutdepth\epmarker}}
\def\endproof{\edproofmarker\vskip10pt}
\def\table#1\endtable{\begin{tabular}#1 \end{tabular}}
\def\smallbullet{\raise.5pt\hbox{\scriptsize$\,\bullet\,$}}

\def\tait{\cite{PGT}}

\def\petersen{\cite{JP}}
\def\tutte{\cite{WTT}}
\def\conway{\cite{JHC}}

\begin{document}
\newtheorem{theorem}{Theorem}[section]
\newtheorem{corollary}{Corollary}[section]
\newtheorem{lemma}{Lemma}[section]
\newtheorem{proposition}{Proposition}[section]
\newtheorem{example}{Example}[section]
{\newtheorem{definition}{Definition}[section]}

\title{O-cycles, vertex-oriented graphs, and the four colour theorem}
\author{Ortho Flint and Stuart Rankin}
\maketitle

\begin{abstract}
In 1880, P. G. Tait showed that the four colour theorem is equivalent to
the assertion that every 3-regular planar graph without cut-edges is
3-edge-colourable, and in 1891, J. Petersen proved that every 3-regular graph
with at most two cut-edges has a 1-factor. In this paper, we introduce the notion
of collapsing all edges of a 1-factor of a 3-regular planar graph, thereby obtaining
what we call a vertex-oriented 4-regular planar graph. We also introduce the notion
of o-colouring a vertex-oriented 4-regular planar graph, and we prove that the four
colour theorem is equivalent to the assertion that every vertex-oriented 4-regular
planar graph without nontransversally oriented cut-vertex (VOGWOC in short) is 3-o-colourable. This work 
proposes an alternative avenue of investigation in the search to find a more
conceptual proof of the four colour theorem, and we are able to prove that
every VOGWOC is o-colourable (although we have not yet been able to prove
3-o-colourability). 
\end{abstract}

\section{Introduction}
In 1880, P. G. Tait \tait\  showed that the four colour theorem is equivalent to
the assertion that every 3-regular planar graph without cut-edges is
3-edge-colourable (and that the latter is true if and only if every 3-regular
3-edge-connected planar graph is 3-edge-colourable). As is well known, Tait
actually felt that he had proven the four colour theorem since he had assumed that every 3-regular
3-edge-connected planar graph was hamiltonian (it being easily seen that hamiltonian
3-regular graphs are 3-edge-colourable), and it was not until 1946 that 
W. Tutte showed in \tutte\ that this is not the case. 

In this paper, we introduce the notion of a vertex-oriented 4-regular planar graph, and 
use it to transform Tait's theorem into another equivalent formulation of the 
four colour theorem. This came about as a result of our wish to provide a more conceptual proof 
of the four colour theorem, and was motivated by our work with 4-regular graphs in the study
of knot theory. In the third section of this paper, we establish that every vertex-oriented 4-regular
planar graph without nontransversally oriented cut-vertex (VOGWOC) is o-colourable (although we are not able to prove
3-o-colourability). It does follow from this result that every vertex-oriented 
4-regular planar graph is an edge-disjoint union of o-cycles (this is of course obvious from
the four colour theorem, but we were unable to prove directly that a given VOGWOC even had
a single o-cycle). We conclude that section with some remarks on how the proof of the
o-colourability result might be improved upon to give 3-o-colourability. We conclude the paper
with a study of the vertex-orientations of a regular projection of the Borromean rings (that is,
the basic polyhedral graph $6^*$).

\section{An equivalent formulation of the four colour theorem}
A vertex $v$ with no incident loop in a 4-regular planar graph $G$ shall
be said to be oriented if  the four edges incident to $v$ have been
partitioned into two cells (called the edge cells at $v$) of two edges
each so that the two edges in each cell are consecutive in the embedding
order at $v$.  If there is exactly one loop $e$ incident to $v$, then if
we denote the other two incident edges by $f$ and $g$, the set of
two subsets $\set e,f\endset$, $\set e,g\endset$ is said to be the
transverse orientation of $v$ (and we shall refer to the sets $\set e,f\endset$
and $\set e,g\endset$ as the edge cells at $v$, even though they are not disjoint),
while the set of two subsets $\set
f,g\endset$, $\set e\endset$ is the nontransverse orientation of $v$.
Finally, if there are two loops $e_1$, $e_2$ incident to $v$, then 
we only define one orientation at $v$; namely $\set\set e_1, e_2\endset\endset$,
and shall refer to this as the transverse orientation of $v$ (if $G$ is connected with two or more vertices,
this situation will never arise).
A vertex that has an incident loop shall be called a loop-anchor in $G$.

For example, if $v$ has incident edges $e$, $f$, $g$, $h$, (or a loop $e$
and incident edges $f$ and $g$), labelled in a
clockwise order, then one orientation of $v$ would be the partition
$\set \set e,f\endset,\set g, h\endset\endset$, while the other
orientation would be $\set \set f,g\endset,\set e, h\endset\endset$ (in
the case of the loop, the transverse orientation of $v$ would be
$\set \set e,f\endset, \set e,g\endset\endset$, while the nontraverse
orientation of $v$ would be $\set \set f,g\endset, \set e\endset\endset$). 
In a plane embedding of $G$, we shall indicate these by a double headed
arrow passing through $v$ in such a way that for each cell, the arrow
separates the two edges in the cell.

\medskip
\begin{figure}[h]
\begin{tabular}{c@{\hskip10pt}c@{\hskip10pt}c@{\hskip10pt}c}
\centering
\hbox{\xy /r20pt/:,
(0,0)="1";(1,1)="3"**\dir{-},
(1,0)="2";(0,1)="4"**\dir{-},
(.5,.5)="0"*{\smallbullet},
"0"*!<0pt,9pt>{v},
"3"*!<-4pt,-4pt>{h},
"2"*!<-4pt,4pt>{g},
"1"*!<4pt,4pt>{f},
"4"*!<4pt,-4pt>{e},
(0,.5)="x";(1,.5)="y"**\dir{-}*\dir{>},"x"*\dir{<},
\endxy}
&
\hbox{\xy /r20pt/:,
(0,0)="1";(1,1)="3"**\dir{-},
(1,0)="2";(0,1)="4"**\dir{-},
(.5,.5)="0"*{\smallbullet},
"0"*!<-9pt,0pt>{v},
"3"*!<-4pt,-4pt>{h},
"2"*!<-4pt,4pt>{g},
"1"*!<4pt,4pt>{f},
"4"*!<4pt,-4pt>{e},
(.5,0)="x";(.5,1)="y"**\dir{-}*\dir{>},"x"*\dir{<},
\endxy}
&
\hbox{\xy /r20pt/:,
(0,0)="1";(1,1)="3",
(1,0)="2";(0,1)="4",
(.5,.5)="o"*{\smallbullet},
"1";"o"**\dir{-},
"4";"o"**\dir{-},
"o";"o"**\crv{"3"+(.25,.25) & (1.75,.5) & "2"+(.25,-.25)},
"o"*!<9pt,0pt>{v},
(1.5,.5)*!<-4pt,0pt>{e},
"4"*!<4pt,-4pt>{f},
"1"*!<4pt,4pt>{g},
(.5,0)="x";(.5,1)="y"**\dir{-}*\dir{>},"x"*\dir{<},
\endxy}
&
\hbox{\xy /r20pt/:,
(0,0)="1";(1,1)="3",
(1,0)="2";(0,1)="4",
(.5,.5)="o"*{\smallbullet},
"1";"o"**\dir{-},
"4";"o"**\dir{-},
"o";"o"**\crv{"3"+(.25,.25) & (1.75,.5) & "2"+(.25,-.25)},
"o"*!<0pt,9pt>{v},
(1.5,.5)*!<-4pt,0pt>{e},
"4"*!<4pt,-4pt>{f},
"1"*!<4pt,4pt>{g},
(0,.5)="x";(1,.5)="y"**\dir{-}*\dir{>},"x"*\dir{<},
\endxy}
\\
\noalign{\vskip4pt}
$\set \set f,g\endset,\set e\endset\endset$ & $\set \set e,f\endset,\set g, h\endset\endset$ & $\set \set f,g\endset,\set e, h\endset\endset$ &
$\set \set e,f\endset,\set e,g\endset\endset$\\
\end{tabular}
\caption{}\label{figure: vertex orientation}
\end{figure}

\medskip
A mapping $\sigma$ such that for each vertex $v$ of $G$, $\sigma(v)$ is
an orientation of $v$, shall be called a vertex-orientation of $G$, and
we say that $G$ has been vertex-oriented by $\sigma$, or that
$(G,\sigma)$ is a vertex-oriented graph. Suppose that $(G,\sigma)$ is a
vertex-oriented 4-regular planar graph. We say that an edge colour
assignment $\varepsilon$ is an o-colouring of $(G,\sigma)$ if for each
$v\in G$, exactly two colours appear on the four edges incident to $v$,
and in each cell of $\sigma(v)$, both colours appear. The colour assignment
$\varepsilon$ is then called an o-colouring of the vertex-oriented graph
$(G,\sigma)$. If at most $k$ colours have been used, then we say that
$(G,\sigma)$ has been $k$-o-coloured. The least $k$ such that there is a
$k$-o-colouring of $(G,\sigma)$ shall be called the o-chromatic index of
$(G,\sigma)$ and denoted by $\chi_o(G)$. Note that
$\chi_o(G)\ge2$ for every vertex oriented 4-regular planar graph
$(G,\sigma)$.

If a vertex-oriented 4-regular planar graph has been o-coloured,
then the set of all edges of a given colour form one or more (vertex
and edge) disjoint cycles in the graph. In particular, if a 
vertex-oriented 4-regular planar graph has an
o-colouring, then every cut-vertex of the graph must be oriented transversely.

\begin{theorem}\label{theorem: 4ct equivalence}
 The 4-colour theorem is equivalent to the assertion that every 
 vertex-orientation of any 4-regular planar graph with no cut-vertex
 can be 3-o-coloured.
\end{theorem}

\proof
 By Tait's result, it suffices to prove that the assertion that every
 3-regular planar graph with no cut-edge can be 3-edge-coloured is
 equivalent to the assertion that every vertex-orientation of any
 4-regular planar graph with no cut-vertex can be
 3-o-coloured.
 
 Suppose that every 3-regular planar graph with no cut-edge can be
 3-edge-coloured, and let $G$ be a 4-regular planar graph with
 no cut-vertex. Further suppose that $\alpha$ is a vertex-orientation of $G$. At
 each vertex $v$ with orientation $\set \set e,f\endset,\set g,
 h\endset\endset$, replace $v$ by a new edge with endpoints $x$ and $y$,
 with $e,f$ incident to $x$ and $g,h$ incident to $y$. The result is a
 3-regular planar graph $H$. If $H$ has a cut-edge, $t$ say, then either
 $t$ is an edge in $G$, in which case each endpoint of $t$ is a
 cut-vertex of $G$, or else $t$ was one of the newly created edges,
 replacing vertex $v$ say, in which case $v$ is a cut-vertex of $G$.
 Since $G$ was without cut-vertices, neither of these situations is
 possible. Thus $H$ has no cut-edge, and so by hypothesis, $H$ can be
 3-edge-coloured. Suppose that $H$ has been 3-edge-coloured.
 Contract all edges of $H$ that were not edges of $G$, thereby
 obtaining $G$, but now each edge of $G$ has been assigned one of three
 colours. Moreover, at each vertex $v$ with orientation $\set \set
 e,f\endset,\set g, h\endset\endset$, the two colours that appear on $e$
 and $f$ are the same as the two colours that appear on $g$ and $h$.
 The result is therefore a 3 o-colouring of $(G,\alpha)$.
 
 Conversely, suppose that every vertex-orientation of any 4-regular
 planar graph with no cut-vertex can  be 3-o-coloured. We prove that
 every 3-regular planar graph with no cut-edge can be 3-edge-coloured by
 induction on the number of vertices. To begin with, we observe that a
 3-regular graph without cut-edge is also without loops. Thus the base
 case consists of the 3-regular planar graphs without cut-edge on two vertices,
 of which there is only one and it can be 3-edge-coloured. Suppose now
 that $n>2$ is an integer such that any 3-regular planar graph without cut-edge
 and fewer than $n$ vertices can be 3-o-coloured, and let $G$ be a 3-regular planar
 graph without cut-edge on $n$ vertices. As observed above, $G$ can't have
 any loops. By our inductive hypothesis, we may assume that $G$ is connected. 
 Furthermore, suppose that $G$ contains a digon. Then we may replace the digon (two vertices
 and the four edges incident to one or the other of the two vertices) by
 a single edge, resulting in a 3-regular planar graph without cut-edge on $n-2$
 vertices, which by our induction hypothesis is 3-edge-colourable. But then $G$ is
 3-edge-colourable. Thus we may further assume that $G$ is without digons.
 Petersen established in \petersen\ that every 3-regular graph with at most two cut-edges has a
 1-factor, so let $F$ be a 1-factor of $G$. Contract each edge $f\in F$, 
 putting the two edges incident to an endpoint of $f$ into a cell. The
 result is an orientation of the vertex formed by contracting $f$, and
 so we have formed a 4-regular planar graph $G'$ and given it a
 vertex-orientation.  
 
 Suppose that $G'$ has a cut-vertex $v$, say.  By
 the handshake lemma, $G'-v$ must consist of two components, and for
 each component,  there are exactly two edges incident to $v$ with
 endpoints in the component. Furthermore, since $G'$ is planar, the
 two edges incident to $v$ with endpoints in the same component of 
 $G'-v$ must be consecutive in the embedding order at $v$. Let $f\in F$
 denote the edge of $G$ that was contracted to form $v$, and let $x$ and
 $y$ denote the endpoints of $f$. Furthermore, let $e_1$ and $e_2$,
 respectively $f_1$ and $f_2$, denote the edges different from $f$ that
 are incident to $x$, respectively $y$. As well, let $x_{1}$ and $x_2$
 denote the non-$x$ endpoints of $e_1$ and $e_2$, respectively, and let
 $y_{1}$ and $y_2$ denote the non-$y$ endpoints of $f_1$ and $f_2$,
 respectively. As $G$ is without digons, it follows that $x_1\ne x_2$ and
 $y_1\ne y_2$. Since $f$ is not a cut-edge of $G$, there is a path in
 $G$ from $x$ to $y$ that does not use $f$, and so there is a path in
 $G'$ from either $x_1$ or $x_2$ to either $y_1$ or  $y_2$ that does not
 use any of $e_1$, $e_2$, $f_1$,  or $f_2$. We may suppose without loss
 of generality that the vertices have been labelled so that there is a
 path in $G'$ from $x_1$ to $y_1$ that does not use any of $e_1$, $e_2$, 
 $f_1$, or $f_2$, so that $x_1$ and $y_1$ belong to the same component of
 $G'-v$, and thus $x_2$ and $y_2$ belong to the other component of 
 $G'-v$. It follows that there exist simple closed curves $S_1$ and $S_2$ 
 (see Figure \ref{figure: cut-vertex}) such that of the edges of
 $G$, $S_1$ meets only $e_1$ and $f_1$ and contains vertices $x_1$ and
 $y_1$ in its interior, while $S_2$ meets only $e_2$ and $f_2$ and
 contains vertices $x_2$ and $y_2$ in its interior. Let $G_1$ and $G_2$
 denote the subgraphs of $G$ that are induced by the vertices of $G$ that lie
 in the interior of $S_1$ and $S_2$, respectively, with an additional edge to
 join $x_1$ to $y_1$ in $G_1$, and an additional edge to join $x_2$ to
 $y_2$ in $G_2$. Then $G_1$ and $G_2$ are 3-regular planar graphs with
 no cut-edge and fewer than $n$ vertices, so by the
 induction hypothesis, we may 3-edge-colour each of $G_1$ and $G_2$. By
 permuting the colours if necessary, we may arrange to have the new edge
 in $G_1$ coloured differently from the new edge in $G_2$, which then
 allows us to extend the colouring to obtain a 3-edge-colouring of $G$.

 \begin{figure}[h]
  \centering
   \table{c}
   $\vcenter{\xy /r17pt/:,
    (0,0)="centre","centre"+(1,0)="s";"centre"+(0,1)="e"**\crv{~*=<2.5pt>{\hbox{\Large.}} "s"+(0,.45) & "centre"+(.72,.72) & "e"+(.45,0)},
    "centre"+(-1,0)="s";"centre"+(0,1)="e"**\crv{~*=<2.5pt>{\hbox{\Large.}} "s"+(0,.45) & "centre"+(-.72,.72) & "e"-(.45,0)},
    "centre"+(1,0)="s";"centre"+(0,-1)="e"**\crv{~*=<2.5pt>{\hbox{\Large.}} "s"+(0,-.45) & "centre"+(.72,-.72) & "e"+(.45,0)},
    "centre"+(-1,0)="s";"centre"+(0,-1)="e"**\crv{~*=<2.5pt>{\hbox{\Large.}} "s"+(0,-.45) & "centre"+(-.72,-.72) & "e"-(.45,0)},
    (4,0)="centre","centre"+(1,0)="s";"centre"+(0,1)="e"**\crv{~*=<2.5pt>{\hbox{\Large.}} "s"+(0,.45) & "centre"+(.72,.72) & "e"+(.45,0)},
    "centre"+(-1,0)="s";"centre"+(0,1)="e"**\crv{~*=<2.5pt>{\hbox{\Large.}} "s"+(0,.45) & "centre"+(-.72,.72) & "e"-(.45,0)},
    "centre"+(1,0)="s";"centre"+(0,-1)="e"**\crv{~*=<2.5pt>{\hbox{\Large.}} "s"+(0,-.45) & "centre"+(.72,-.72) & "e"+(.45,0)},
    "centre"+(-1,0)="s";"centre"+(0,-1)="e"**\crv{~*=<2.5pt>{\hbox{\Large.}} "s"+(0,-.45) & "centre"+(-.72,-.72) & "e"-(.45,0)},
    (2,0)+(0,.25)="y"*{\smallbullet},   
    (2,0)+(0,-.25)="x"*{\smallbullet},   
    "x";"y"**\dir{-},
    (.6,-.5)="x1";"x"**\dir{-},
    "y";(3.4,.5)="y2"**\dir{-},
    (.6,.5)="y1";"y"**\dir{-},
    "x";(3.4,-.5)="x2"**\dir{-},   
    "x1"*!<7pt,2pt>{x_1},
    "y1"*!<7pt,-2pt>{y_1},
    "x2"*!<-8pt,2pt>{x_2},
    "y2"*!<-8pt,-2pt>{y_2},
    "x1"*{\smallbullet},
    "y1"*{\smallbullet},
    "x2"*{\smallbullet},
    "y2"*{\smallbullet},
    "x"*!<0pt,7pt>{x},
    "y"*!<0pt,-7pt>{y},
    (0,0)*!(1.1,-1){S_1},
    (4,0)*!(-1,-1){S_2},   
    (0,-1.8)*{},
    (0,1.7)*{},
    (9,0)="t",
    (0,0)+"t"="centre","centre"+(1,0)="s";"centre"+(0,1)="e"**\crv{~*=<2.5pt>{\hbox{\Large.}} "s"+(0,.45) & "centre"+(.72,.72) & "e"+(.45,0)},
    "centre"+(-1,0)="s";"centre"+(0,1)="e"**\crv{~*=<2.5pt>{\hbox{\Large.}} "s"+(0,.45) & "centre"+(-.72,.72) & "e"-(.45,0)},
    "centre"+(1,0)="s";"centre"+(0,-1)="e"**\crv{~*=<2.5pt>{\hbox{\Large.}} "s"+(0,-.45) & "centre"+(.72,-.72) & "e"+(.45,0)},
    "centre"+(-1,0)="s";"centre"+(0,-1)="e"**\crv{~*=<2.5pt>{\hbox{\Large.}} "s"+(0,-.45) & "centre"+(-.72,-.72) & "e"-(.45,0)},
    (4,0)+"t"="centre","centre"+(1,0)="s";"centre"+(0,1)="e"**\crv{~*=<2.5pt>{\hbox{\Large.}} "s"+(0,.45) & "centre"+(.72,.72) & "e"+(.45,0)},
    "centre"+(-1,0)="s";"centre"+(0,1)="e"**\crv{~*=<2.5pt>{\hbox{\Large.}} "s"+(0,.45) & "centre"+(-.72,.72) & "e"-(.45,0)},
    "centre"+(1,0)="s";"centre"+(0,-1)="e"**\crv{~*=<2.5pt>{\hbox{\Large.}} "s"+(0,-.45) & "centre"+(.72,-.72) & "e"+(.45,0)},
    "centre"+(-1,0)="s";"centre"+(0,-1)="e"**\crv{~*=<2.5pt>{\hbox{\Large.}} "s"+(0,-.45) & "centre"+(-.72,-.72) & "e"-(.45,0)},
    (2,0)+(0,.25)+"t"="y",
    (2,0)+(0,-.25)+"t"="x",
    (.6,-.5)+"t"="x1",
    (3.4,.5)+"t"="y2",
    (.6,.5)+"t"="y1",
    (3.4,-.5)+"t"="x2",
    "x1"*!<7pt,2pt>{x_1},
    "y1"*!<7pt,-2pt>{y_1},
    "x2"*!<-8pt,2pt>{x_2},
    "y2"*!<-8pt,-2pt>{y_2},
    "x1"*{\smallbullet},
    "y1"*{\smallbullet},
    "x2"*{\smallbullet},
    "y2"*{\smallbullet},
    "x1";"y1"**\dir{-},
    "x2";"y2"**\dir{-},   
    (0,0)+"t"*!(1.1,-1){G_1},
    (4,0)+"t"*!(-1,-1){G_2},   
    (2,-1.75)*{\hbox{(a)}},
    (2,-1.75)+"t"*{\hbox{(b)}},
  \endxy}$
 \endtable\caption{}\label{figure: cut-vertex}
\end{figure}
  
 \noindent 
  We may therefore assume that $G'$ has no cut-vertex; that is, $G'$ is a 
  vertex-oriented 4-regular planar graph without cut-vertex, with vertex-orientation
  $\alpha$ say, and by assumption, every such graph may be 3-o-coloured.
  Suppose then that $(G',\alpha)$ has been 3-o-coloured. Give
  each edge of $G$ the colour it has in $G'$, so that the only edges of $G$ that
  have not been coloured are those of $F$. Let $f\in F$, and let $x$ and
  $y$ denote the endpoints of $f$. Then the two edges incident to $x$ in
  $G'$ will be coloured with two different colours, say $c_1$ and $c_2$,
  and the two edges incident to $y$ will be coloured with the same two
  colours, one with $c_1$ and the other with $c_2$. Thus $f$ can be
  coloured with the third colour. The result is a 3-edge-colouring of
  $G$. This completes the proof of the inductive step, and so the result
  follows by induction.
 \endproof 

\section{O-colourings and o-cycles}
A walk $v_0,e_1,v_1,\ldots,e_n,v_n$ of length $n\ge 2$ in a vertex-oriented 
4-regular planar graph $(G,\sigma)$ shall be
called an o-walk if for each $i=1,2,\ldots,n-1$, $e_i$ and $e_{i+1}$ 
belong to different cells of $\sigma(v_i)$. An o-trail (respectively
o-circuit, o-cycle) is an o-walk that is a trail (respectively circuit,
cycle). If $(G,\sigma)$ has been o-coloured, then for each assigned
colour, the set of edges of $G$ that have been assigned that colour
forms a set of o-cycles with the property that no two have a vertex in
common. Thus an o-colouring of $(G,\sigma)$ provides a decomposition of
the edge set of $G$ into o-cycles, each of which has only edges of one
colour and such that any two cycles of the same colour have no vertex in
common.

Let $G_1$ and $G_2$ be (disjoint) graphs. Choose edges $e$ in $G_1$ and $f$ in
$G_2$ and remove them. Then join one end point of $e$ to one endpoint of
$f$, and join the other endpoint of $e$ to the other endpoint of $f$.
Denote the result by $G_1\#_{e,f}G_2$, or simply $G_1\# G_2$ when the
edges $e$ and $f$ are understood (there are two ways to carry out this construction,
but for convenience, we shall refer to both graphs -- in general, nonisomorphic --
by the same notation). Note that if $G_1$ and $G_2$ are
4-regular graphs, then $G_1\# G_2$ is also 4-regular, and if both $G_1$
and $G_2$ are planar, then $G_1\# G_2$ is planar. Conversely, suppose
that $G$ is a 4-regular graph. By the handshake lemma, it is  not
possible for $G$ to have a cut-edge. However, $G$ might have a cut-set
of size 2.  Suppose that $\set e,f\endset$ is in fact a cut-set for $G$.
Then again by the handshake lemma, $G-\set e,f\endset$ must have exactly
two connected components. Let $G_1$ denote the graph obtained from one
of these two components by creating an edge joining the endpoints of $e$
and $f$ that belong to the component (so the new edge is a loop if these
two endpoints are equal). Let $G_2$ denote the graph obtained by
applying the same procedure to the second component. Then $G=G_1\# G_2$ (that is
to say, one of the two ways to carry out the construction yields $G$).
Moreover, if $G$ is planar, then so are $G_1$ and $G_2$. Finally,
observe that there is a natural way to obtain vertex-orientations
$\sigma_1$ of $G_1$ and  $\sigma_2$ of $G_2$ from a vertex-orientation
$\sigma$ of $G$, and vice-versa, and we shall say that $\sigma$ is
compatible with $\sigma_1$ and $\sigma_2$ and vice-versa.

\begin{lemma}\label{lemma: prime is ok}
 Let $G$, $G_1$, and $G_2$ be 4-regular planar graphs such that $G=G_1\#
 G_2$. Suppose further that $G$ is vertex-oriented by $\sigma$, and give
 $G_1$ and $G_2$ the induced vertex-orientations $\sigma_1$ and
 $\sigma_2$, respectively. For every positive integer $k$, if 
 $(G_1,\sigma_1)$ and $(G_2,\sigma_2)$ can
 be $k$-o-coloured, then $(G,\sigma)$ can be $k$-o-coloured.
\end{lemma}

\proof
Suppose that $G_1$ and $G_2$ have been $k$-o-coloured. If the new edges
in $G_1$ and $G_2$ have been coloured differently, then we may permute
the colours in the colouring of $G_2$ to arrange that the two new edges
have been coloured the same, say with colour $c_1$. Then assign $e$ and
$f$ colour $c_1$ to obtain a $k$-o-colouring of $(G,\sigma)$.
\endproof 

Now suppose that $G$ is a 4-regular graph with a cut-vertex $v$. As we
have seen in the proof of Theorem \ref{theorem: 4ct equivalence}, $G-v$ must
consist of two components, and for each component, there are exactly two
edges incident to $v$ with endpoints in the component. Furthermore, since
$G$ is planar, the two edges incident to $v$ with endpoints in the
same component of  $G-v$ must be consecutive in the embedding order at
$v$. Thus in any plane embedding of $G$, there exists a simple closed
curve $S_1$ that meets exactly two edges incident to $v$ and no other
edges of $G$ and contains one of the components of  $G-v$ in its
interior, and a simple closed curve $S_2$ that meets the other two edges
incident to $v$ and no other edges of $G$ and contains the other
component of $G-v$ in its interior. Let $G_1$ be the graph formed from
one of the components of $G-v$ by creating a new edge whose endpoints
are those of the two edges incident to $v$ that meet the component in
question, and let $G_2$ be the graph constructed by the same process but
applied to the other component of $G-v$. We shall use the notation
$G=G_1\#_v G_2$ to denote this situation. Moreover, there is a natural
way to associate two different  vertex-orientations of $G$ corresponding
to a vertex-orientation of each of $G_1$ and $G_2$, depending on the
orientation assigned to $v$. We shall let
$\connsumparallelvertex{G_1}{G_2}{v}$ indicate the choice of
orientation at $v$ whose cells are the pairs of edges incident to $G_1$,
respectively $G_2$, and we shall call this the nontransverse orientation at
$v$. The other orientation, called the transverse orientation at $v$, shall be denoted by
$\connsumtransversevertex{G_1}{G_2}{v}$. 

If $G$ is a vertex-oriented 4-regular planar graph, then for any vertex
$v$ that is not a loop-anchor, form a new 4-regular planar graph by
removing $v$ and identifying each edge $e$ in an edge cell at $v$ with
the unique edge in the other edge cell at $v$ that is adjacent to $e$ in
the embedding order at $v$ (see Figure \ref{figure: smoothing figure}). If $v$ is
a loop-anchor, oriented transversely or non-transversely, smoothing $v$
is achieved by removing the loop and $v$ and identifying the other two
edges incident to $v$. The resulting graph $G'$  is  vertex-oriented,
and shall be said to have been obtained from $G$ by smoothing $v$.

\begin{figure}[h]
\centering
 \table{c@{\hskip60pt}c}\\
 $\vcenter{\xy /r25pt/:,
  (0,0);(1,1)**\dir{-},
  (0,1);(1,0)**\dir{-},  
  (.5,.5)*{\smallbullet}*!<-8pt,0pt>{v},
  (.5,0);(.5,1)**\dir{-}*\dir{>},(.5,0)*\dir{<},
 \endxy}$ 
 &
 $\vcenter{\xy /r25pt/:,
  (0,0);(0,1)**\crv{(.4,.5)},
  (1,0);(1,1)**\crv{(.6,.5)},
 \endxy}$ 
 \\
 \noalign{\vskip6pt}
 (a) & (b)
 \endtable\caption{}\label{figure: smoothing figure}
\end{figure}
 
\begin{lemma}\label{lemma: transverse cut vertex can be o-coloured}
 Let $(G,\sigma)$ be a vertex-oriented 4-regular planar graph with a
 cut-vertex $v$ transversely oriented, so that 
 $G=\connsumtransversevertex{G_1}{G_2}{v}$ for some vertex-oriented
 4-regular planar graphs $(G_1,\sigma_1)$ and $(G_2,\sigma_2)$ such that
 $\sigma_1$ and $\sigma_2$ are consistent with $\sigma$. If
 $(G_1,\sigma_1)$ and $(G_2,\sigma_2)$ can be $k$-o-coloured, then
 $(G,\sigma)$ can be $k$-o-coloured.
\end{lemma}

\proof
 Embed $G$ in the plane as shown in Figure \ref{figure: transverse cut figure}
 (a), where each of the closed curves $S_1$ and $S_2$ contain at least
 one vertex in their respective interiors, and then smooth $v$,
 obtaining 4-regular planar graphs $G_1$ and $G_2$ as shown in
 Figure \ref{figure: transverse cut figure} (b). By assumption, we may
 o-colour each of $G_1$ and $G_2$ with $k\ge2$ colours. Suppose that colour $c_1$ appears on $e$, and choose
 a second colour $c_2$. By permuting the colours in $G_2$ if
 necessary, we can arrange to have $f$ coloured with $c_2$. Then colour
 every edge of $G$ that is an edge in either $G_1$ or $G_2$ with the
 colour it has in the respective graphs, and colour the edges incident
 to $v$ as shown in Figure \ref{figure: transverse coloured figure}. 
 The result is a $k$-o-colouring for $(G,\sigma)$.

 \begin{figure}[h]
  \centering\table{c}
  $\vcenter{\xy /r20pt/:,
   (0,0)="centre","centre"+(1,0)="s";"centre"+(0,1)="e"**\crv{~*=<2.5pt>{\hbox{\Large.}} "s"+(0,.45) & "centre"+(.72,.72) & "e"+(.45,0)},
   "centre"+(-1,0)="s";"centre"+(0,1)="e"**\crv{~*=<2.5pt>{\hbox{\Large.}} "s"+(0,.45) & "centre"+(-.72,.72) & "e"-(.45,0)},
   "centre"+(1,0)="s";"centre"+(0,-1)="e"**\crv{~*=<2.5pt>{\hbox{\Large.}} "s"+(0,-.45) & "centre"+(.72,-.72) & "e"+(.45,0)},
   "centre"+(-1,0)="s";"centre"+(0,-1)="e"**\crv{~*=<2.5pt>{\hbox{\Large.}} "s"+(0,-.45) & "centre"+(-.72,-.72) & "e"-(.45,0)},
   (4,0)="centre","centre"+(1,0)="s";"centre"+(0,1)="e"**\crv{~*=<2.5pt>{\hbox{\Large.}} "s"+(0,.45) & "centre"+(.72,.72) & "e"+(.45,0)},
   "centre"+(-1,0)="s";"centre"+(0,1)="e"**\crv{~*=<2.5pt>{\hbox{\Large.}} "s"+(0,.45) & "centre"+(-.72,.72) & "e"-(.45,0)},
   "centre"+(1,0)="s";"centre"+(0,-1)="e"**\crv{~*=<2.5pt>{\hbox{\Large.}} "s"+(0,-.45) & "centre"+(.72,-.72) & "e"+(.45,0)},
   "centre"+(-1,0)="s";"centre"+(0,-1)="e"**\crv{~*=<2.5pt>{\hbox{\Large.}} "s"+(0,-.45) & "centre"+(-.72,-.72) & "e"-(.45,0)},
   (.6,-.5);(3.4,.5)**\crv{(1.5,-.5) & (2,0) & (2.5,.5)},
   (.6,.5);(3.4,-.5)**\crv{(1.5,.5) & (2,0) & (2.5,-.5)},   
   (2,0)*{\smallbullet}*!<-8pt,0pt>{v},   
   (2,-.45);(2,.45)**\dir{-}*\dir{>},(2,-.45)*\dir{<},
   (0,0)*!(1,-1){S_1},
   (4,0)*!(-1,-1){S_2},   
   (0,-1.8)*{},
   (0,1.7)*{},
   (8,0)="x",
   (0,0)+"x"="centre","centre"+(1,0)="s";"centre"+(0,1)="e"**\crv{~*=<2.5pt>{\hbox{\Large.}} "s"+(0,.45) & "centre"+(.72,.72) & "e"+(.45,0)},
   "centre"+(-1,0)="s";"centre"+(0,1)="e"**\crv{~*=<2.5pt>{\hbox{\Large.}} "s"+(0,.45) & "centre"+(-.72,.72) & "e"-(.45,0)},
   "centre"+(1,0)="s";"centre"+(0,-1)="e"**\crv{~*=<2.5pt>{\hbox{\Large.}} "s"+(0,-.45) & "centre"+(.72,-.72) & "e"+(.45,0)},
   "centre"+(-1,0)="s";"centre"+(0,-1)="e"**\crv{~*=<2.5pt>{\hbox{\Large.}} "s"+(0,-.45) & "centre"+(-.72,-.72) & "e"-(.45,0)},
   (4,0)+"x"="centre","centre"+(1,0)="s";"centre"+(0,1)="e"**\crv{~*=<2.5pt>{\hbox{\Large.}} "s"+(0,.45) & "centre"+(.72,.72) & "e"+(.45,0)},
   "centre"+(-1,0)="s";"centre"+(0,1)="e"**\crv{~*=<2.5pt>{\hbox{\Large.}} "s"+(0,.45) & "centre"+(-.72,.72) & "e"-(.45,0)},
   "centre"+(1,0)="s";"centre"+(0,-1)="e"**\crv{~*=<2.5pt>{\hbox{\Large.}} "s"+(0,-.45) & "centre"+(.72,-.72) & "e"+(.45,0)},
   "centre"+(-1,0)="s";"centre"+(0,-1)="e"**\crv{~*=<2.5pt>{\hbox{\Large.}} "s"+(0,-.45) & "centre"+(-.72,-.72) & "e"-(.45,0)},
   (.6,-.5)+"x";(.6,.5)+"x"**\crv{(1.3,-.5)+"x" & (1.5,0)+"x" & (1.3,.5)+"x"},   
   (1.5,0)+"x"*!<-3pt,0pt>{e},
   (2.5,0)+"x"*!<3pt,0pt>{f},   
   (3.4,-.5)+"x";(3.4,.5)+"x"**\crv{(2.8,-.5)+"x" & (2.5,0)+"x" & (2.8,.5)+"x"},   
   (0,0)+"x"*!(0,1.6){G_1},
   (4,0)+"x"*!(0,1.6){G_2},   
   (2,-2.5)*{\hbox{(a)}},
   (2,-2.5)+"x"*{\hbox{(b)}},
  \endxy}$
 \endtable\caption{}\label{figure: transverse cut figure}
\end{figure}

\begin{figure}[h]
 \centering\table{c}  
   $\vcenter{\xy /r25pt/:,
   (0,0)="centre","centre"+(1,0)="s";"centre"+(0,1)="e"**\crv{~*=<2.5pt>{\hbox{\Large.}} "s"+(0,.45) & "centre"+(.72,.72) & "e"+(.45,0)},
   "centre"+(-1,0)="s";"centre"+(0,1)="e"**\crv{~*=<2.5pt>{\hbox{\Large.}} "s"+(0,.45) & "centre"+(-.72,.72) & "e"-(.45,0)},
   "centre"+(1,0)="s";"centre"+(0,-1)="e"**\crv{~*=<2.5pt>{\hbox{\Large.}} "s"+(0,-.45) & "centre"+(.72,-.72) & "e"+(.45,0)},
   "centre"+(-1,0)="s";"centre"+(0,-1)="e"**\crv{~*=<2.5pt>{\hbox{\Large.}} "s"+(0,-.45) & "centre"+(-.72,-.72) & "e"-(.45,0)},
   (4,0)="centre","centre"+(1,0)="s";"centre"+(0,1)="e"**\crv{~*=<2.5pt>{\hbox{\Large.}} "s"+(0,.45) & "centre"+(.72,.72) & "e"+(.45,0)},
   "centre"+(-1,0)="s";"centre"+(0,1)="e"**\crv{~*=<2.5pt>{\hbox{\Large.}} "s"+(0,.45) & "centre"+(-.72,.72) & "e"-(.45,0)},
   "centre"+(1,0)="s";"centre"+(0,-1)="e"**\crv{~*=<2.5pt>{\hbox{\Large.}} "s"+(0,-.45) & "centre"+(.72,-.72) & "e"+(.45,0)},
   "centre"+(-1,0)="s";"centre"+(0,-1)="e"**\crv{~*=<2.5pt>{\hbox{\Large.}} "s"+(0,-.45) & "centre"+(-.72,-.72) & "e"-(.45,0)},
   (.6,-.5)="a";(3.4,.5)="b"**\crv{(1.5,-.5) & (2,0) & (2.5,.5)},
   (.6,.5)="c";(3.4,-.5)="d"**\crv{(1.5,.5) & (2,0) & (2.5,-.5)},   
   (2,0)*{\smallbullet}*!<-8pt,0pt>{v},   
   (2,-.45);(2,.45)**\dir{-}*\dir{>},(2,-.45)*\dir{<},
   "a"*!<-15pt,5pt>{c_1},
   "c"*!<-15pt,-4pt>{c_1},   
   "b"*!<15pt,-4pt>{c_2},
   "d"*!<15pt,5pt>{c_2},   
   \endxy}$
 \endtable\caption{}\label{figure: transverse coloured figure}
\end{figure}
 
\endproof

\begin{lemma}\label{lemma: flyping parallel}
 If $(G,\sigma)$ is a vertex-oriented 4-regular planar graph of the form
 as shown in Figure \ref{figure: flyping parallel figure} (a), (where it is not
 intended that the endpoints of the edges entering $S_1$, repectively
 $S_2$, need be distinct), and each of $S_1$ and $S_2$ contain at least
 one vertex in their interior, and each of the compatibly
 vertex-oriented 4-regular planar graphs $(G_1,\sigma_1)$ and
 $(G_2,\sigma_2)$ in Figure \ref{figure: flyping parallel figure} (b) can be
 $k$-o-coloured, then $(G,\sigma)$ can be $k$-o-coloured.
\end{lemma}

 \begin{figure}[h]
  \centerline{\table{c}
   $\vcenter{\xy /r15pt/:,
   (0,2.2)*{},
   (-.5,0)="centre",{(0,0);(1.5,0):,"centre"+(1,0)="s";"centre"+(0,1)="e"**\crv{~*=<2.5pt>{\hbox{\Large.}} "s"+(0,.45) & "centre"+(.72,.72) & "e"+(.45,0)},
   "centre"+(-1,0)="s";"centre"+(0,1)="e"**\crv{~*=<2.5pt>{\hbox{\Large.}} "s"+(0,.45) & "centre"+(-.72,.72) & "e"-(.45,0)},
   "centre"+(1,0)="s";"centre"+(0,-1)="e"**\crv{~*=<2.5pt>{\hbox{\Large.}} "s"+(0,-.45) & "centre"+(.72,-.72) & "e"+(.45,0)},
   "centre"+(-1,0)="s";"centre"+(0,-1)="e"**\crv{~*=<2.5pt>{\hbox{\Large.}} "s"+(0,-.45) & "centre"+(-.72,-.72) & "e"-(.45,0)}},
   (4.5,0)="centre",{(0,0);(1.5,0):,"centre"+(1,0)="s";"centre"+(0,1)="e"**\crv{~*=<2.5pt>{\hbox{\Large.}} "s"+(0,.45) & "centre"+(.72,.72) & "e"+(.45,0)},
   "centre"+(-1,0)="s";"centre"+(0,1)="e"**\crv{~*=<2.5pt>{\hbox{\Large.}} "s"+(0,.45) & "centre"+(-.72,.72) & "e"-(.45,0)},
   "centre"+(1,0)="s";"centre"+(0,-1)="e"**\crv{~*=<2.5pt>{\hbox{\Large.}} "s"+(0,-.45) & "centre"+(.72,-.72) & "e"+(.45,0)},
   "centre"+(-1,0)="s";"centre"+(0,-1)="e"**\crv{~*=<2.5pt>{\hbox{\Large.}} "s"+(0,-.45) & "centre"+(-.72,-.72) & "e"-(.45,0)}},
   (.6,-.5)="a";(3.4,.5)="b"**\crv{(1.5,-.5) & (2,0) & (2.5,.5)},
   (.6,.5)="c";(3.4,-.5)="d"**\crv{(1.5,.5) & (2,0) & (2.5,-.5)},   
   (2,0)*{\smallbullet}*!<0pt,-8pt>{v},   
   (1.5,0);(2.5,0)**\dir{-}*\dir{>},(1.5,0)*\dir{<},
   "a"+(-.2,-.5);"d"+(.2,-.5)**\dir{-};
   "b"+(.2,.5);"c"+(-.2,.5)**\dir{-};   
   (-.5,0)*!(1.5,-1.5){S_1},
   (4.5,0)*!(-1.5,-1.5){S_2},   
   (2,-2.5)*{\hbox{(a)}},
   (10,0)="x",
   (1,0)="y",
   (-.5,0)+"x"="centre",{(0,0);(1.5,0):,"centre"+(1,0)="s";"centre"+(0,1)="e"**\crv{~*=<2.5pt>{\hbox{\Large.}} "s"+(0,.45) & "centre"+(.72,.72) & "e"+(.45,0)},
   "centre"+(-1,0)="s";"centre"+(0,1)="e"**\crv{~*=<2.5pt>{\hbox{\Large.}} "s"+(0,.45) & "centre"+(-.72,.72) & "e"-(.45,0)},
   "centre"+(1,0)="s";"centre"+(0,-1)="e"**\crv{~*=<2.5pt>{\hbox{\Large.}} "s"+(0,-.45) & "centre"+(.72,-.72) & "e"+(.45,0)},
   "centre"+(-1,0)="s";"centre"+(0,-1)="e"**\crv{~*=<2.5pt>{\hbox{\Large.}} "s"+(0,-.45) & "centre"+(-.72,-.72) & "e"-(.45,0)}},
   (4.5,0)+"x"+"y"="centre",{(0,0);(1.5,0):,"centre"+(1,0)="s";"centre"+(0,1)="e"**\crv{~*=<2.5pt>{\hbox{\Large.}} "s"+(0,.45) & "centre"+(.72,.72) & "e"+(.45,0)},
   "centre"+(-1,0)="s";"centre"+(0,1)="e"**\crv{~*=<2.5pt>{\hbox{\Large.}} "s"+(0,.45) & "centre"+(-.72,.72) & "e"-(.45,0)},
   "centre"+(1,0)="s";"centre"+(0,-1)="e"**\crv{~*=<2.5pt>{\hbox{\Large.}} "s"+(0,-.45) & "centre"+(.72,-.72) & "e"+(.45,0)},
   "centre"+(-1,0)="s";"centre"+(0,-1)="e"**\crv{~*=<2.5pt>{\hbox{\Large.}} "s"+(0,-.45) & "centre"+(-.72,-.72) & "e"-(.45,0)}},
   (.6,-.5)+"x"="a",
   "a"+(-.2,-.5)="bl", 
   (3.4,.5)+"x"+"y"="b",
   (.6,.5)+"x"="c",
   (3.4,-.5)+"x"+"y"="d",
   "a";"c"**\dir{}?(.5)="m",
   "m"+(1.1,0)="arrow"*{\smallbullet},
   "arrow"-(.45,0)="le";"arrow"+(.45,0)**\dir{-}*\dir{>},"le"*\dir{<},
   "d"+(.2,-.5)="br",
   "b"+(.2,.5)="tr",
   "c"+(-.2,.5)="tl",
   "tl";"a"**\crv{"tl"+(3,0) & "a"+(.8,0)},   
   "bl";"c"**\crv{"bl"+(3,0) & "c"+(.8,0)},      
   "tr";"d"**\crv{"tr"+(-3,0) & "d"+(-.8,0)},   
   "br";"b"**\crv{"br"+(-3,0) & "b"+(-.8,0)},      
   "b";"d"**\dir{}?(.5)="n",
   "n"+(-1.1,0)="arrow"*{\smallbullet},
   "arrow"-(.45,0)="le";"arrow"+(.45,0)**\dir{-}*\dir{>},"le"*\dir{<},
   (-.5,0)+"x"*!(1.5,-1.5){S_1},
   (4.5,0)+"x"+"y"*!(-1.5,-1.5){S_2},   
   (-.5,0)+"x"*!(0,2){G_1},
   (4.5,0)+"x"+"y"*!(0,2){G_2},   
   (2.5,-2.5)+"x"*{\hbox{(b)}},   
   \endxy}$
  \endtable}\caption{}\label{figure: flyping parallel figure}
 \end{figure}
 
 \proof
  By hypothesis, both $(G_1,\sigma_1)$ and $(G_2,\sigma_2)$ can be
  $k$-o-coloured. Label the (necessarily distinct) colours on the top
  and bottom edges incident to the copy of $v$ in $G_1$ as $c_1$ and
  $c_2$, and label the colours on the other two edges incident to that
  vertex with $x$ and $y$, so that $\set x,y\endset=\set
  c_1,c_2\endset$. By permutating the colours in $G_2$ if necessary, we
  can ensure that the $k$-o-colouring of $(G_2,\sigma_2)$ is as shown
  in Figure \ref{figure: new graphs}, where $\set r,s\endset=\set c_1,c_2\endset$. 

  \medskip
  \begin{figure}[h]
 \centering\table{c}
   $\vcenter{\xy /r20pt/:,
    (0,0)="x",
    (,0)="y",
    (-.5,0)+"x"="centre",{(0,0);(1.5,0):,"centre"+(1,0)="s";"centre"+(0,1)="e"**\crv{~*=<2.5pt>{\hbox{\Large.}} "s"+(0,.45) & "centre"+(.72,.72) & "e"+(.45,0)},
    "centre"+(-1,0)="s";"centre"+(0,1)="e"**\crv{~*=<2.5pt>{\hbox{\Large.}} "s"+(0,.45) & "centre"+(-.72,.72) & "e"-(.45,0)},
    "centre"+(1,0)="s";"centre"+(0,-1)="e"**\crv{~*=<2.5pt>{\hbox{\Large.}} "s"+(0,-.45) & "centre"+(.72,-.72) & "e"+(.45,0)},
    "centre"+(-1,0)="s";"centre"+(0,-1)="e"**\crv{~*=<2.5pt>{\hbox{\Large.}} "s"+(0,-.45) & "centre"+(-.72,-.72) & "e"-(.45,0)}},
    (4.5,0)+"x"+"y"="centre",{(0,0);(1.5,0):,"centre"+(1,0)="s";"centre"+(0,1)="e"**\crv{~*=<2.5pt>{\hbox{\Large.}} "s"+(0,.45) & "centre"+(.72,.72) & "e"+(.45,0)},
    "centre"+(-1,0)="s";"centre"+(0,1)="e"**\crv{~*=<2.5pt>{\hbox{\Large.}} "s"+(0,.45) & "centre"+(-.72,.72) & "e"-(.45,0)},
    "centre"+(1,0)="s";"centre"+(0,-1)="e"**\crv{~*=<2.5pt>{\hbox{\Large.}} "s"+(0,-.45) & "centre"+(.72,-.72) & "e"+(.45,0)},
    "centre"+(-1,0)="s";"centre"+(0,-1)="e"**\crv{~*=<2.5pt>{\hbox{\Large.}} "s"+(0,-.45) & "centre"+(-.72,-.72) & "e"-(.45,0)}},
    (.6,-.5)+"x"="a",
    "a"+(-.2,-.5)="bl", 
    (3.4,.5)+"x"+"y"="b",
    (.6,.5)+"x"="c",
    (3.4,-.5)+"x"+"y"="d",
    "a";"c"**\dir{}?(.5)="m",
    "m"+(1.1,0)="arrow"*{\smallbullet},
    "arrow"-(.45,0)="le";"arrow"+(.45,0)**\dir{-}*\dir{>},"le"*\dir{<},
    "d"+(.2,-.5)="br",
    "b"+(.2,.5)="tr",
    "c"+(-.2,.5)="tl",
    "tl";"a"**\crv{"tl"+(3,0) & "a"+(.8,0)},   
    "bl";"c"**\crv{"bl"+(3,0) & "c"+(.8,0)},      
    "tr";"d"**\crv{"tr"+(-3,0) & "d"+(-.8,0)},   
    "br";"b"**\crv{"br"+(-3,0) & "b"+(-.8,0)},      
    "b";"d"**\dir{}?(.5)="n",
    "n"+(-1.1,0)="arrow"*{\smallbullet},
    "arrow"-(.45,0)="le";"arrow"+(.45,0)**\dir{-}*\dir{>},"le"*\dir{<},
    (-.5,0)+"x"*!(1.5,-1.5){S_1},
    (4.5,0)+"x"+"y"*!(-1.5,-1.5){S_2},   
    (-.5,0)+"x"*!(0,2){G_1},
    (4.5,0)+"x"+"y"*!(0,2){G_2},   
    "tl"+(1,0)*!<-8pt,-2pt>{c_1},
    "tr"+(-1,0)*!<8pt,-2pt>{c_1},   
    "bl"+(1,0)*!<-8pt,2pt>{c_2},
    "br"+(-1,0)*!<8pt,2pt>{c_2},   
    "c"+(.25,0)*!<-9pt,0pt>{x},
    "b"+(-.25,0)*!<9pt,0pt>{r},   
    "a"+(.25,0)*!<-9pt,1pt>{y},
    "d"+(-.25,0)*!<9pt,1pt>{s},   
   \endxy}$
 \endtable\caption{}\label{figure: new graphs}
 \end{figure}
 
 \medskip
 \noindent Now assign to each edge of $G$ that is also an edge of either
 $G_1$ or $G_2$ the colour it has been assigned in the $k$-o-colouring
 of the respective graphs, and complete the colouring of the edges
 incident to $v$ as shown in Figure \ref{figure: another one}.

\medskip
 \begin{figure}[h]
 \centering\table{c}
   $\vcenter{\xy /r20pt/:,
   (0,2.2)*{},
    (-.5,0)="centre",{(0,0);(1.5,0):,"centre"+(1,0)="s";"centre"+(0,1)="e"**\crv{~*=<2.5pt>{\hbox{\Large.}} "s"+(0,.45) & "centre"+(.72,.72) & "e"+(.45,0)},
    "centre"+(-1,0)="s";"centre"+(0,1)="e"**\crv{~*=<2.5pt>{\hbox{\Large.}} "s"+(0,.45) & "centre"+(-.72,.72) & "e"-(.45,0)},
    "centre"+(1,0)="s";"centre"+(0,-1)="e"**\crv{~*=<2.5pt>{\hbox{\Large.}} "s"+(0,-.45) & "centre"+(.72,-.72) & "e"+(.45,0)},
    "centre"+(-1,0)="s";"centre"+(0,-1)="e"**\crv{~*=<2.5pt>{\hbox{\Large.}} "s"+(0,-.45) & "centre"+(-.72,-.72) & "e"-(.45,0)}},
    (4.5,0)="centre",{(0,0);(1.5,0):,"centre"+(1,0)="s";"centre"+(0,1)="e"**\crv{~*=<2.5pt>{\hbox{\Large.}} "s"+(0,.45) & "centre"+(.72,.72) & "e"+(.45,0)},
    "centre"+(-1,0)="s";"centre"+(0,1)="e"**\crv{~*=<2.5pt>{\hbox{\Large.}} "s"+(0,.45) & "centre"+(-.72,.72) & "e"-(.45,0)},
    "centre"+(1,0)="s";"centre"+(0,-1)="e"**\crv{~*=<2.5pt>{\hbox{\Large.}} "s"+(0,-.45) & "centre"+(.72,-.72) & "e"+(.45,0)},
    "centre"+(-1,0)="s";"centre"+(0,-1)="e"**\crv{~*=<2.5pt>{\hbox{\Large.}} "s"+(0,-.45) & "centre"+(-.72,-.72) & "e"-(.45,0)}},
   (.6,-.5)="a";(3.4,.5)="b"**\crv{(1.5,-.5) & (2,0) & (2.5,.5)},
   (.6,.5)="c";(3.4,-.5)="d"**\crv{(1.5,.5) & (2,0) & (2.5,-.5)},   
   (2,0)="v"*{\smallbullet}*!<0pt,-8pt>{v},   
   "v"+(0,1.3)*{c_1},
   "v"+(0,-1.3)*{c_2},   
   (1.5,0);(2.5,0)**\dir{-}*\dir{>},(1.5,0)*\dir{<},
   "a"+(-.2,-.5);"d"+(.2,-.5)**\dir{-};
   "b"+(.2,.5);"c"+(-.2,.5)**\dir{-};   
   "a"*!<-18pt,3pt>{y},
   "c"*!<-18pt,-2pt>{x},   
   "b"*!<18pt,-2pt>{r},
   "d"*!<18pt,3pt>{s},   
   (-.5,0)*!(1.5,-1.5){S_1},
   (4.5,0)*!(-1.5,-1.5){S_2},   
   \endxy}$
 \endtable
 \caption{}\label{figure: another one}
 \end{figure}
\medskip
 \noindent The result is a $k$-o-colouring of $G$.
 \endproof
 
 Before continuing on to the main theorem, we introduce one final bit of
 terminology. We say that an edge-colouring of a 4-regular planar graph
 $G$ is alternating at $v$ if exactly two colours appear on the edges
 incident to $v$, and they appear in alternating order as we examine the
 edges in the embedding order. If an edge-colouring of $G$ is
 alternating at $v$, then it is compatible with either of the two
 possible vertex-orientations at $v$. 
 
\begin{theorem}\label{theorem: o-colouring theorem}
 Every vertex-oriented 4-regular planar graph $(G,\sigma)$ in which each
 cut-vertex or loop-anchor is oriented transversely can be o-coloured.
\end{theorem}

\proof
 The proof is by induction on the number of vertices. There is only one
 such graph on a single vertex, and  two such graphs
 on two vertices. O-colourings for each are shown in Figure \ref{figure: base
 case figure}. Note that in Figure \ref{figure: base case figure} (b), we have
 given an edge-colouring that is alternating at each vertex, and is 
 therefore an o-colouring for any vertex-orientation of the graph.
 
 \begin{figure}[h]
 \centerline{\table{c@{\hskip40pt}c@{\hskip40pt}c}
  $\vcenter{\xy /r20pt/:,
   (0,0)="v"*{\smallbullet},  
   "v"-(0,.45)="le";"v"+(0,.45)**\dir{-}*\dir{>},"le"*\dir{<},   
   "v";"v"**\crv{ "v"+(.7,-.7) & "v"+(1.2,0) & "v"+(.7,.7)},
   "v";"v"**\crv{ "v"+(-.7,-.7) & "v"+(-1.2,0) & "v"+(-.7,.7)},   
   "v"+(-1,.6)*{c_1},
   "v"+(1,.6)*{c_2},   
  \endxy}$ 
  &
  $\vcenter{\xy /r20pt/:,
   (0,0)="v"*{\smallbullet},  
   (0,1)="w"*{\smallbullet},
   "v";"w"**\crv{ (-.6,.5)},
   "v";"w"**\crv{ (.6,.5)},   
   "v";"w"**\crv{ "v"+(-1,0) & (-2,.5) & "w"+(-1,0)},
   "v";"w"**\crv{ "v"+(1,0) & (2,.5) & "w"+(1,0)},   
   "v"+(-1.9,.7)*{c_1},
   "v"+(1.9,.7)*{c_2},   
   "v"+(-.7,.5)*{c_2},
   "v"+(.7,.5)*{c_1},   
  \endxy}$ 
  &
  $\vcenter{\xy /r20pt/:,
   (0,0)="v"*{\smallbullet},  
   "v"-(0,.45)="le";"v"+(0,.45)**\dir{-}*\dir{>},"le"*\dir{<},   
   "v";"v"**\crv{ "v"+(-.7,-.7) & "v"+(-1.2,0) & "v"+(-.7,.7)},   
   "v"+(-1,.6)*{c_1},
   (1.25,0)="w"*{\smallbullet},  
   "w"-(0,.45)="wle";"w"+(0,.45)**\dir{-}*\dir{>},"wle"*\dir{<},   
   "w";"w"**\crv{ "w"+(.7,-.7) & "w"+(1.2,0) & "w"+(.7,.7)},
   "w"+(1,.6)*{c_1},  
   "v";"w"**\dir{}?(.5)="x", 
   "v";"w"**\crv{"x"+(0,.75)},
   "v";"w"**\crv{"x"+(0,-.75)},   
   "x"+(0,.75)*{c_2},
   "x"-(0,.75)*{c_2},   
  \endxy}$ 
  \\
  \noalign{\vskip1pt}
  (a) & (b) & (c)
  \endtable}\caption{}\label{figure: base case figure}
 \end{figure}
 
 Suppose now that $n>2$ is an integer such that every vertex-oriented
 4-regular planar graph on fewer than $n$ vertices for which any
 cut-vertex or loop-anchor is oriented transversely can be o-coloured,
 and let $(G,\sigma)$ be a  vertex-oriented 4-regular planar graph on
 $n$ vertices in which any cut-vertex or loop-anchor has been oriented
 transversely. By Lemma \ref{lemma: prime is ok}, we may suppose that $G$ is
 3-edge-connected. 
 
 Suppose first of all that $G$ does have a cut-vertex $v$, so that $G$
 is as shown in Figure \ref{figure: transverse cut figure} (a). If either of
 $G_1$ or $G_2$ as shown in Figure \ref{figure: transverse cut figure} (b)
 contains a cut-vertex that is oriented nontransversely, then that
 vertex is a cut-vertex of $G$ oriented nontransversely, which is not
 possible. If either of $G_1$ or $G_2$ contains a loop-anchor $w$ that
 is oriented nontransversely, then in $G$, $w$ is either a cut-vertex or
 a loop-anchor that is oriented nontransversely, neither of which is possible. Thus by
 our inductive hypothesis, each of $G_1$ and $G_2$, with the
 vertex-orientations induced by $\sigma$, can be o-coloured, and then by
 Lemma \ref{lemma: transverse cut vertex can be o-coloured}, $G$ can be
 o-coloured. If $G$ contains a loop-anchor $v$, then $v$ is oriented
 transversely, in which case we can o-colour the vertex-oriented graph that
 is obtained from $(G,\sigma)$ by smoothing $v$, and consequently we can o-colour
 $(G,\sigma)$. Thus we may assume that $G$ has no loops or cut-vertices.
 
\noindent Case 1: $G$ contains a vertex $v$ such that $(G,\sigma)$ is of 
 the form shown in Figure \ref{figure: flyping parallel figure}
 (a). Smooth $v$ to form the vertex-oriented graphs $(G_1,\sigma_1)$ and
 $(G_2,\sigma_2)$ as shown in Figure \ref{figure: flyping parallel figure} (b).
 Neither can contain a cut-vertex or a loop-anchor, so by our induction
 hypothesis, each can be o-coloured. Then by Lemma \ref{lemma: flyping
 parallel}, $G$ can be o-coloured.
 
 We may therefore suppose that Case 1 does not occur.
 
\noindent Case 2: $G$ contains a vertex $v$ such that $(G,\sigma)$ is of
the form shown in Figure \ref{figure: transverse 4-tangle figure}  (a), where
each of the closed curves $S_1$ and $S_2$ contain at least one vertex in
their respective interiors. Smooth $v$ to form the vertex-oriented graph
$(G^{(1)},\sigma_1)$ as shown in Figure \ref{figure: transverse 4-tangle figure}
(b),  where the marked colours are for later reference.

 \begin{figure}[h]
  \centering\table{c}
  $\vcenter{\xy /r20pt/:,
    (0,0)="centre","centre"+(1,0)="s";"centre"+(0,1)="e"**\crv{~*=<2.5pt>{\hbox{\Large.}} "s"+(0,.45) & "centre"+(.72,.72) & "e"+(.45,0)},
    "centre"+(-1,0)="s";"centre"+(0,1)="e"**\crv{~*=<2.5pt>{\hbox{\Large.}} "s"+(0,.45) & "centre"+(-.72,.72) & "e"-(.45,0)},
    "centre"+(1,0)="s";"centre"+(0,-1)="e"**\crv{~*=<2.5pt>{\hbox{\Large.}} "s"+(0,-.45) & "centre"+(.72,-.72) & "e"+(.45,0)},
    "centre"+(-1,0)="s";"centre"+(0,-1)="e"**\crv{~*=<2.5pt>{\hbox{\Large.}} "s"+(0,-.45) & "centre"+(-.72,-.72) & "e"-(.45,0)},
    (4,0)="centre","centre"+(1,0)="s";"centre"+(0,1)="e"**\crv{~*=<2.5pt>{\hbox{\Large.}} "s"+(0,.45) & "centre"+(.72,.72) & "e"+(.45,0)},
    "centre"+(-1,0)="s";"centre"+(0,1)="e"**\crv{~*=<2.5pt>{\hbox{\Large.}} "s"+(0,.45) & "centre"+(-.72,.72) & "e"-(.45,0)},
    "centre"+(1,0)="s";"centre"+(0,-1)="e"**\crv{~*=<2.5pt>{\hbox{\Large.}} "s"+(0,-.45) & "centre"+(.72,-.72) & "e"+(.45,0)},
    "centre"+(-1,0)="s";"centre"+(0,-1)="e"**\crv{~*=<2.5pt>{\hbox{\Large.}} "s"+(0,-.45) & "centre"+(-.72,-.72) & "e"-(.45,0)},
   (.6,-.5)="dl";(3.4,.5)="ur"**\crv{(1.5,-.5) & (2,0) & (2.5,.5)},
   (.6,.5)="ul";(3.4,-.5)="dr"**\crv{(1.5,.5) & (2,0) & (2.5,-.5)},   
   "ul"+(-.3,.2)="tl",
   "ur"+(.3,.2)="tr",   
   "tl";"tr"**\crv{"tl"+(.2,.2) & "tr"+(-.2,.2)},
   "dl"+(-.3,-.2)="bl",
   "dr"+(.3,-.2)="br",   
   "bl";"br"**\crv{"bl"+(.2,-.2) & "br"+(-.2,-.2)},
   (2,0)*{\smallbullet}*!<-8pt,0pt>{v},   
   (2,-.45);(2,.45)**\dir{-}*\dir{>},(2,-.45)*\dir{<},
   (0,0)*!(1,-1){S_1},
   (4,0)*!(-1,-1){S_2},   
   (0,-1.8)*{},
   (0,1.7)*{},
   (8,0)="x",
    (0,0)+"x"="centre","centre"+(1,0)="s";"centre"+(0,1)="e"**\crv{~*=<2.5pt>{\hbox{\Large.}} "s"+(0,.45) & "centre"+(.72,.72) & "e"+(.45,0)},
    "centre"+(-1,0)="s";"centre"+(0,1)="e"**\crv{~*=<2.5pt>{\hbox{\Large.}} "s"+(0,.45) & "centre"+(-.72,.72) & "e"-(.45,0)},
    "centre"+(1,0)="s";"centre"+(0,-1)="e"**\crv{~*=<2.5pt>{\hbox{\Large.}} "s"+(0,-.45) & "centre"+(.72,-.72) & "e"+(.45,0)},
    "centre"+(-1,0)="s";"centre"+(0,-1)="e"**\crv{~*=<2.5pt>{\hbox{\Large.}} "s"+(0,-.45) & "centre"+(-.72,-.72) & "e"-(.45,0)},
    (4,0)+"x"="centre","centre"+(1,0)="s";"centre"+(0,1)="e"**\crv{~*=<2.5pt>{\hbox{\Large.}} "s"+(0,.45) & "centre"+(.72,.72) & "e"+(.45,0)},
    "centre"+(-1,0)="s";"centre"+(0,1)="e"**\crv{~*=<2.5pt>{\hbox{\Large.}} "s"+(0,.45) & "centre"+(-.72,.72) & "e"-(.45,0)},
    "centre"+(1,0)="s";"centre"+(0,-1)="e"**\crv{~*=<2.5pt>{\hbox{\Large.}} "s"+(0,-.45) & "centre"+(.72,-.72) & "e"+(.45,0)},
    "centre"+(-1,0)="s";"centre"+(0,-1)="e"**\crv{~*=<2.5pt>{\hbox{\Large.}} "s"+(0,-.45) & "centre"+(-.72,-.72) & "e"-(.45,0)},
   (.6,-.5)+"x"="dl";(.6,.5)+"x"="ul"**\crv{(1.3,-.5)+"x" & (1.5,0)+"x" & (1.3,.5)+"x"},   
   (1.5,0)+"x"*!<-4pt,-2pt>{c_1},
   (2.5,0)+"x"*!<3pt,-2pt>{c_2},   
    (3.4,-.5)+"x"="dr";(3.4,.5)+"x"="ur"**\crv{(2.8,-.5)+"x" & (2.5,0)+"x" & (2.8,.5)+"x"},   
   "ul"+(-.3,.2)="tl",
   "ur"+(.3,.2)="tr",   
   "tl";"tr"**\crv{"tl"+(.2,.2) & "tr"+(-.2,.2)},
   "dl"+(-.3,-.2)="bl",
   "dr"+(.3,-.2)="br",   
   "bl";"br"**\crv{"bl"+(.2,-.2) & "br"+(-.2,-.2)},
   "tl";"tr"**\dir{}?(.5)="labeltop",
   "labeltop"*!<0pt,-8pt>{c_3},
   "bl";"br"**\dir{}?(.5)="labelbot",
   "labelbot"*!<0pt,9pt>{c_3},
   (2,-2.5)*{\hbox{(a) }},
   (2,-2.3)+"x"*{\hbox{(b) $(G^{(1)},\sigma_1)$}},
  \endxy}$
  \endtable\caption{}\label{figure: transverse 4-tangle figure}
 \end{figure}
 
 \noindent If $(G^{(1)},\sigma_1)$ contains a cut-vertex or a
 loop-anchor $w$ oriented nontraversely, then in $G$, $w$ provides a
 Case 1 scenario, and we have assumed that there are no such vertices in
 $G$. Thus by our inductive hypothesis, there is an o-colouring of
 $(G^{(1)},\sigma_1)$, as shown in Figure \ref{figure: transverse 4-tangle
 figure} (b). Note that the top and bottom edges between the subgraphs
 enclosed by closed curves $S_1$ and $S_2$ must belong to the same
 o-cycle, and thus have the same colour, labelled $c_3$. There are three
 subcases to consider.
 
\noindent Case 2 (i): $c_1\ne c_2$. Then give each edge of $G$ that is also
an edge of $G^{(1)}$ the colour it received in the o-colouring
of $(G^{(1)},\sigma_1)$, and colour the edges incident to $v$ as shown
in Figure \ref{figure: transverse original coloured figure}. The result is an
o-colouring of $(G,\sigma)$.
 
 \begin{figure}[h]
  \centering\table{c}
  $\vcenter{\xy /r20pt/:,
    (0,0)="centre","centre"+(1,0)="s";"centre"+(0,1)="e"**\crv{~*=<2.5pt>{\hbox{\Large.}} "s"+(0,.45) & "centre"+(.72,.72) & "e"+(.45,0)},
    "centre"+(-1,0)="s";"centre"+(0,1)="e"**\crv{~*=<2.5pt>{\hbox{\Large.}} "s"+(0,.45) & "centre"+(-.72,.72) & "e"-(.45,0)},
    "centre"+(1,0)="s";"centre"+(0,-1)="e"**\crv{~*=<2.5pt>{\hbox{\Large.}} "s"+(0,-.45) & "centre"+(.72,-.72) & "e"+(.45,0)},
    "centre"+(-1,0)="s";"centre"+(0,-1)="e"**\crv{~*=<2.5pt>{\hbox{\Large.}} "s"+(0,-.45) & "centre"+(-.72,-.72) & "e"-(.45,0)},
    (4,0)="centre","centre"+(1,0)="s";"centre"+(0,1)="e"**\crv{~*=<2.5pt>{\hbox{\Large.}} "s"+(0,.45) & "centre"+(.72,.72) & "e"+(.45,0)},
    "centre"+(-1,0)="s";"centre"+(0,1)="e"**\crv{~*=<2.5pt>{\hbox{\Large.}} "s"+(0,.45) & "centre"+(-.72,.72) & "e"-(.45,0)},
    "centre"+(1,0)="s";"centre"+(0,-1)="e"**\crv{~*=<2.5pt>{\hbox{\Large.}} "s"+(0,-.45) & "centre"+(.72,-.72) & "e"+(.45,0)},
    "centre"+(-1,0)="s";"centre"+(0,-1)="e"**\crv{~*=<2.5pt>{\hbox{\Large.}} "s"+(0,-.45) & "centre"+(-.72,-.72) & "e"-(.45,0)},
   (.6,-.5)="dl";(3.4,.5)="ur"**\crv{(1.5,-.5) & (2,0) & (2.5,.5)},
   (.6,.5)="ul";(3.4,-.5)="dr"**\crv{(1.5,.5) & (2,0) & (2.5,-.5)},   
   "ul"+(-.3,.2)="tl",
   "ur"+(.3,.2)="tr",   
   "tl";"tr"**\crv{"tl"+(.2,.2) & "tr"+(-.2,.2)},
   "dl"+(-.3,-.2)="bl",
   "dr"+(.3,-.2)="br",   
   "bl";"br"**\crv{"bl"+(.2,-.2) & "br"+(-.2,-.2)},
   (2,0)*{\smallbullet}*!<-8pt,0pt>{v},   
   (2,-.45);(2,.45)**\dir{-}*\dir{>},(2,-.45)*\dir{<},
   (0,0)*!(1,-1){S_1},
   (4,0)*!(-1,-1){S_2},   
   "tl";"tr"**\dir{}?(.5)="labeltop",
   "labeltop"*!<0pt,-8pt>{c_3},
   "bl";"br"**\dir{}?(.5)="labelbot",
   "labelbot"*!<0pt,9pt>{c_3},
   "ul"+(.25,0)*!<-13pt,-1pt>{c_1},
   "ur"+(-.25,0)*!<14pt,-1.5pt>{c_2},   
   "dl"+(.25,0)*!<-13pt,2.5pt>{c_1},
   "dr"+(-.25,0)*!<14pt,1.5pt>{c_2},   
  \endxy}$
  \endtable\caption{}\label{figure: transverse original coloured figure}
 \end{figure}
  
\noindent Case 2 (ii): $c_1=c_2\ne c_3$. As shown in
Figure \ref{figure: transverse g2 figure} (a), form the vertex-oriented
4-regular planar graph $(G^{(2)},\sigma_2)$, where the
vertex-orientation is that induced  by $\sigma_1$, and give it the
o-colouring obtained from that of $(G^{(1)},\sigma_1)$ as shown in the
figure. Choose a third colour $c\ne c_2,c_3$ (so $c$ is a new colour if
the o-colouring of $(G^{(1)},\sigma_1)$ used only two colours), and in
this o-colouring of $(G^{(2)},\sigma_2)$, swap $c$ and $c_2$, so that
now $(G^{(2)},\sigma_2)$ is o-coloured as shown in Figure \ref{figure: transverse
g2 figure} (b).

 \begin{figure}[h]
  \centering\table{c@{\hskip60pt}c}
   $\vcenter{\xy /r20pt/:,
   (-2,0)="x",
   (1,0)="y",
   (4.5,2.1)+"x"+"y"*{},    
   (4.5,0)+"x"+"y"="centre",{(0,0);(1.5,0):,"centre"+(1,0)="s";"centre"+(0,1)="e"**\crv{~*=<2.5pt>{\hbox{\Large.}} "s"+(0,.45) & "centre"+(.72,.72) & "e"+(.45,0)},
    "centre"+(-1,0)="s";"centre"+(0,1)="e"**\crv{~*=<2.5pt>{\hbox{\Large.}} "s"+(0,.45) & "centre"+(-.72,.72) & "e"-(.45,0)},
    "centre"+(1,0)="s";"centre"+(0,-1)="e"**\crv{~*=<2.5pt>{\hbox{\Large.}} "s"+(0,-.45) & "centre"+(.72,-.72) & "e"+(.45,0)},
    "centre"+(-1,0)="s";"centre"+(0,-1)="e"**\crv{~*=<2.5pt>{\hbox{\Large.}} "s"+(0,-.45) & "centre"+(-.72,-.72) & "e"-(.45,0)}},
   (3.4,.5)+"x"+"y"="b",
   (3.4,-.5)+"x"+"y"="d",
   "d"+(.2,-.5)="br",
   "b"+(.2,.5)="tr",
   "tr";"br"**\crv{"tr"+(-2.5,0) & "br"+(-2.5,0)}?(.5)="t",   
   "d";"b"**\crv{"d"+(-1.2,0) & "b"+(-1.2,0)}?(.5)="c",      
   (4.5,0)+"x"+"y"*!(-1.5,-1.5){S_2},   
   "t"*!<6pt,0pt>{c_3},
   "c"*!<6pt,0pt>{c_2},
   (4.5,-2.3)+"x"+"y"*{\hbox{(a)\quad $(G^{(2)},\sigma_2)$}},   
   \endxy}$
  &
   $\vcenter{\xy /r20pt/:,
    (-2,0)="x",
    (1,0)="y",
    (4.5,0)+"x"+"y"="centre",{(0,0);(1.5,0):,"centre"+(1,0)="s";"centre"+(0,1)="e"**\crv{~*=<2.5pt>{\hbox{\Large.}} "s"+(0,.45) & "centre"+(.72,.72) & "e"+(.45,0)},
    "centre"+(-1,0)="s";"centre"+(0,1)="e"**\crv{~*=<2.5pt>{\hbox{\Large.}} "s"+(0,.45) & "centre"+(-.72,.72) & "e"-(.45,0)},
    "centre"+(1,0)="s";"centre"+(0,-1)="e"**\crv{~*=<2.5pt>{\hbox{\Large.}} "s"+(0,-.45) & "centre"+(.72,-.72) & "e"+(.45,0)},
    "centre"+(-1,0)="s";"centre"+(0,-1)="e"**\crv{~*=<2.5pt>{\hbox{\Large.}} "s"+(0,-.45) & "centre"+(-.72,-.72) & "e"-(.45,0)}},
   (3.4,.5)+"x"+"y"="b",
   (3.4,-.5)+"x"+"y"="d",
   "d"+(.2,-.5)="br",
   "b"+(.2,.5)="tr",
   "tr";"br"**\crv{"tr"+(-2.5,0) & "br"+(-2.5,0)}?(.5)="t",   
   "d";"b"**\crv{"d"+(-1.2,0) & "b"+(-1.2,0)}?(.5)="c",      
   (4.5,0)+"x"+"y"*!(-1.5,-1.5){S_2},   
   "t"*!<6pt,0pt>{c_3},
   "c"*!<6pt,0pt>{c},
   (4.5,-2.3)+"x"+"y"*{\hbox{(b)\quad $(G^{(2)},\sigma_2)$}},   
   \endxy}$
  \endtable\caption{}\label{figure: transverse g2 figure}
 \end{figure}

\noindent Lift this colouring back to $(G^{(1)},\sigma_1)$ as shown in
Figure \ref{figure: transverse rebuilt figure}. Then we are back in Case 2 (i),
and so $(G,\sigma)$ can be o-coloured.

 \begin{figure}[h]
  \centering\table{c}
  $\vcenter{\xy /r20pt/:,
   (0,0)="x",
   (0,0)+"x"="centre","centre"+(1,0)="s";"centre"+(0,1)="e"**\crv{~*=<2.5pt>{\hbox{\Large.}} "s"+(0,.45) & "centre"+(.72,.72) & "e"+(.45,0)},
   "centre"+(-1,0)="s";"centre"+(0,1)="e"**\crv{~*=<2.5pt>{\hbox{\Large.}} "s"+(0,.45) & "centre"+(-.72,.72) & "e"-(.45,0)},
   "centre"+(1,0)="s";"centre"+(0,-1)="e"**\crv{~*=<2.5pt>{\hbox{\Large.}} "s"+(0,-.45) & "centre"+(.72,-.72) & "e"+(.45,0)},
   "centre"+(-1,0)="s";"centre"+(0,-1)="e"**\crv{~*=<2.5pt>{\hbox{\Large.}} "s"+(0,-.45) & "centre"+(-.72,-.72) & "e"-(.45,0)},
   (4,0)+"x",{\ellipse(1){.}},
   (4,0)+"x"="centre","centre"+(1,0)="s";"centre"+(0,1)="e"**\crv{~*=<2.5pt>{\hbox{\Large.}} "s"+(0,.45) & "centre"+(.72,.72) & "e"+(.45,0)},
   "centre"+(-1,0)="s";"centre"+(0,1)="e"**\crv{~*=<2.5pt>{\hbox{\Large.}} "s"+(0,.45) & "centre"+(-.72,.72) & "e"-(.45,0)},
   "centre"+(1,0)="s";"centre"+(0,-1)="e"**\crv{~*=<2.5pt>{\hbox{\Large.}} "s"+(0,-.45) & "centre"+(.72,-.72) & "e"+(.45,0)},
   "centre"+(-1,0)="s";"centre"+(0,-1)="e"**\crv{~*=<2.5pt>{\hbox{\Large.}} "s"+(0,-.45) & "centre"+(-.72,-.72) & "e"-(.45,0)},
   (.6,-.5)+"x"="dl";(.6,.5)+"x"="ul"**\crv{(1.3,-.5)+"x" & (1.5,0)+"x" & (1.3,.5)+"x"},   
   (1.5,0)+"x"*!<-4pt,-2pt>{c_1},
   (2.5,0)+"x"*!<3pt,-2pt>{c},   
    (3.4,-.5)+"x"="dr";(3.4,.5)+"x"="ur"**\crv{(2.8,-.5)+"x" & (2.5,0)+"x" & (2.8,.5)+"x"},   
   "ul"+(-.3,.2)="tl",
   "ur"+(.3,.2)="tr",   
   "tl";"tr"**\crv{"tl"+(.2,.2) & "tr"+(-.2,.2)},
   "dl"+(-.3,-.2)="bl",
   "dr"+(.3,-.2)="br",   
   "bl";"br"**\crv{"bl"+(.2,-.2) & "br"+(-.2,-.2)},
   "tl";"tr"**\dir{}?(.5)="labeltop",
   "labeltop"*!<0pt,-8pt>{c_3},
   "bl";"br"**\dir{}?(.5)="labelbot",
   "labelbot"*!<0pt,9pt>{c_3},
   (0,0)*!(1,-1){S_1},
   (4,0)*!(-1,-1){S_2},   
  \endxy}$
  \endtable\caption{}\label{figure: transverse rebuilt figure}
 \end{figure}

\noindent Case 2 (iii): $c_1=c_2=c_3$. Suppose first that in the o-colouring
of $(G^{(1)},\sigma_1)$ as shown in Figure \ref{figure: transverse 4-tangle figure} (b),
the edges $e$ and $f$ with colour labels $c_1$ and $c_2$, respectively, do not
belong to the same o-cycle. Then the two edges $u$ (up) and $d$ (down) shown in 
Figure \ref{figure: transverse 4-tangle figure} (b) with colour
label $c_3$ must belong to the same o-cycle, $O$ say, and not both 
$e$ and $f$ can belong to $O$. Without loss of generality, suppose 
that $f$ is not in $O$. Then we may permute the colours of the edges
that appear in the interior of $S_2$ other than those that belong to
$O$ in such a way that $f$ is not coloured with colour $c_1$ (if the
edges of $S_2$ had been coloured with only two colours, then a third
colour would need to be introduced). This would then place us in the 
context of Case 2 (i), and so $(G,\sigma)$ is o-colourable.

Suppose now that $e$ and $f$ belong to the same o-cycle, which
is then the o-cycle $O$ that contains $u$ and $d$. Since at least one
vertex of $G$ is contained in the interior of $S_2$, $S_2$ must contain
at least one o-cycle in addition to $O$, so let $C$ be an o-cycle
contained in the interior of $S_2$. Since $G$ contains no loop-anchors,
$C$ must pass through at least two vertices. Remove the edges of $C$ from $G$.
Now each vertex of $C$ has two incident edges, and
both have the same colour, so we may remove the vertex and identify the
two edges, giving this new edge the common colour of the two that have
been identified. Denote the resulting vertex-oriented graph by $(G^{(3)},\sigma_3)$,
and note that the o-colouring of $(G^{(1)},\sigma_1)$ provides an o-colouring
of $(G^{(3)},\sigma_3)$. As a result, any cut-vertex of $G^{(3)}$ must be
oriented transversely by $\sigma_3$. It follows therefore that if we modify
$G^{(3)}$ by re-introducing $v$, calling the vertex-oriented result 
$(G'',\sigma'')$, then the only vertex of $G''$ that could possibly
be a cut-vertex oriented nontransversely is $v$. Suppose that in fact,
$v$ is a nontransversely oriented cut-vertex of $(G'',\sigma'')$. 
Then there exist simple closed curves $U_1$ and $U_2$, as shown in Figure 
\ref{figure: modified G'}, such that one of the two components of $G''-v$
is contained within $U_1$ and the other component is contained within $U_2$.

  \begin{figure}[h]
  \centering\table{c}
  $\vcenter{\xy /r20pt/:,
   (0,0)="centre","centre"+(1,0)="s";"centre"+(0,1)="e"**\crv{~*=<2.5pt>{\hbox{\Large.}} "s"+(0,.45) & "centre"+(.72,.72) & "e"+(.45,0)},
   "centre"+(-1,0)="s";"centre"+(0,1)="e"**\crv{~*=<2.5pt>{\hbox{\Large.}} "s"+(0,.45) & "centre"+(-.72,.72) & "e"-(.45,0)},
   "centre"+(1,0)="s";"centre"+(0,-1)="e"**\crv{~*=<2.5pt>{\hbox{\Large.}} "s"+(0,-.45) & "centre"+(.72,-.72) & "e"+(.45,0)},
   "centre"+(-1,0)="s";"centre"+(0,-1)="e"**\crv{~*=<2.5pt>{\hbox{\Large.}} "s"+(0,-.45) & "centre"+(-.72,-.72) & "e"-(.45,0)},
   (4,0)="centre","centre"+(1,0)="s";"centre"+(0,1)="e"**\crv{~*=<2.5pt>{\hbox{\Large.}} "s"+(0,.45) & "centre"+(.72,.72) & "e"+(.45,0)},
   "centre"+(-1,0)="s";"centre"+(0,1)="e"**\crv{~*=<2.5pt>{\hbox{\Large.}} "s"+(0,.45) & "centre"+(-.72,.72) & "e"-(.45,0)},
   "centre"+(1,0)="s";"centre"+(0,-1)="e"**\crv{~*=<2.5pt>{\hbox{\Large.}} "s"+(0,-.45) & "centre"+(.72,-.72) & "e"+(.45,0)},
   "centre"+(-1,0)="s";"centre"+(0,-1)="e"**\crv{~*=<2.5pt>{\hbox{\Large.}} "s"+(0,-.45) & "centre"+(-.72,-.72) & "e"-(.45,0)},
   (.6,-.5)="dl";(3.4,.5)="ur"**\crv{(1.5,-.5) & (2,0) & (2.5,.5)},
   (.6,.5)="ul";(3.4,-.5)="dr"**\crv{(1.5,.5) & (2,0) & (2.5,-.5)},   
   "ul"+(-.3,.2)="tl",
   "ur"+(.3,.2)="tr",   
   "tl";"tr"**\crv{"tl"+(.2,.2) & "tr"+(-.2,.2)},
   "dl"+(-.3,-.2)="bl",
   "dr"+(.3,-.2)="br",   
   "bl";"br"**\crv{"bl"+(.2,-.2) & "br"+(-.2,-.2)},
   (2,0)*{\smallbullet}*!<-8pt,0pt>{v},   
   (2,-.45);(2,.45)**\dir{-}*\dir{>},(2,-.45)*\dir{<},
   (0,0)*!(1,-1){S_1},
   (4,0)*!(-1,-1){S_2},   
   "tl";"tr"**\dir{}?(.5)="labeltop",
   "bl";"br"**\dir{}?(.5)="labelbot",
   "ur"+(.85,0)="a";"dr"+(.85,0)="b"**\crv{~*=<2.5pt>{\hbox{\LARGE.}} "a"+(.8,0) & "b"+(.8,0) },
   "a";"b"**\crv{~*=<2.5pt>{\hbox{\LARGE.}} "a"+(-.8,0) & "b"+(-.8,0) },   
   "a"*!<0pt,5pt>{C},
   (0,0)+(0,.25)="s";(4,0)+(0,.25)="e"**\crv{~*=<2.5pt>{\hbox{\LARGE.}} "s"+(-2,0) & (0,0)+(-2,2) & (4,0)+(2,2) & "e"+(2,0)},
   (2,0)+(0,2.4)*{U_1},
   "s";"e"**\crv{~*=<2.5pt>{\hbox{\LARGE.}} "s"+(1,0) & (2,0)+(0,1) & "e"+(-1,0)},
   (0,0)+(0,-.25)="s";(4,0)+(0,-.25)="e"**\crv{~*=<2.5pt>{\hbox{\LARGE.}} "s"+(-2,0) & (0,0)+(-2,-2) & (4,0)+(2,-2) & "e"+(2,0)},
   (2,0)+(0,-2.4)*{U_2},
   "s";"e"**\crv{~*=<2.5pt>{\hbox{\LARGE.}} "s"+(1,0) & (2,0)+(0,-1) & "e"+(-1,0)},
  \endxy}$
  \endtable\caption{$(G'',\sigma'')$}\label{figure: modified G'}
 \end{figure}

 \noindent Since $S_1$ does contain vertices of $G$, we have a contradiction to
 the fact that $G$ is 3-edge-connected. Thus in $G''$, $v$ is oriented transversely
 by $\sigma''$.  As $C$ contained at 
 least two vertices, the number of vertices in $G'$ is at least two fewer than the number of 
 vertices in $G$ and thus $G''$ contains fewer vertices than $G$. We may therefore apply the 
 induction hypothesis to $(G'',\sigma'')$ to obtain an o-colouring 
 of $(G'',\sigma'')$. Finally, reintroduce the vertices and edges of the o-cycle $C$, colouring
 the edges of $C$ with a new colour if necessary. The result is an o-colouring of $(G,\sigma)$.
 
 \noindent Case 3: no vertex in $(G,\sigma)$ is of the type in either 
 Case 1 or Case 2. In particular, $G$  must be simple (it was loopfree,
 and since we are not in Case 1 or Case 2, there are no multiple edges).
 Furthermore, each vertex of $G$ can be smoothed without creating either
 a cut-vertex or a loop-anchor oriented nontransversely (since $G$ was loopfree, such a
 vertex would establish that $(G,\sigma)$ belonged in Case 1 or Case
 2). Choose any vertex $v$ and smooth it, thereby obtaining a 4-regular
 planar vertex-oriented graph $(G',\sigma')$  on $n-1$ vertices with no
 cut-vertices or loop-anchors, so by hypothesis, there is an o-colouring
 for this graph. Consider a particular o-colouring of this graph. If the
 two edges that resulted from the smoothing of $v$ belong to different
 o-cycles, then the o-colouring lifts to an o-colouring of $(G,\sigma)$
 (they may be coloured the same, but since they are different o-cycles,
 we may then change the colour of one, possibly requiring a new colour).
 Thus we may assume that the two edges that resulted from the smoothing
 belong to the same o-cycle, which we shall denote by $C_0$. Every other
 o-cycle of this o-colouring of $(G',\sigma')$ is an o-cycle of
 $(G,\sigma)$, while the edges of $C_0$ other than the two edges of the
 smoothing, together with the four edges incident to $v$, form two
 cycles in $G$, $C_1$ and $C_1'$ say, that meet only at $v$, and which
 meet the o-cycle requirement at every vertex except $v$. If the removal
 of any one of the o-cycles other than $C_0$ from $G$ results in a graph
 with no cut-vertex or loop-anchor oriented nontransversely, then we
 could o-colour the result and reinsert the o-cycle, giving it a new
 colour if necessary, thereby obtaining an o-colouring of $(G,\sigma)$.
 Suppose then that the removal of any of these o-cycles other than $C_0$
 from $G$ results in a cut-vertex or loop-anchor oriented
 nontransversely. Since the removal of the same o-cycle from $G'$ does
 not result in such a vertex (since $(G',\sigma')$ was o-coloured), we
 see that the vertex that has become a   nontransversely oriented
 cut-vertex or loop-anchor is $v$. Consider the abstract graph whose
 vertices are the cycles  $C_1$, $C_1'$, and the o-cycles of the
 o-colouring of $(G',\sigma')$ other than $C_0$. Two vertices of this
 graph are to be joined by an edge if they have a vertex of $G$ in
 common. This graph is connected with the property that  every vertex
 other than $C_1$ and $C_1'$ lies on every path in this graph from $C_1$
 to $C_1'$. Thus this graph is a chain with endpoints $C_1$ and $C_1'$.
 Note that there is at least one o-cycle in this chain. Begin at $C_1$
 and follow this chain, labelling each vertex on the chain (o-cycle, or
 at the end, $C_1'$) as $C_i$, $i=1,2,\ldots,m+1$, where 
 $C_{m+1}=C_1'$. Then for any plane embedding of $G$, there exist simple
 closed curves $S_1$, $S_2,\ldots,S_{m}$ such that for each
 $i=1,2,\ldots,m$, all vertices in common to $C_i$ and $C_{i+1}$ lie
 within $S_i$ and no other vertices of $G$ lie within $S_i$,  and any
 edge joining two vertices of $G$ that lie within $S_i$ also lies within
 $S_i$. Now, every vertex of $G$ other than $v$ lies within one and only
 one $S_i$, and for each $i=1,2,\ldots,m$, we shall let $G_i$ denote the
 subgraph of $G$ that is induced by the vertices lying within $S_i$.
 Additionally, let $G_0=G_{m+1}$ denote the null graph whose only vertex
 is $v$. Note that any edge not contained within any of these simple
 closed curves must join a vertex in $G_i$ to a vertex in $G_{i+1}$ for
 some $i$. 
 
 We shall demonstrate that there is at least one vertex $v$ such that,
 when smoothed, $m=2$; that is, there are two subgraphs $G_1$ and $G_2$, and three
 cycles $C_1$, $C_2$, and $C_3$, to use the notation introduced above.
 To do this, we shall examine the 3-faces in $G$, of which there must
 be at least eight.
 
 Suppose first that $G$ has at least one 3-face with orientation as shown
 in Figure \ref{figure: 3-faces} (a). Choose any vertex $v$ not belonging
 to the 3-face boundary and smooth it. By hypothesis, the resulting graph
 can be o-coloured. Since no two edges of the 3-face can be coloured the
 same, it follows that no two of the edges of the 3-face belong to the same
 o-cycle. As observed above, this means that for some $i$, the o-cycles
 $C_{i-1}$, $C_i$, and $C_{i+1}$ each have an edge on the 3-face. But this
 means that $C_{i-1}$ and $C_{i+1}$ have a vertex in common, which is not
 possible. Thus no 3-face of $G$ can be as in Figure \ref{figure: 3-faces} (a).
 
 \vskip10pt
 \begin{figure}[h]
  \centering\table{c@{\hskip20pt}c@{\hskip20pt}c@{\hskip20pt}c}
   \xy /r30pt/:,
   (0,0)="u"*{\smallbullet},
   (1,0)="w"*{\smallbullet},
   (.5,.866)="v"*{\smallbullet},
   "u";"w"**\dir{-};"v"**\dir{-};"u"**\dir{-},
   "u";"u"+a(30)**\dir{}?(.35)="u1",
   "u";"u1"**\dir{-}*\dir{>},
   "u";"u"+"u"-"u1"**\dir{-}*\dir{>},
   "w";"w"+a(150)**\dir{}?(.35)="w1",
   "w";"w1"**\dir{-}*\dir{>},
   "w";"w"+"w"-"w1"**\dir{-}*\dir{>},
   "v";"v"+(0,.35)**\dir{-}*\dir{>},
   "v";"v"-(0,.35)**\dir{-}*\dir{>},   
\endxy
 &
\xy /r30pt/:,
   (0,0)="u"*{\smallbullet},
   (1,0)="w"*{\smallbullet},
   (.5,.866)="v"*{\smallbullet},
   "u";"w"**\dir{-};"v"**\dir{-};"u"**\dir{-},
   "u";"u"+a(30)**\dir{}?(.35)="u1",
   "u";"u1"**\dir{-}*\dir{>},
   "u";"u"+"u"-"u1"**\dir{-}*\dir{>},
   "w";"w"+a(150)**\dir{}?(.35)="w1",
   "w";"w1"**\dir{-}*\dir{>},
   "w";"w"+"w"-"w1"**\dir{-}*\dir{>},
   "v";"v"+(.35,0)**\dir{-}*\dir{>},
   "v";"v"-(.35,0)**\dir{-}*\dir{>},   
   "v"*!<0pt,-7pt>{v},
\endxy
&
\xy /r30pt/:,
   (0,0)="u"*{\smallbullet},
   (1,0)="w"*{\smallbullet},
   (.5,.866)="v"*{\smallbullet},
   "u";"w"**\dir{-};"v"**\dir{-};"u"**\dir{-},
   "u";"u"+a(-60)**\dir{}?(.35)="u1",
   "u";"u1"**\dir{-}*\dir{>},
   "u";"u"+"u"-"u1"**\dir{-}*\dir{>},
   "w";"w"+a(150)**\dir{}?(.35)="w1",
   "w";"w1"**\dir{-}*\dir{>},
   "w";"w"+"w"-"w1"**\dir{-}*\dir{>},
   "v";"v"+(.35,0)**\dir{-}*\dir{>},
   "v";"v"-(.35,0)**\dir{-}*\dir{>},   
\endxy
&
\xy /r30pt/:,
   (0,0)="u"*{\smallbullet},
   (1,0)="w"*{\smallbullet},
   (.5,.866)="v"*{\smallbullet},
   "u";"w"**\dir{-};"v"**\dir{-};"u"**\dir{-},
   "u";"u"+a(-60)**\dir{}?(.35)="u1",
   "u";"u1"**\dir{-}*\dir{>},
   "u";"u"+"u"-"u1"**\dir{-}*\dir{>},
   "w";"w"+a(60)**\dir{}?(.35)="w1",
   "w";"w1"**\dir{-}*\dir{>},
   "w";"w"+"w"-"w1"**\dir{-}*\dir{>},
   "v";"v"+(.35,0)**\dir{-}*\dir{>},
   "v";"v"-(.35,0)**\dir{-}*\dir{>},   
\endxy\\
  \noalign{\vskip4pt}
  (a) & (b) & (c) & (d)
  \endtable\caption{}\label{figure: 3-faces}
 \end{figure}

 Suppose now that $G$ has a 3-face as in Figure \ref{figure: 3-faces} (b). Choose vertex $v$
 as shown in (b) and smooth it. By hypothesis, the resulting graph may be o-coloured. The Case 3
 restrictions mandate that the two new edges that resulted from smoothing $v$ are necessarily on
 the same o-cycle $C_i$, as shown in Figure \ref{figure: 3-face b} (a), where $C_i\ne C_j$. But
 then we may exchange the colours on the two arcs as shown in Figure \ref{figure: 3-face b} (b),
 thereby obtaining an o-colouring of the graph in which the two new arcs that resulted from smoothing
 $v$ have different colours, which is not possible. Thus no 3-face of $G$ can be as in Figure
 \ref{figure: 3-faces} (b), which means that every 3-face of $G$ is of the form shown in Figure
 \ref{figure: 3-faces} (c) or (d).

 \vskip10pt
 \begin{figure}[h]
  \centering\table{c@{\hskip20pt}c}
 \xy /r30pt/:,
   (-1.3,0)*{\null}, 
   (-.3,.5)="v1",(.3,.5)="v2",
   "v1";"v2"**\crv{(0,.2)},
   (-.6,-.2)="u1",(.6,-.2)="u2",
   "u1";"u2"**\dir{-},
   "u1"+(0,-.4)="w1",
   "u2"+(0,-.4)="w2",   
   "w1";"w2"**\crv{(0,1)},
   "u1"+(.18,0)="t1"*{\smallbullet},
   "u2"-(.18,0)="t2"*{\smallbullet},   
   "v2"*!<-8pt,0pt>{C_i},
   "u2"*!<0pt,-10pt>{C_i},   
   "u1";"u2"**\dir{}?(.5)="x",
   "x"*!<0pt,10pt>{C_j},
   "t1";"t1"+a(30)**\dir{}?(.35)="t1a",
   "t1";"t1a"**\dir{-}*\dir{>},
   "t1";"t1"+"t1"-"t1a"**\dir{-}*\dir{>},
   "t2";"t2"+a(-30)**\dir{}?(.35)="t2a",
   "t2";"t2a"**\dir{-}*\dir{>},
   "t2";"t2"+"t2"-"t2a"**\dir{-}*\dir{>},
\endxy
&
 \xy /r30pt/:,
   (-1.15,0)*{\null},
   (-.3,.5)="v1",(.3,.5)="v2",
   "v1";"v2"**\crv{(0,.2)},
   (-.6,-.2)="u1",(.6,-.2)="u2",
   "u1";"u2"**\dir{-},
   "u1"+(0,-.4)="w1",
   "u2"+(0,-.4)="w2",   
   "w1";"w2"**\crv{(0,1)},
   "u1"+(.18,0)="t1"*{\smallbullet},
   "u2"-(.18,0)="t2"*{\smallbullet},   
   "v2"*!<-8pt,0pt>{C_i},
   "u2"*!<0pt,-10pt>{C_j},   
   "u1";"u2"**\dir{}?(.5)="x",
   "x"*!<0pt,10pt>{C_i},
   "t1";"t1"+a(30)**\dir{}?(.35)="t1a",
   "t1";"t1a"**\dir{-}*\dir{>},
   "t1";"t1"+"t1"-"t1a"**\dir{-}*\dir{>},
   "t2";"t2"+a(-30)**\dir{}?(.35)="t2a",
   "t2";"t2a"**\dir{-}*\dir{>},
   "t2";"t2"+"t2"-"t2a"**\dir{-}*\dir{>},
\endxy\\
  \noalign{\vskip4pt}
  (a) & (b) 
  \endtable\caption{}\label{figure: 3-face b}
 \end{figure}

 Choose a 3-face as shown in Figure \ref{figure: G_1 and G_2}, with either orientation at 
 $w$ (it is in fact possible to prove that there can be no 3-face of the type shown in Figure
 \ref{figure: 3-faces} (d), but this is not necessary for our argument), and smooth $v$. Then
 $u$ is in $G_1$, and $w$ is in $G_m$. As $u$ and $w$ are adjacent, it follows that $m=2$, as
 desired. We have o-cycle $C_1$ contained entirely within $G_1$, except for the two edges incident to
 $v$, one of which is $e=vu$, o-cycle $C_3$ contained entirely within $G_2$ except for the two edges incident
 to $v$, one of which is $vw$, and o-cycle $C_2$, which has edges incident to vertices of $G_1$ and
 to vertices of $G_2$. 

 \vskip10pt
 \begin{figure}[h]
  \centering\table{c}
\xy /r30pt/:,
   (0,0)="u"*{\smallbullet},
   (1,0)="w"*{\smallbullet},
   (.5,.866)="v"*{\smallbullet},
   "u";"w"**\dir{-};"v"**\dir{-};"u"**\dir{-},
   "u";"u"+a(-60)**\dir{}?(.35)="u1",
   "u";"u1"**\dir{-}*\dir{>},
   "u";"u"+"u"-"u1"**\dir{-}*\dir{>},
   "v";"v"+(.35,0)**\dir{-}*\dir{>},
   "v";"v"-(.35,0)**\dir{-}*\dir{>},   
   "v"*!<0pt,-7pt>{v},
   "u"*!<6pt,6pt>{u},   
   "w"*!<-6pt,6pt>{w},   
\endxy
  \endtable\caption{}\label{figure: G_1 and G_2}
 \end{figure}

 Let $\it O$ denote the set of all triples $(P_1,P_2,P_3)$, where $P_1$ is a 
 subpath of $C_1$ with initial vertex $v$ and initial edge $e$, $P_2$ is a subpath of $C_2$, and $P_3$
 is a subpath of $C_3$ such that the terminal vertex of $P_1$ is the initial vertex of $P_2$, the
 terminal vertex of $P_2$ is the initial vertex of $P_3$, and $P_1+P_2+P_3$ is an o-cycle. We show first
 that $\it O$ is not empty. Let $R_1$ be the o-path of length 1 with initial vertex $v$, initial edge
 $e$, and terminal vertex $u$. Let $R_2$ denote the o-path of length 1 with initial vertex $u$ and initial
 edge $uw$, so $R_2$ has terminal vertex $w$. Note that $R_2$ is a subpath of $C_2$ and that $R_1+R_2$ is
 defined and is an o-path. Finally, let $R_3$ denote the o-path with initial vertex $w$ and which follows
 $C_3$ in the direction which will make $R_2+R_3$ an o-path (this is uniquely determined), terminating at $v$.
 Thus $R_1+R_2+R_3$ is defined and is an o-cycle, so $\it O_1=(R_1,R_2,R_3)\in {\it O}$.
 
 For o-paths $P$ and $Q$, we shall say that $P\le Q$ if $P$ is a subpath of $Q$. This defines a partial order
 relation on the set of all o-paths in $G$. Now consider the lexical order relation on $\it O$ that is defined
 by this partial order relation on o-paths; that is, for $(P_1,P_2,P_3),(Q_1,Q_2,Q_3)\in \it O$,  we have
 $(P_1,P_2,P_3)<(Q_1,Q_2,Q_3)$ if $P_1<Q_1$, or else $P_1=Q_1$ and $P_2<Q_2$ (we note that if $P_1=Q_1$ and
 $P_2=Q_3$, then necessarily $P_3=Q_3$). We claim that this is a total order relation on $\it O$. For
 let $(P_1,P_2,P_3),(Q_1,Q_2,Q_3)\in \it O$, and suppose without loss of generality that $P_1\le Q_1$.
 If $P_1<Q_1$, then $(P_1,P_2,P_3)<(Q_1,Q_2,Q_3)$, so suppose that $P_1=Q_1$. Then both $P_2$ and $Q_2$
 have the same initial vertex, and travel along $C_2$ in the same direction. Thus we have exactly one of
 $P_2=Q_2$, $P_2<Q_2$, or $Q_2<P_2$. In every case, $(P_1,P_2,P_3)$ and $(Q_1,Q_2,Q_3)$ are comparable.
 
 Thus $\it O$ is a finite chain, in fact with minimum element $\it O_1$ defined above. Suppose that there
 are $t$ elements in the chain. Label the remaining $t-1$ as $\it O_2,\ldots,\it O_t$, so that for any
 $i$ and $j$ with $1\le i<j\le t$, we have $\it O_i<\it O_j$. For each $\it O_i=(P_1,P_2,P_3)$, let
 ${\it \hat{O}}_i$ denote the o-cycle $P_1+P_2+P_3$.
 
 If for some $i$, $G-E(\it \hat{O}_i)$ (where we shall think of the vertices of degree 2 as having been removed by an elementary
 subdivision operation) has no non-transversal cut-vertex, then by our induction hypothesis, $G-E(\it \hat{O}_i)$ may be
 o-coloured, in which case we may re-introduce the edges of
 $\it \hat{O}_i$, and colour them with a colour that is different from that used at
 any vertex of $\it\hat{O}_i$ (this may require introducing a new colour). The
 result is an o-colouring of $G$. Suppose to the contrary that for every $i$, $G-E(\it\hat{O}_i)$ has 
 at least one non-transversally oriented cut-vertex. 

 We shall say that $(P_1,P_2,P_3)\in \it O$ satisfies Condition A if $P_1$ meets $C_2$ only at vertices
 of one of the two arcs of $C_2$ that are determined by
 $u$ and the terminal vertex of $P_1$, and $C_1-P_1$ does not meet $C_2$ at any vertex of this arc.
 We note that $\it O_1$ trivially satisfies Condition A.

 Suppose now that $\it O_i=(P_1,P_2,P_3)\in \it O$ satisfies Condition A, and that $G-E(\it \hat{O}_i)$
 has a non-transversally oriented cut-vertex $z$ belonging to $G_1$. Let
 $U_1$ and $U_2$ denote the two components of $G-E(\it \hat O_i)-z$. Then the $C_1$ arc and the $C_2$ arc
 determined by one orientation cell at $z$ enter $U_1$, while the $C_1$ arc and the $C_2$ arc
 determined by the other orientation cell at $z$ enter $U_2$. One of the $C_2$ arcs must meet the terminal vertex
 of $P_1$, $x$ say, and we shall suppose that $U_1$ and $U_2$ are labelled so that $x$ is in $U_1$. Thus the $C_1$ 
 and $C_2$ arcs entering $U_2$ must meet $v$ and the terminal vertex of $P_2$, $y$ say, respectively. As $C_1-P_1$ meets
 $C_2$ at $z$, it follows by Condition A that $P_1$ can only meet $C_2$ at vertices on the arc of $C_2$ between
 $u$ and $x$ which does not contain $z$, and that $C_1-P_1$ does not meet this arc of $C_2$ (see Figure \ref{figure: bad guy 1}
 for a schematic diagram of this situation, with very few actual crossings depicted). 

 \vskip10pt
 \begin{figure}[h]
  \centering\table{c}
\hbox{\xy /r30pt/:,
   (-1,0)="x"*{\smallbullet}*!<6pt,0pt>{x},
   (1,0)="z"*{\smallbullet}*!<0pt,-6pt>{z},
   (2.3,.566)="v"*{\smallbullet}*!<0pt,-6pt>{v},
   (3,1)="u"*{\smallbullet}*!<-2pt,-6pt>{u},
   (3.4,-.2)="y"*{\smallbullet}*!<-5pt,7pt>{y},
   "x"+(0,-.1);"v"+(.1,.1)**\crv{~*=<8pt>{}~**!/1pt/\dir{-}  "x"+(-.25,-1) & (2.5,-1.5) & "u"+(1.5,-1.25) & "u"+(1.5,.3) & "u"+(0,.15)},
   "x"+(2.5,-1.5)*{P_1},
   "x"+(.0,-.1);"y"+(-.05,.05)**\crv{~*=<8pt>{}~**!/1pt/\dir{-}  "x"+(.25,-.5) & "x"+(2.4,-1) & "y"+(-.5,.2) "y"+(-.1,.1)},
   "x"+(2.45,-.95)*{P_2},
   "x";"z"**\crv{ "x"+(0,.3) & "x"+(.75,.75) & "z"+(-.3,-.5)}, 
   "x"+(.75,.85)*{C_2},
   "x";"z"**\crv{ "x"+(0,.75) & "x"+(.75,1)},  
   "x"+(.95,0)*{C_1},
   "z";"u"**\crv{ "z"+(.4,.2) &"z"+(.5,.75) & "u"+(-.7,.3)},
   "z"+(.4,.8)*{C_2},   
   "z";"v"**\crv{"z"+(.25,-.5) & "v"+(-.3,-.2)},
   "u";"y"**\crv{"u"+(.5,-.2) & "u"+(1.4,-.5) & "y"+(.4,-.2)},
   "v"+(.1,.01);"y"**\crv{~*=<7pt>{}~**!/1pt/\dir{-} "v"+(.3,.1) & "v"+(1.3,0) & "y"+(.5,.3)},   
   "v"+(.95,-.25)*{P_3},
   "v";"y"**\crv{ "v"+(-.3,-.4) & "v"+(-.6,-1) & "y"+(-.3,-.3)},
   "v"+(.05,-.45)*{C_3},   
   "z";"z"+<9pt,0pt>**\dir{-}*\dir{>},
   "z";"z"-<9pt,0pt>**\dir{-}*\dir{>},   
   (0,0);(1.5,0):,
   (-.7,0)="x",
   (0,0)+"x"="centre","centre"+(1,0)="s";"centre"+(0,1)="e"**\crv{~*=<2.5pt>{\hbox{\Large.}} "s"+(0,.45) & "centre"+(.72,.72) & "e"+(.45,0)},
   "centre"+(-1,0)="s";"centre"+(0,1)="e"**\crv{~*=<2.5pt>{\hbox{\Large.}} "s"+(0,.45) & "centre"+(-.72,.72) & "e"-(.45,0)},
   "centre"+(1,0)="s";"centre"+(0,-1)="e"**\crv{~*=<2.5pt>{\hbox{\Large.}} "s"+(0,-.45) & "centre"+(.72,-.72) & "e"+(.45,0)},
   "centre"+(-1,0)="s";"centre"+(0,-1)="e"**\crv{~*=<2.5pt>{\hbox{\Large.}} "s"+(0,-.45) & "centre"+(-.72,-.72) & "e"-(.45,0)},
   (-1.2,1.1)*{U_1},
   (2.1,0)="x",
   (0,0)+"x"="centre","centre"+(1,0)="s";"centre"+(0,1)="e"**\crv{~*=<2.5pt>{\hbox{\Large.}} "s"+(0,.45) & "centre"+(.72,.72) & "e"+(.45,0)},
   "centre"+(-1,0)="s";"centre"+(0,1)="e"**\crv{~*=<2.5pt>{\hbox{\Large.}} "s"+(0,.45) & "centre"+(-.72,.72) & "e"-(.45,0)},
   "centre"+(1,0)="s";"centre"+(0,-1)="e"**\crv{~*=<2.5pt>{\hbox{\Large.}} "s"+(0,-.45) & "centre"+(.72,-.72) & "e"+(.45,0)},
   "centre"+(-1,0)="s";"centre"+(0,-1)="e"**\crv{~*=<2.5pt>{\hbox{\Large.}} "s"+(0,-.45) & "centre"+(-.72,-.72) & "e"-(.45,0)},
   (2.5,1.1)*{U_2},
\endxy}
  \endtable\caption{}\label{figure: bad guy 1}
 \end{figure}

 \noindent Let $P_1^{(1)}$ denote the o-path obtained by extending $P_1$ (following $C_1$) to $z$, let $P_2^{(1)}$ denote the
 o-path obtained by following $C_2$ from $z$ into $U_2$, stopping at the first encountered vertex, $r$ say, on $C_3$, and
 let $P_3^{(1)}$ denote the o-path obtained by following $C_3$ from $r$ to $v$, travelling in the correct direction on $C_3$ in
 order that $P_2^{(1)}+P_3^{(1)}$ meets the o-path criteria at $r$. It then follows from our construction of
 $P_2^{(1)}$ and $P_3^{(1)}$ that $P_2^{(1)}+P_3^{(1)}$ is an o-path from $z$ to $v$. Moreover, 
 $P_1^{(1)}+P_2^{(1)}+P_3^{(1)}$ meets the o-path criteria at $z$ since $z$ was a non-transversely oriented cut-vertex
 and $P_1^{(1)}$ is coming out of $U_1$ while $P_2^{(1)}$ is entering $U_2$. Finally, $P_2^{(1)}$ lies on the 
 arc of $C_2$ that did not meet $P_1$, while the extension of $P_1$ (except for the edge to $z$) was contained 
 within $U_1$ and $P_2^{(1)}$ is
 contained within $U_2$, so $P_1^{(1)}$ does not meet $P_2^{(1)}$ other than at $z$. Thus $P_1^{(1)}+P_2^{(1)}+P_3^{(1)}$
 is an o-cycle with $P_1^{(1)}$ lying on $C_1$, $P_2^{(1)}$ lying on $C_2$, and $P_3^{(1)}$ lying on $C_3$, so
 $(P_1^{(1)},P_2^{(1)},P_3^{(1)})\in \it O$ and $P_1<P_1^{(1)}$, so $(P_1,P_2,P_3)<(P_1^{(1)},P_2^{(1)},P_3^{(1)})$.
 We claim that $(P_1^{(1)},P_2^{(1)},P_3^{(1)})$ satisfies Condition $A$. Let $I$ denote the arc
 of $C_2$ from $u$ to $z$ which passes through $x$. Now $P_1$ only meets $C_2$ at vertices on the arc of $C_2$ between $u$
 and $x$ that does not contain $z$, which is a subpath of $I$, so $P_1$ only meets $C_2$ at vertices of $I$.
 As well, the extension of $P_1$ can only meet $C_2$ at vertices of $U_1$ or $z$, hence only at vertices of $I$. Thus
 $P_1^{(1)}$ only meets $C_2$ at vertices of $I$. It remains to prove
 that $C_1-P_1^{(1)}$ does not meet $C_2$ at vertices of $I$. As the vertices of $C_1-P_1^{(1)}$ form a subset of
 the set of vertices of $C_1-P_1$, and $C_1-P_1$ could only meet $C_2$ on the arc of $C_2$ from $x$ to $u$ that 
 passes through $z$, it follows that $C_1-P_1^{(1)}$ can only meet $C_2$ at vertices on the arc of $C_2$
 from $x$ to $u$ that passes through $z$. However, $C_1-P_1^{(1)}$ lies in $U_2$ and thus does not pass through
 any vertex of the arc of $C_2$ from $x$ to $z$ that lies in $U_1$, so it follows that $C_1-P_1^{(1)}$ does not
 meet $C_2$ at any vertex of $I$. Thus $(P_1^{(1)},P_2^{(1)},P_3^{(1)})$ satisfies Condition A.
 
 We have now established that for every $\it O_i\in \it O$ that satisfies Condition A and is such that
 $G-E(\it \hat{O}_i)$ has a non-transversely oriented cut-vertex belonging to $G_1$, there is a larger element of $\it O$ that also
 satisfies Condition A. Since $\it O_1$ satisfies Condition A and $\it O$ is finite, it follows that there is a greatest element $\it O_i=(P_1,P_2,P_3)$
 of $\it O$ that can be produced by applying this construction to an element of $\it O$ that satisfies Condition
 A. Thus for any $\it O_j\in \it O$ that satisfies Condition $A$ and is greater than or equal to $\it O_i$, 
 $G-E(\it\hat{O}_j)$ does not have a non-transversely oriented cut-vertex in $G_1$. Since $\it O_i$ was constructed
 by an application of the procedure described above, we know that the terminal vertex of $P_2$, $s$ say, is the only
 vertex on $P_2$ that lies on $C_3$. We shall say that an element $(Q_1,Q_2,Q_3)\in \it O$ satisfies Condition
 B if $Q_1=P_1$ and $Q_2$ meets $C_3$ only at vertices of one of the arcs of $C_3$ determined by $s$ and the
 terminal vertex of $P_2$, while $C_2-P_2$ does not meet $C_3$ at any vertex of this arc. In particular, 
 $(P_1,P_2,P_3)$ satisfies Condition B as well as Condition A.
 By assumption, every element $\it O_j$ of $\it O$ is such that $G-E(\it\hat O_j)$ contains a non-transversely
 oriented cut-vertex, and thus in particular, $G-E(\it\hat O_i)$ has a non-transversely oriented cut-vertex,
 necessarily in $G_2$.
 For any $\it O_j=(Q_1,Q_2,Q_3)\in \it O$ that satisfies Condition B (so $Q_1=P_1$ and thus it must satisfy Condition A)
 for which $\it O_j\ge \it O_i$ and $G-E(\it\hat O_j)$ has a non-transversely oriented cut-vertex (necessarily in $G_2$),
 we may carry out a procedure completely analagous to that described above for $C_1$ and
 $C_2$ to obtain an element $(Q_1^{(1)},Q_2^{(1)},Q_3^{(1)})$ that satisfies condition B (and thus A),
 and which is greater than $(Q_1,Q_2,Q_3)$. Again, since $\it O$ is finite, there is a maximum such element of
 $\it O$, which we shall denote by $M$. Thus $M$ satisfies both Condtions A and B, and $G-E(\hat M)$ can
 not contain a non-transversely oriented cut-vertex in either $G_1$ or $G_2$, which contradicts our assumption
 that every element $\it O_k$ of $\it O$ was such that $G-E(\it\hat O_k)$ contains a non-transversely oriented
 cut-vertex.

 This  completes the proof of the inductive step, and so the result follows.
 \endproof 

Of course, the goal is to obtain an alternative proof of the four colour
theorem. This would be accomplished if we could sharpen our theorem above
to say that every vertex-oriented 4-regular planar graph with all
cut-vertices oriented transversely can be 3-o-coloured. There are three places
in our proof of the inductive step where the number of colours used to colour
$(G,\sigma)$ may be increased over the number used to colour the smaller graph.
Two of these situations involve the removal of a cycle, o-colouring 
the result, and finally reinserting the cycle, possibly needing an additional 
colour for it, while the other appears in a simplification step during the
proof of Case 3, where after smoothing a vertex and o-colouring the resulting
graph, if the two new edges that resulted from the smoothing were coloured the
same but belonged to different o-cycles, then we observed that one of the o-cycles
could have its colour changed, possibly requiring a new colour. It might in fact
be possible to argue that the case itself never happens, in the sense that it may be impossible
that each and every vertex of $G$ can result in the scenario of Case 3. If that can be established, the remaining problem occurs in
Case 2 (iii), and is potentially the more intractable one. We offer an example
below of the situation that may occur. 

\medskip
 \begin{figure}[h]
  \centering\table{c@{\hskip10pt}c@{\hskip10pt}c}
   $\vcenter{
\xy /r.2pt/:,(0,0);(1,0):(0,.65)::,
(75.5228,154.173)="1";
(116.084,171.016)="flex18"; 
(119.936,250.534)="flex19"; 
(120.93,289.768)="3";
(116.522,329.991)="flex20"; 
(5.10846,250.711)="flex21"; 
(88.9067,194.31)="flex22"; 
(120.71,211.271)="2";
(151.901,227.879)="flex23"; 
(199.713,271.97)="flex24"; 
(225.692,282.866)="6";
(264.805,283.61)="flex25"; 
(330.95,321.421)="flex26"; 
(338.496,367.79)="8";
(318.885,409.793)="flex27"; 
(153.17,451.891)="flex28"; 
(76.0992,347.012)="4";
(89.2568,306.86)="flex29"; 
(152.001,273.021)="flex30"; 
(176.354,250.415)="5";
(199.636,228.786)="flex31"; 
(264.731,216.865)="flex32"; 
(300.752,209.342)="11";
(330.683,178.682)="flex33"; 
(318.059,90.2756)="flex34"; 
(279.213,52.0458)="13";
(151.924,48.9551)="flex35"; 
(256.33,98.0217)="flex36"; 
(260.53,143.002)="14";
(281.013,176.771)="flex37"; 
(312.455,250.121)="flex38"; 
(300.85,290.936)="7";
(281.268,323.608)="flex39"; 
(257.148,402.418)="flex40"; 
(280.319,448.173)="9";
(394.1,450.64)="flex41"; 
(415.699,249.837)="flex42"; 
(432.676,162.643)="17";
(393.072,48.9368)="flex43"; 
(384.486,144.128)="flex44"; 
(337.948,132.242)="12";
(297.861,132.409)="flex45"; 
(230.283,173.92)="flex46"; 
(225.586,217.8)="10";
(225.063,250.34)="flex47"; 
(230.589,326.721)="flex48"; 
(260.99,357.459)="15";
(298.383,367.841)="flex49"; 
(384.99,355.672)="flex50"; 
(433.097,336.894)="16";
(494.916,249.628)="flex51"; 
"1";"flex18"**\crv{ (91.4283,152.433) & (107.83,157.582) },
"flex18";"2"*[o]=(0,0){\,}**\crv{ (123.364,182.864) & (121.636,197.43) },
"2"*[o]=(0,0){\,};"flex19"**\crv{ (119.835,224.337) & (119.899,237.44) },
"flex19";"3"**\crv{ (119.973,263.62) & (119.984,276.715) },
"3";"flex20"**\crv{ (121.932,303.585) & (123.734,318.123) },
"flex20";"4"*[o]=(0,0){\,}**\crv{ (108.345,343.447) & (91.9916,348.664) },
"4"*[o]=(0,0){\,};"flex21"**\crv{ (31.7432,342.402) & (5.21707,297.779) },
"flex21";"1"**\crv{ (5,203.711) & (31.2529,159.017) },
"1"*[o]=(0,0){\,};"flex22"**\crv{ (75.1843,168.89) & (78.5694,183.887) },
"flex22";"2"**\crv{ (97.489,202.963) & (109.426,206.768) },
"2";"flex23"**\crv{ (131.691,215.654) & (142.348,220.903) },
"flex23";"5"*[o]=(0,0){\,}**\crv{ (160.874,234.432) & (168.69,242.378) },
"5"*[o]=(0,0){\,};"flex24"**\crv{ (183.685,258.103) & (190.956,265.942) },
"flex24";"6"**\crv{ (207.534,277.353) & (216.351,281.152) },
"6";"flex25"**\crv{ (238.57,285.229) & (251.732,283.512) },
"flex25";"7"*[o]=(0,0){\,}**\crv{ (277.218,283.703) & (289.693,285.47) },
"7"*[o]=(0,0){\,};"flex26"**\crv{ (314.008,297.383) & (324.238,308.397) },
"flex26";"8"**\crv{ (338.316,335.713) & (341.242,351.949) },
"8";"flex27"**\crv{ (335.82,383.223) & (327.981,397.063) },
"flex27";"9"*[o]=(0,0){\,}**\crv{ (308.245,424.682) & (295.754,438.32) },
"9"*[o]=(0,0){\,};"flex28"**\crv{ (242.33,472.421) & (194.276,470.992) },
"flex28";"4"**\crv{ (110.378,432.007) & (77.4163,393.639) },
"4";"flex29"**\crv{ (75.6841,332.316) & (78.9845,317.32) },
"flex29";"3"*[o]=(0,0){\,}**\crv{ (97.7839,298.177) & (109.686,294.323) },
"3"*[o]=(0,0){\,};"flex30"**\crv{ (131.876,285.333) & (142.49,280.03) },
"flex30";"5"**\crv{ (160.938,266.434) & (168.723,258.472) },
"5";"flex31"**\crv{ (183.658,242.705) & (190.9,234.842) },
"flex31";"10"*[o]=(0,0){\,}**\crv{ (207.442,223.376) & (216.25,219.551) },
"10"*[o]=(0,0){\,};"flex32"**\crv{ (238.466,215.384) & (251.647,217.033) },
"flex32";"11"**\crv{ (277.149,216.706) & (289.62,214.872) },
"11";"flex33"**\crv{ (313.877,202.823) & (324.047,191.749) },
"flex33";"12"*[o]=(0,0){\,}**\crv{ (337.967,164.339) & (340.789,148.076) },
"12"*[o]=(0,0){\,};"flex34"**\crv{ (335.176,116.797) & (327.244,102.977) },
"flex34";"13"**\crv{ (307.317,75.4227) & (294.727,61.8322) },
"13";"flex35"**\crv{ (241.056,27.9758) & (192.963,29.6086) },
"flex35";"1"*[o]=(0,0){\,}**\crv{ (109.286,69.0555) & (76.5946,107.584) },
"13"*[o]=(0,0){\,};"flex36"**\crv{ (266.709,64.3593) & (259.212,80.7113) },
"flex36";"14"**\crv{ (253.81,113.153) & (254.86,128.746) },
"14";"flex37"**\crv{ (265.418,155.291) & (273.498,165.909) },
"flex37";"11"*[o]=(0,0){\,}**\crv{ (288.243,187.223) & (295.003,198.006) },
"11"*[o]=(0,0){\,};"flex38"**\crv{ (307.236,222.128) & (312.442,235.769) },
"flex38";"7"**\crv{ (312.468,264.475) & (307.285,278.125) },
"7";"flex39"**\crv{ (295.144,302.296) & (288.446,313.118) },
"flex39";"15"*[o]=(0,0){\,}**\crv{ (273.817,334.5) & (265.791,345.148) },
"15"*[o]=(0,0){\,};"flex40"**\crv{ (255.422,371.738) & (254.512,387.324) },
"flex40";"9"**\crv{ (260.164,419.679) & (267.758,435.955) },
"9";"flex41"**\crv{ (311.117,478.128) & (359.816,477.674) },
"flex41";"16"*[o]=(0,0){\,}**\crv{ (427.737,424.117) & (439.589,379.605) },
"16"*[o]=(0,0){\,};"flex42"**\crv{ (428.634,307.525) & (415.777,279.614) },
"flex42";"17"**\crv{ (415.622,220.042) & (428.349,192.05) },
"17";"flex43"**\crv{ (438.975,119.841) & (426.894,75.3203) },
"flex43";"13"*[o]=(0,0){\,}**\crv{ (358.611,22.0551) & (309.855,21.8718) },
"17"*[o]=(0,0){\,};"flex44"**\crv{ (417.338,154.763) & (400.894,149.406) },
"flex44";"12"**\crv{ (369.215,139.216) & (353.862,134.337) },
"12";"flex45"**\crv{ (324.644,130.491) & (311.171,130.717) },
"flex45";"14"*[o]=(0,0){\,}**\crv{ (284.937,134.053) & (272.151,137.092) },
"14"*[o]=(0,0){\,};"flex46"**\crv{ (247.256,149.754) & (235.915,160.143) },
"flex46";"10"**\crv{ (224.646,187.712) & (225.583,202.944) },
"10";"flex47"**\crv{ (225.587,228.649) & (225.046,239.491) },
"flex47";"6"*[o]=(0,0){\,}**\crv{ (225.079,261.185) & (225.652,272.021) },
"6"*[o]=(0,0){\,};"flex48"**\crv{ (225.745,297.722) & (224.88,312.959) },
"flex48";"15"**\crv{ (236.29,340.462) & (247.684,350.785) },
"15";"flex49"**\crv{ (272.644,363.305) & (285.448,366.27) },
"flex49";"8"*[o]=(0,0){\,}**\crv{ (311.71,369.459) & (325.192,369.611) },
"8"*[o]=(0,0){\,};"flex50"**\crv{ (354.404,365.613) & (369.739,360.663) },
"flex50";"16"**\crv{ (401.377,350.309) & (417.796,344.858) },
"16";"flex51"**\crv{ (467.498,318.989) & (495,287.832) },
"flex51";"17"*[o]=(0,0){\,}**\crv{ (494.832,211.414) & (467.181,180.371) },
"5"*!<0pt,7pt>{\hbox{$v$}},
"5";"5"+(40,0)**\dir{-}*\dir{>},
"5";"5"-(40,0)**\dir{-}*\dir{>},
"1"*!<0pt,0pt>{\smallbullet},
"2"*!<0pt,0pt>{\smallbullet},
"3"*!<0pt,0pt>{\smallbullet},
"4"*!<0pt,0pt>{\smallbullet},
"5"*!<0pt,0pt>{\smallbullet},
"6"*!<0pt,0pt>{\smallbullet},
"7"*!<0pt,0pt>{\smallbullet},
"8"*!<0pt,0pt>{\smallbullet},
"9"*!<0pt,0pt>{\smallbullet},
"10"*!<0pt,0pt>{\smallbullet},
"11"*!<0pt,0pt>{\smallbullet},
"12"*!<0pt,0pt>{\smallbullet},
"13"*!<0pt,0pt>{\smallbullet},
"14"*!<0pt,0pt>{\smallbullet},
"15"*!<0pt,0pt>{\smallbullet},
"16"*!<0pt,0pt>{\smallbullet},
"17"*!<0pt,0pt>{\smallbullet},
\endxy}$
 &
 $\vcenter{
\xy /r.2pt/:,(0,0);(1,0):(0,.65)::,
(75.5228,154.173)="1";
(116.084,171.016)="flex18"; 
(119.936,250.534)="flex19"; 
(120.93,289.768)="3";
(116.522,329.991)="flex20"; 
(5.10846,250.711)="flex21"; 
(88.9067,194.31)="flex22"; 
(120.71,211.271)="2";
(151.901,227.879)="flex23"; 
(199.713,271.97)="flex24"; 
(225.692,282.866)="6";
(264.805,283.61)="flex25"; 
(330.95,321.421)="flex26"; 
(338.496,367.79)="8";
(318.885,409.793)="flex27"; 
(153.17,451.891)="flex28"; 
(76.0992,347.012)="4";
(89.2568,306.86)="flex29"; 
(152.001,273.021)="flex30"; 
(176.354,250.415)="5";
(199.636,228.786)="flex31"; 
(264.731,216.865)="flex32"; 
(300.752,209.342)="11";
(330.683,178.682)="flex33"; 
(318.059,90.2756)="flex34"; 
(279.213,52.0458)="13";
(151.924,48.9551)="flex35"; 
(256.33,98.0217)="flex36"; 
(260.53,143.002)="14";
(281.013,176.771)="flex37"; 
(312.455,250.121)="flex38"; 
(300.85,290.936)="7";
(281.268,323.608)="flex39"; 
(257.148,402.418)="flex40"; 
(280.319,448.173)="9";
(394.1,450.64)="flex41"; 
(415.699,249.837)="flex42"; 
(432.676,162.643)="17";
(393.072,48.9368)="flex43"; 
(384.486,144.128)="flex44"; 
(337.948,132.242)="12";
(297.861,132.409)="flex45"; 
(230.283,173.92)="flex46"; 
(225.586,217.8)="10";
(225.063,250.34)="flex47"; 
(230.589,326.721)="flex48"; 
(260.99,357.459)="15";
(298.383,367.841)="flex49"; 
(384.99,355.672)="flex50"; 
(433.097,336.894)="16";
(494.916,249.628)="flex51"; 
"1";"flex18"**\crv{ (91.4283,152.433) & (107.83,157.582) },
"flex18";"2"*[o]=(0,0){\,}**\crv{ (123.364,182.864) & (121.636,197.43) },
"2"*[o]=(0,0){\,};"flex19"**\crv{ (119.835,224.337) & (119.899,237.44) },
"flex19";"3"**\crv{ (119.973,263.62) & (119.984,276.715) },
"3";"flex20"**\crv{ (121.932,303.585) & (123.734,318.123) },
"flex20";"4"*[o]=(0,0){\,}**\crv{ (108.345,343.447) & (91.9916,348.664) },
"4"*[o]=(0,0){\,};"flex21"**\crv{ (31.7432,342.402) & (5.21707,297.779) },
"flex21";"1"**\crv{ (5,203.711) & (31.2529,159.017) },
"1"*[o]=(0,0){\,};"flex22"**\crv{ (75.1843,168.89) & (78.5694,183.887) },
"flex22";"2"**\crv{ (97.489,202.963) & (109.426,206.768) },
"2";"flex23"**\crv{ (131.691,215.654) & (142.348,220.903) },
"flex23";
"flex23";"flex30"**\dir{}?(.5)="x",
"flex30"**\crv{"flex23"+<2pt,2pt> & "x"+<2pt,0pt> & "flex30"+<2pt,-2pt> },
"flex24";"6"**\crv{ (207.534,277.353) & (216.351,281.152) },
"6";"flex25"**\crv{ (238.57,285.229) & (251.732,283.512) },
"flex25";"7"*[o]=(0,0){\,}**\crv{ (277.218,283.703) & (289.693,285.47) },
"7"*[o]=(0,0){\,};"flex26"**\crv{ (314.008,297.383) & (324.238,308.397) },
"flex26";"8"**\crv{ (338.316,335.713) & (341.242,351.949) },
"8";"flex27"**\crv{ (335.82,383.223) & (327.981,397.063) },
"flex27";"9"*[o]=(0,0){\,}**\crv{ (308.245,424.682) & (295.754,438.32) },
"9"*[o]=(0,0){\,};"flex28"**\crv{ (242.33,472.421) & (194.276,470.992) },
"flex28";"4"**\crv{ (110.378,432.007) & (77.4163,393.639) },
"4";"flex29"**\crv{ (75.6841,332.316) & (78.9845,317.32) },
"flex29";"3"*[o]=(0,0){\,}**\crv{ (97.7839,298.177) & (109.686,294.323) },
"3"*[o]=(0,0){\,};"flex30"**\crv{ (131.876,285.333) & (142.49,280.03) },
"flex24";"flex31"**\dir{}?(.5)="y",
"flex24";"flex31"**\crv{"flex24"+<-2pt,-2pt> & "y"+<-2pt,0pt> & "flex31"+<-2pt,2pt> },
"flex31";"10"*[o]=(0,0){\,}**\crv{ (207.442,223.376) & (216.25,219.551) },
"10"*[o]=(0,0){\,};"flex32"**\crv{ (238.466,215.384) & (251.647,217.033) },
"flex32";"11"**\crv{ (277.149,216.706) & (289.62,214.872) },
"11";"flex33"**\crv{ (313.877,202.823) & (324.047,191.749) },
"flex33";"12"*[o]=(0,0){\,}**\crv{ (337.967,164.339) & (340.789,148.076) },
"12"*[o]=(0,0){\,};"flex34"**\crv{ (335.176,116.797) & (327.244,102.977) },
"flex34";"13"**\crv{ (307.317,75.4227) & (294.727,61.8322) },
"13";"flex35"**\crv{ (241.056,27.9758) & (192.963,29.6086) },
"flex35";"1"*[o]=(0,0){\,}**\crv{ (109.286,69.0555) & (76.5946,107.584) },
"13"*[o]=(0,0){\,};"flex36"**\crv{ (266.709,64.3593) & (259.212,80.7113) },
"flex36";"14"**\crv{ (253.81,113.153) & (254.86,128.746) },
"14";"flex37"**\crv{ (265.418,155.291) & (273.498,165.909) },
"flex37";"11"*[o]=(0,0){\,}**\crv{ (288.243,187.223) & (295.003,198.006) },
"11"*[o]=(0,0){\,};"flex38"**\crv{ (307.236,222.128) & (312.442,235.769) },
"flex38";"7"**\crv{ (312.468,264.475) & (307.285,278.125) },
"7";"flex39"**\crv{ (295.144,302.296) & (288.446,313.118) },
"flex39";"15"*[o]=(0,0){\,}**\crv{ (273.817,334.5) & (265.791,345.148) },
"15"*[o]=(0,0){\,};"flex40"**\crv{ (255.422,371.738) & (254.512,387.324) },
"flex40";"9"**\crv{ (260.164,419.679) & (267.758,435.955) },
"9";"flex41"**\crv{ (311.117,478.128) & (359.816,477.674) },
"flex41";"16"*[o]=(0,0){\,}**\crv{ (427.737,424.117) & (439.589,379.605) },
"16"*[o]=(0,0){\,};"flex42"**\crv{ (428.634,307.525) & (415.777,279.614) },
"flex42";"17"**\crv{ (415.622,220.042) & (428.349,192.05) },
"17";"flex43"**\crv{ (438.975,119.841) & (426.894,75.3203) },
"flex43";"13"*[o]=(0,0){\,}**\crv{ (358.611,22.0551) & (309.855,21.8718) },
"17"*[o]=(0,0){\,};"flex44"**\crv{ (417.338,154.763) & (400.894,149.406) },
"flex44";"12"**\crv{ (369.215,139.216) & (353.862,134.337) },
"12";"flex45"**\crv{ (324.644,130.491) & (311.171,130.717) },
"flex45";"14"*[o]=(0,0){\,}**\crv{ (284.937,134.053) & (272.151,137.092) },
"14"*[o]=(0,0){\,};"flex46"**\crv{ (247.256,149.754) & (235.915,160.143) },
"flex46";"10"**\crv{ (224.646,187.712) & (225.583,202.944) },
"10";"flex47"**\crv{ (225.587,228.649) & (225.046,239.491) },
"flex47";"6"*[o]=(0,0){\,}**\crv{ (225.079,261.185) & (225.652,272.021) },
"6"*[o]=(0,0){\,};"flex48"**\crv{ (225.745,297.722) & (224.88,312.959) },
"flex48";"15"**\crv{ (236.29,340.462) & (247.684,350.785) },
"15";"flex49"**\crv{ (272.644,363.305) & (285.448,366.27) },
"flex49";"8"*[o]=(0,0){\,}**\crv{ (311.71,369.459) & (325.192,369.611) },
"8"*[o]=(0,0){\,};"flex50"**\crv{ (354.404,365.613) & (369.739,360.663) },
"flex50";"16"**\crv{ (401.377,350.309) & (417.796,344.858) },
"16";"flex51"**\crv{ (467.498,318.989) & (495,287.832) },
"flex51";"17"*[o]=(0,0){\,}**\crv{ (494.832,211.414) & (467.181,180.371) },
"1"*!<0pt,0pt>{\smallbullet},
"2"*!<0pt,0pt>{\smallbullet},
"3"*!<0pt,0pt>{\smallbullet},
"4"*!<0pt,0pt>{\smallbullet},
"6"*!<0pt,0pt>{\smallbullet},
"7"*!<0pt,0pt>{\smallbullet},
"8"*!<0pt,0pt>{\smallbullet},
"9"*!<0pt,0pt>{\smallbullet},
"10"*!<0pt,0pt>{\smallbullet},
"11"*!<0pt,0pt>{\smallbullet},
"12"*!<0pt,0pt>{\smallbullet},
"13"*!<0pt,0pt>{\smallbullet},
"14"*!<0pt,0pt>{\smallbullet},
"15"*!<0pt,0pt>{\smallbullet},
"16"*!<0pt,0pt>{\smallbullet},
"17"*!<0pt,0pt>{\smallbullet},
\endxy}$
&
 $\vcenter{
\xy /r.2pt/:,(0,0);(1,0):(0,.65)::,
(75.5228,154.173)="1";
(116.084,171.016)="flex18"; 
(119.936,250.534)="flex19"; 
(120.93,289.768)="3";
(116.522,329.991)="flex20"; 
(5.10846,250.711)="flex21"; 
(88.9067,194.31)="flex22"; 
(120.71,211.271)="2";
(151.901,227.879)="flex23"; 
(199.713,271.97)="flex24"; 
(225.692,282.866)="6";
(264.805,283.61)="flex25"; 
(330.95,321.421)="flex26"; 
(338.496,367.79)="8";
(318.885,409.793)="flex27"; 
(153.17,451.891)="flex28"; 
(76.0992,347.012)="4";
(89.2568,306.86)="flex29"; 
(152.001,273.021)="flex30"; 
(176.354,250.415)="5";
(199.636,228.786)="flex31"; 
(264.731,216.865)="flex32"; 
(300.752,209.342)="11";
(330.683,178.682)="flex33"; 
(318.059,90.2756)="flex34"; 
(279.213,52.0458)="13";
(151.924,48.9551)="flex35"; 
(256.33,98.0217)="flex36"; 
(260.53,143.002)="14";
(281.013,176.771)="flex37"; 
(312.455,250.121)="flex38"; 
(300.85,290.936)="7";
(281.268,323.608)="flex39"; 
(257.148,402.418)="flex40"; 
(280.319,448.173)="9";
(394.1,450.64)="flex41"; 
(415.699,249.837)="flex42"; 
(432.676,162.643)="17";
(393.072,48.9368)="flex43"; 
(384.486,144.128)="flex44"; 
(337.948,132.242)="12";
(297.861,132.409)="flex45"; 
(230.283,173.92)="flex46"; 
(225.586,217.8)="10";
(225.063,250.34)="flex47"; 
(230.589,326.721)="flex48"; 
(260.99,357.459)="15";
(298.383,367.841)="flex49"; 
(384.99,355.672)="flex50"; 
(433.097,336.894)="16";
(494.916,249.628)="flex51"; 
"1";"flex18"**\crv{ (91.4283,152.433) & (107.83,157.582) },
"flex18";"2"*[o]=(0,0){\,}**\crv{ (123.364,182.864) & (121.636,197.43) },
"2"*[o]=(0,0){\,};"flex19"**\crv{ (119.835,224.337) & (119.899,237.44) },
"flex19";"3"**\crv{ (119.973,263.62) & (119.984,276.715) },
"3";"flex20"**\crv{ (121.932,303.585) & (123.734,318.123) },
"flex20";"4"*[o]=(0,0){\,}**\crv{ (108.345,343.447) & (91.9916,348.664) },
"4"*[o]=(0,0){\,};"flex21"**\crv{ (31.7432,342.402) & (5.21707,297.779) },
"flex21";"1"**\crv{ (5,203.711) & (31.2529,159.017) },
"1"*[o]=(0,0){\,};"flex22"**\crv{ (75.1843,168.89) & (78.5694,183.887) },
"flex22";"2"**\crv{ (97.489,202.963) & (109.426,206.768) },
"2";"flex23"**\crv{ (131.691,215.654) & (142.348,220.903) },
"flex23";"5"*[o]=(0,0){\,}**\crv{ (160.874,234.432) & (168.69,242.378) },
"5"*[o]=(0,0){\,};"flex24"**\crv{ (183.685,258.103) & (190.956,265.942) },
"flex24";"6"**\crv{ (207.534,277.353) & (216.351,281.152) },
"6";"flex25"**\crv{ (238.57,285.229) & (251.732,283.512) },
"flex25";"7"*[o]=(0,0){\,}**\crv{ (277.218,283.703) & (289.693,285.47) },
"7"*[o]=(0,0){\,};"flex26"**\crv{ (314.008,297.383) & (324.238,308.397) },
"flex26";"8"**\crv{ (338.316,335.713) & (341.242,351.949) },
"8";"flex27"**\crv{ (335.82,383.223) & (327.981,397.063) },
"flex27";"9"*[o]=(0,0){\,}**\crv{ (308.245,424.682) & (295.754,438.32) },
"9"*[o]=(0,0){\,};"flex28"**\crv{ (242.33,472.421) & (194.276,470.992) },
"flex28";"4"**\crv{ (110.378,432.007) & (77.4163,393.639) },
"4";"flex29"**\crv{ (75.6841,332.316) & (78.9845,317.32) },
"flex29";"3"*[o]=(0,0){\,}**\crv{ (97.7839,298.177) & (109.686,294.323) },
"3"*[o]=(0,0){\,};"flex30"**\crv{ (131.876,285.333) & (142.49,280.03) },
"flex30";"5"**\crv{ (160.938,266.434) & (168.723,258.472) },
"5";"flex31"**\crv{ (183.658,242.705) & (190.9,234.842) },
"flex31";"10"*[o]=(0,0){\,}**\crv{ (207.442,223.376) & (216.25,219.551) },
"10"*[o]=(0,0){\,};"flex32"**\crv{ (238.466,215.384) & (251.647,217.033) },
"flex32";"11"**\crv{ (277.149,216.706) & (289.62,214.872) },
"11";"flex33"**\crv{ (313.877,202.823) & (324.047,191.749) },
"flex33";"12"*[o]=(0,0){\,}**\crv{ (337.967,164.339) & (340.789,148.076) },
"12"*[o]=(0,0){\,};"flex34"**\crv{ (335.176,116.797) & (327.244,102.977) },
"flex34";"13"**\crv{ (307.317,75.4227) & (294.727,61.8322) },
"13";"flex35"**\crv{ (241.056,27.9758) & (192.963,29.6086) },
"flex35";"1"*[o]=(0,0){\,}**\crv{ (109.286,69.0555) & (76.5946,107.584) },
"13"*[o]=(0,0){\,};"flex36"**\crv{ (266.709,64.3593) & (259.212,80.7113) },
"flex36";"14"**\crv{ (253.81,113.153) & (254.86,128.746) },
"14";"flex37"**\crv{ (265.418,155.291) & (273.498,165.909) },
"flex37";"11"*[o]=(0,0){\,}**\crv{ (288.243,187.223) & (295.003,198.006) },
"11"*[o]=(0,0){\,};"flex38"**\crv{ (307.236,222.128) & (312.442,235.769) },
"flex38";"7"**\crv{ (312.468,264.475) & (307.285,278.125) },
"7";"flex39"**\crv{ (295.144,302.296) & (288.446,313.118) },
"flex39";"15"*[o]=(0,0){\,}**\crv{ (273.817,334.5) & (265.791,345.148) },
"15"*[o]=(0,0){\,};"flex40"**\crv{ (255.422,371.738) & (254.512,387.324) },
"flex40";"9"**\crv{ (260.164,419.679) & (267.758,435.955) },
"9";"flex41"**\crv{ (311.117,478.128) & (359.816,477.674) },
"flex41";"16"*[o]=(0,0){\,}**\crv{ (427.737,424.117) & (439.589,379.605) },
"16"*[o]=(0,0){\,};"flex42"**\crv{ (428.634,307.525) & (415.777,279.614) },
"flex42";"17"**\crv{ (415.622,220.042) & (428.349,192.05) },
"17";"flex43"**\crv{ (438.975,119.841) & (426.894,75.3203) },
"flex43";"13"*[o]=(0,0){\,}**\crv{ (358.611,22.0551) & (309.855,21.8718) },
"17"*[o]=(0,0){\,};"flex44"**\crv{~*=<2.5pt>{\hbox{\LARGE.}} (417.338,154.763) & (400.894,149.406) },
"flex44";"12"**\crv{~*=<2.5pt>{\hbox{\LARGE.}} (369.215,139.216) & (353.862,134.337) },
"12";"flex45"**\crv{~*=<2.5pt>{\hbox{\LARGE.}} (324.644,130.491) & (311.171,130.717) },
"flex45";"14"*[o]=(0,0){\,}**\crv{~*=<2.5pt>{\hbox{\LARGE.}} (284.937,134.053) & (272.151,137.092) },
"14"*[o]=(0,0){\,};"flex46"**\crv{~*=<2.5pt>{\hbox{\LARGE.}} (247.256,149.754) & (235.915,160.143) },
"flex46";"10"**\crv{~*=<2.5pt>{\hbox{\LARGE.}} (224.646,187.712) & (225.583,202.944) },
"10";"flex47"**\crv{~*=<2.5pt>{\hbox{\LARGE.}} (225.587,228.649) & (225.046,239.491) },
"flex47";"6"*[o]=(0,0){\,}**\crv{~*=<2.5pt>{\hbox{\LARGE.}} (225.079,261.185) & (225.652,272.021) },
"6"*[o]=(0,0){\,};"flex48"**\crv{~*=<2.5pt>{\hbox{\LARGE.}} (225.745,297.722) & (224.88,312.959) },
"flex48";"15"**\crv{~*=<2.5pt>{\hbox{\LARGE.}} (236.29,340.462) & (247.684,350.785) },
"15";"flex49"**\crv{~*=<2.5pt>{\hbox{\LARGE.}} (272.644,363.305) & (285.448,366.27) },
"flex49";"8"*[o]=(0,0){\,}**\crv{~*=<2.5pt>{\hbox{\LARGE.}} (311.71,369.459) & (325.192,369.611) },
"8"*[o]=(0,0){\,};"flex50"**\crv{~*=<2.5pt>{\hbox{\LARGE.}} (354.404,365.613) & (369.739,360.663) },
"flex50";"16"**\crv{~*=<2.5pt>{\hbox{\LARGE.}} (401.377,350.309) & (417.796,344.858) },
"16";"flex51"**\crv{~*=<2.5pt>{\hbox{\LARGE.}} (467.498,318.989) & (495,287.832) },
"flex51";"17"*[o]=(0,0){\,}**\crv{~*=<2.5pt>{\hbox{\LARGE.}} (494.832,211.414) & (467.181,180.371) },
"5"*!<0pt,7pt>{\hbox{$v$}},
"5";"5"+(40,0)**\dir{-}*\dir{>},
"5";"5"-(40,0)**\dir{-}*\dir{>},
"1"*!<0pt,0pt>{\smallbullet},
"2"*!<0pt,0pt>{\smallbullet},
"3"*!<0pt,0pt>{\smallbullet},
"4"*!<0pt,0pt>{\smallbullet},
"5"*!<0pt,0pt>{\smallbullet},
"7"*!<0pt,0pt>{\smallbullet},
"9"*!<0pt,0pt>{\smallbullet},
"11"*!<0pt,0pt>{\smallbullet},
"13"*!<0pt,0pt>{\smallbullet},
\endxy}$\\
  \noalign{\vskip4pt}
  (a) & (b) & (c)
  \endtable\caption{}\label{figure: our difficulties figure}
 \end{figure}

\noindent  In Figure \ref{figure: our difficulties figure} (a), we have not shown the orientations of any vertex
other than $v$, but it is intended that the graph in (b) (in which vertex $v$ has been smoothed) 
has been o-coloured in such a way that the four simple smooth curves are o-cycles. As in the proof of 
Case 2 (iii), we choose an o-cycle to remove, and our choice is the curve $C$ shown dotted in (c). Now
o-colour $G-E(C)$. No matter what orientations had been assigned to the vertices of $G$ (other than $v$, which
is to be oriented as shown), $C$ will meet o-cycles of $G-E(C)$ of three different colours, and so the edges of
$C$ must be assigned a new colour.

We conclude this section with a brief discussion of vertex-orientation for
arbitrary 4-regular graphs. By an orientation of a vertex $v$, we mean a partition of the four incident edges
into two cells of size 2 (where we treat each loop at $v$ as two incident edges). Then define o-colouring
of a vertex-oriented 4-regular graph just as was done for planar vertex-oriented 4-regular graphs. If a
vertex-oriented 4-regular graph can be o-coloured, then its edge set can be decomposed into a collection of
edge-disjoint cycles (each an o-cycle of the vertex-oriented graph). 

Note that given a 3-regular graph and a 1-factor of the graph, we may obtain a vertex-oriented 4-regular graph
by collapsing the edges of the 1-factor (with the orientation of each vertex determined by the 1-factor edge
that gave rise to the vertex, just as was done in the proof of Theorem \ref{theorem: 4ct equivalence}).
In an initial examination of snarks, we observed that many snarks had the property that there was at least one 
1-factor of the snark that gave rise to a non-o-colourable vertex-oriented 4-regular graph, and frequently,
this was true for every 1-factor of the snark. This appears to be an interesting avenue of exploration.

\section{Examples}

In many of the early examples of vertex-oriented 4-regular planar graphs that we had examined, it was noticed that 
there was at least one 3-o-colouring in which there is one colour and
exactly one o-cycle component of the subgraph induced by the edges of
that colour, and that o-cycle meets all other o-cycles determined by the
3-o-colouring. Often, this o-cycle has maximum
length over all o-cycles determined by the o-colouring. Our first example 
to demonstrate that it is possible to have an o-cycle of maximum length
and which meets every o-cycle, yet the o-cycle does not participate in
any o-colouring, is a vertex-orientation of a link projection of the Whitehead link.

\medspace
\noindent{\bf Example 4.1.}
We have assigned a vertex-orientation to the Whitehead link as shown
below. In (a), we have shown a 2-o-colouring, where the dotted curve is
an o-cycle of maximum length (four). In (b), for the same vertex-orientation,
we show an o-cycle of maximum length (dotted) which is not an o-cycle
for any o-colouring of the graph. Thus not every o-cycle of maximum
length is necessarily an o-cycle for some o-colouring.

\vskip8pt
\centerline{\table{c@{\hskip70pt}c}
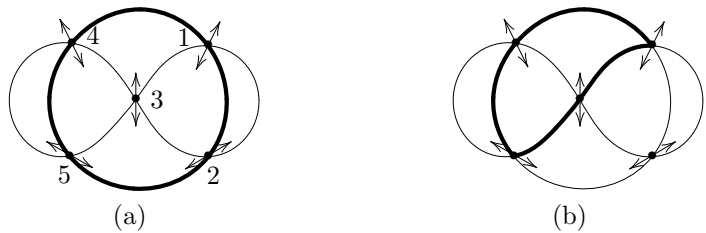
\hbox{
\xy /r.2pt/:,
(391.91488791533, 352.3413638559693)="1";
(488.23963523795584, 250.7659748989001)="flex6"; 
(304.63185851208146, 179.87197583238355)="flex7"; 
(255.21548123737043, 250.45696392221237)="3";
(217.7750076748735, 303.9474653663747)="flex8"; 
(16.54168002921108, 259.5628135887432)="flex9"; 
(130.07853913415153, 143.5014649620283)="5";
(209.87045061497656, 196.28196062290476)="flex10"; 
(307.14380480112897, 317.5545833255335)="flex11"; 
(427.17900728907483, 251.27837781359727)="flex12"; 
(392.2309370041537, 142.25732773272796)="2";
(251.42071093860625, 81.96490741602179)="flex13"; 
(91.92866235553595, 254.6097743849504)="flex14"; 
(134.3575638708463, 357.2138212500434)="4";
(263.9883631944421, 421.1817969437918)="flex15"; 
"1";"flex6"**\crv{ (445.34326907488906, 349.45561157773267) & (486.56031604159915, 304.72088971862195) },
"flex6";"2"*[o]=(0,0){\,}**\crv{ (490.0, 194.20714380008747) & (448.1102077772812, 145.47224681727857) },
"2"*[o]=(0,0){\,};"flex7"**\crv{ (359.0690074005095, 140.3494123485753) & (328.06252036965066, 156.1999502981455) },
"flex7";"3"**\crv{ (284.3314086465915, 200.38154446913785) & (270.2574733679984, 225.87474382493315) },
"3";"flex8"**\crv{ (243.84077093663188, 269.0459666106394) & (231.7982171140436, 287.25928530147655) },
"flex8";"4"*[o]=(0,0){\,}**\crv{ (195.77447028784547, 330.12898472767824) & (168.19341986481308, 352.62245744914907) },
"4"*[o]=(0,0){\,};"flex9"**\crv{ (75.12058734742762, 365.25199559776934) & (22.333924386248313, 320.02705403474664) },
"flex9";"5"**\crv{ (10.0, 191.27534333958482) & (65.15707992221807, 133.10145207054518) },
"5";"flex10"**\crv{ (162.52191950281167, 148.6986921681496) & (187.51546846253603, 172.21151232245245) },
"flex10";"3"*[o]=(0,0){\,}**\crv{ (225.90932145756702, 213.5516148226535) & (241.3109492912824, 231.4239015481662) },
"3"*[o]=(0,0){\,};"flex11"**\crv{ (271.9379099677385, 273.34727265786063) & (286.57430558515404, 297.99559277726166) },
"flex11";"1"**\crv{ (330.25163864549603, 339.527208862377) & (360.1381918399302, 354.0576740099409) },
"1"*[o]=(0,0){\,};"flex12"**[thicker][thicker]\crv{ (413.5408772608842, 322.9452190236518) & (425.9242918944621, 287.75123436850686) },
"flex12";"2"**[thicker][thicker]\crv{ (428.5332488010918, 211.9124334563536) & (416.81391758050387, 173.0244299991158) },
"2";"flex13"**[thicker][thicker]\crv{ (358.45789920296585, 99.98830731640845) & (305.47728776929796, 78.34366147388141) },
"flex13";"5"*[o]=(0,0){\,}**[thicker][thicker]\crv{ (204.27690609597215, 85.12306729207455) & (160.40651882788342, 107.27327845453624) },
"5"*[o]=(0,0){\,};"flex14"**[thicker][thicker]\crv{ (104.03989213273773, 174.6058434763985) & (89.8725403309002, 214.1070463173396) },
"flex14";"4"**[thicker][thicker]\crv{ (93.8554646006644, 292.56508367080056) & (109.87496835320914, 328.1437560715922) },
"4";"flex15"**[thicker][thicker]\crv{ (166.99456697019883, 395.96624056130514) & (213.445657049406, 421.6563385261186) },
"flex15";"1"*[o]=(0,0){\,}**[thicker][thicker]\crv{ (315.0547785053475, 420.70233829636896) & (361.5829581705556, 393.5714683184209) },
"1"*!<9pt,-3pt>{1},
"2"*!<-2pt,7pt>{2},
"3"*!<-8pt,0pt>{3},
"4"*!<-8pt,-3pt>{4},
"5"*!<2pt,7pt>{5},
"1"*!<0pt,0pt>{\smallbullet},
"2"*!<0pt,0pt>{\smallbullet},
"3"*!<0pt,0pt>{\smallbullet},
"4"*!<0pt,0pt>{\smallbullet},
"5"*!<0pt,0pt>{\smallbullet},
a(63);a(243)**\dir{},
"1"+/10pt/="1x";"1"+/-10pt/="1y"**\dir{-},
"1x"*\dir{<},"1y"*\dir{>},
a(30);a(210)**\dir{},
"2"+<0pt,.5pt>+/10pt/="2x";"2"+<-0pt,.5pt>+/-10pt/="2y"**\dir{-},
"2x"*\dir{<},"2y"*\dir{>},
a(-90);a(90)**\dir{},
"3"+/10pt/="3x";"3"+/-10pt/="3y"**\dir{-},
"3x"*\dir{<},"3y"*\dir{>},
a(-63);a(-243)**\dir{},
"4"+/10pt/="4x";"4"+/-10pt/="4y"**\dir{-},
"4x"*\dir{<},"4y"*\dir{>},
a(-30);a(-210)**\dir{},
"5"+<0pt,.5pt>+/10pt/="5x";"5"+<-0pt,.5pt>+/-10pt/="5y"**\dir{-},
"5x"*\dir{<},"5y"*\dir{>},
\endxy}
&
\hbox{\xy /r.2pt/:,
(391.91488791533, 352.3413638559693)="1";
(488.23963523795584, 250.7659748989001)="flex6"; 
(304.63185851208146, 179.87197583238355)="flex7"; 
(255.21548123737043, 250.45696392221237)="3";
(217.7750076748735, 303.9474653663747)="flex8"; 
(16.54168002921108, 259.5628135887432)="flex9"; 
(130.07853913415153, 143.5014649620283)="5";
(209.87045061497656, 196.28196062290476)="flex10"; 
(307.14380480112897, 317.5545833255335)="flex11"; 
(427.17900728907483, 251.27837781359727)="flex12"; 
(392.2309370041537, 142.25732773272796)="2";
(251.42071093860625, 81.96490741602179)="flex13"; 
(91.92866235553595, 254.6097743849504)="flex14"; 
(134.3575638708463, 357.2138212500434)="4";
(263.9883631944421, 421.1817969437918)="flex15"; 
"1";"flex6"**\crv{ (445.34326907488906, 349.45561157773267) & (486.56031604159915, 304.72088971862195) },
"flex6";"2"*[o]=(0,0){\,}**\crv{ (490.0, 194.20714380008747) & (448.1102077772812, 145.47224681727857) },
"2"*[o]=(0,0){\,};"flex7"**\crv{ (359.0690074005095, 140.3494123485753) & (328.06252036965066, 156.1999502981455) },
"flex7";"3"**\crv{ (284.3314086465915, 200.38154446913785) & (270.2574733679984, 225.87474382493315) },
"3";"flex8"**\crv{ (243.84077093663188, 269.0459666106394) & (231.7982171140436, 287.25928530147655) },
"flex8";"4"*[o]=(0,0){\,}**\crv{ (195.77447028784547, 330.12898472767824) & (168.19341986481308, 352.62245744914907) },
"4"*[o]=(0,0){\,};"flex9"**\crv{ (75.12058734742762, 365.25199559776934) & (22.333924386248313, 320.02705403474664) },
"flex9";"5"**\crv{ (10.0, 191.27534333958482) & (65.15707992221807, 133.10145207054518) },
"5";"flex10"**[thicker][thicker]\crv{ (162.52191950281167, 148.6986921681496) & (187.51546846253603, 172.21151232245245) },
"flex10";"3"*[o]=(0,0){\,}**[thicker][thicker]\crv{ (225.90932145756702, 213.5516148226535) & (241.3109492912824, 231.4239015481662) },
"3"*[o]=(0,0){\,};"flex11"**[thicker][thicker]\crv{ (271.9379099677385, 273.34727265786063) & (286.57430558515404, 297.99559277726166) },
"flex11";"1"**[thicker][thicker]\crv{ (330.25163864549603, 339.527208862377) & (360.1381918399302, 354.0576740099409) },
"1"*[o]=(0,0){\,};"flex12"**\crv{ (413.5408772608842, 322.9452190236518) & (425.9242918944621, 287.75123436850686) },
"flex12";"2"**\crv{ (428.5332488010918, 211.9124334563536) & (416.81391758050387, 173.0244299991158) },
"2";"flex13"**\crv{ (358.45789920296585, 99.98830731640845) & (305.47728776929796, 78.34366147388141) },
"flex13";"5"*[o]=(0,0){\,}**\crv{ (204.27690609597215, 85.12306729207455) & (160.40651882788342, 107.27327845453624) },
"5"*[o]=(0,0){\,};"flex14"**[thicker][thicker]\crv{ (104.03989213273773, 174.6058434763985) & (89.8725403309002, 214.1070463173396) },
"flex14";"4"**[thicker][thicker]\crv{ (93.8554646006644, 292.56508367080056) & (109.87496835320914, 328.1437560715922) },
"4";"flex15"**[thicker][thicker]\crv{ (166.99456697019883, 395.96624056130514) & (213.445657049406, 421.6563385261186) },
"flex15";"1"*[o]=(0,0){\,}**[thicker][thicker]\crv{ (315.0547785053475, 420.70233829636896) & (361.5829581705556, 393.5714683184209) },
"1"*!<0pt,0pt>{\smallbullet},
"2"*!<0pt,0pt>{\smallbullet},
"3"*!<0pt,0pt>{\smallbullet},
"4"*!<0pt,0pt>{\smallbullet},
"5"*!<0pt,0pt>{\smallbullet},
a(63);a(243)**\dir{},
"1"+/10pt/="1x";"1"+/-10pt/="1y"**\dir{-},
"1x"*\dir{<},"1y"*\dir{>},
a(30);a(210)**\dir{},
"2"+<0pt,.5pt>+/10pt/="2x";"2"+<-0pt,.5pt>+/-10pt/="2y"**\dir{-},
"2x"*\dir{<},"2y"*\dir{>},
a(-90);a(90)**\dir{},
"3"+/10pt/="3x";"3"+/-10pt/="3y"**\dir{-},
"3x"*\dir{<},"3y"*\dir{>},
a(-63);a(-243)**\dir{},
"4"+/10pt/="4x";"4"+/-10pt/="4y"**\dir{-},
"4x"*\dir{<},"4y"*\dir{>},
a(-30);a(-210)**\dir{},
"5"+<0pt,.5pt>+/10pt/="5x";"5"+<-0pt,.5pt>+/-10pt/="5y"**\dir{-},
"5x"*\dir{<},"5y"*\dir{>},
\endxy}\\
\noalign{\vskip-15pt}
(a) & (b)
\endtable}
\medskip

\noindent{\bf Example 4.2.}
This next example, a vertex-orientation of one of the basic polyhedra ($8^*$ in Figure 6 of \conway), is 
interesting in that it contains an o-cycle of maximum length that does
not participate in any o-colouring.. For this vertex-orientation, there
were a total of twelve o-cycles, of lengths 3,4,5,6, and 7, and there were 4 different ways to
decompose the edge-set as an edge-disjoint union of o-cycles (what we have called an o-colouring,
although we have not assigned any colours to the o-cycles). There were two o-cycles of length 7,
and neither participated in any of the four o-colourings (listed in Table 1).

\medskip
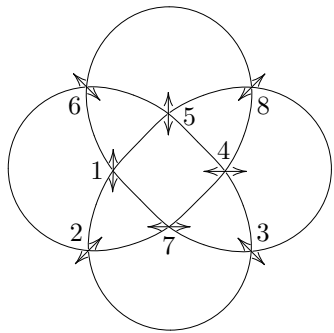
\table{l@{\hskip0pt}c}
\table[t]{lll}
O-cycle & Length & In o-colourings\\
1,2,3,4,5,8,6,1 & 7 &-\\
1,2,6,1 & 3 & 1,2,3\\
1,2,6,5,1& 4& 4\\
1,5,8,3,7,1& 5& 1\\
1,5,8,3,4,7,1& 6& 2\\
1,6,8,4,7,1& 5& 4\\
1,5,6,8,4,7,1& 6& 3\\
2,6,5,4,8,3,7,2& 7 &-\\
2,3,4,7,2& 4& 1\\
2,3,7,2& 3& 2,3,4\\
3,4,5,8,3& 4& 3,4\\
4,5,6,8,4& 4& 1,2\\
\endtable
&
\vbox to 0pt{\vskip.25in\hbox to 1in{\hskip0in
\xy /r.25pt/:,
(164,246.883)="1";
(138.363,195.579)="flex9"; 
(250.956,5.35538)="flex10"; 
(372.97,124.618)="3";
(360.934,194.119)="flex11"; 
(291.806,290.621)="flex12"; 
(247.555,333.211)="5";
(194.847,361.28)="flex13"; 
(5.40639,250.364)="flex14"; 
(126.818,125.817)="2";
(196.94,136.654)="flex15"; 
(291.991,201.795)="flex16"; 
(333.411,245.967)="4";
(361.597,299.41)="flex17"; 
(248.868,494.957)="flex18"; 
(123.927,373.58)="6";
(135.326,300.193)="flex19"; 
(205.009,202.868)="flex20"; 
(248.457,161.851)="7";
(301.274,134.724)="flex21"; 
(494.369,248.9)="flex22"; 
(373.625,373.095)="8";
(300.715,360.973)="flex23"; 
(203.676,290.969)="flex24"; 
"1";"flex9"**\crv{ (152.845,231.234) & (144.725,213.705) },
"flex9";"2"*[o]=(5,5){\,}**\crv{ (130.48,173.117) & (125.273,149.591) },
"2"*[o]=(5,5){\,};"flex10"**\crv{ (131.173,58.8004) & (184.622,5) },
"flex10";"3"**\crv{ (316.302,5.70547) & (368.725,58.6935) },
"3";"flex11"**\crv{ (374.499,148.354) & (369.2,171.835) },
"flex11";"4"*[o]=(5,5){\,}**\crv{ (354.094,212.558) & (345.224,230.236) },
"4"*[o]=(5,5){\,};"flex12"**\crv{ (321.173,262.264) & (306.065,276.077) },
"flex12";"5"**\crv{ (277.455,305.26) & (263.852,320.752) },
"5";"flex13"**\crv{ (231.65,345.371) & (213.611,354.322) },
"flex13";"6"*[o]=(5,5){\,}**\crv{ (172.101,369.714) & (148.148,375.237) },
"6"*[o]=(5,5){\,};"flex14"**\crv{ (57.878,369.063) & (5.24994,315.931) },
"flex14";"2"**\crv{ (5.56503,183.878) & (59.6302,130.642) },
"2";"flex15"**\crv{ (150.656,124.105) & (174.309,129.03) },
"flex15";"7"*[o]=(5,5){\,}**\crv{ (215.151,142.789) & (232.801,150.692) },
"7"*[o]=(5,5){\,};"flex16"**\crv{ (264.534,173.311) & (278.02,187.866) },
"flex16";"4"**\crv{ (306.304,216.064) & (321.257,229.803) },
"4";"flex17"**\crv{ (345.582,262.156) & (354.661,280.387) },
"flex17";"8"*[o]=(5,5){\,}**\crv{ (370.217,323.053) & (375.521,347.985) },
"8"*[o]=(5,5){\,};"flex18"**\crv{ (368.535,440.525) & (315.44,495) },
"flex18";"6"**\crv{ (182.537,494.914) & (129.66,440.751) },
"6";"flex19"**\crv{ (121.795,348.608) & (126.7,323.71) },
"flex19";"1"*[o]=(5,5){\,}**\crv{ (142.313,281.145) & (151.712,263.033) },
"1"*[o]=(5,5){\,};"flex20"**\crv{ (176.162,230.899) & (190.913,217.161) },
"flex20";"7"**\crv{ (219.015,188.667) & (232.479,173.813) },
"7";"flex21"**\crv{ (264.402,149.915) & (282.433,141.163) },
"flex21";"3"*[o]=(5,5){\,}**\crv{ (324.36,126.834) & (348.657,122.436) },
"3"*[o]=(5,5){\,};"flex22"**\crv{ (439.853,130.62) & (493.989,182.898) },
"flex22";"8"**\crv{ (494.75,315.17) & (440.712,368.319) },
"8";"flex23"**\crv{ (348.767,374.865) & (324.139,369.427) },
"flex23";"5"*[o]=(5,5){\,}**\crv{ (281.808,354.148) & (263.594,345.343) },
"5"*[o]=(5,5){\,};"flex24"**\crv{ (231.325,320.936) & (217.796,305.602) },
"flex24";"1"**\crv{ (189.928,276.722) & (175.5,263.015) },
"1"*!<6pt,0pt>{\hbox{1}},
"1";"1"+(0,32)**\dir{-}*\dir{>},
"1";"1"+(-,-32)**\dir{-}*\dir{>},
"2"*!<4.5pt,-7pt>{\hbox{2}},
"2";"2"+(20,20)**\dir{-}*\dir{>},
"2";"2"+(-20,-20)**\dir{-}*\dir{>},
"3"*!<-4.5pt,-7pt>{\hbox{3}},
"3";"3"+(20,-20)**\dir{-}*\dir{>},
"3";"3"+(-20,20)**\dir{-}*\dir{>},
"4"*!<.5pt,-8pt>{\hbox{4}},
"4";"4"+(32,0)**\dir{-}*\dir{>},
"4";"4"+(-32,0)**\dir{-}*\dir{>},
"5"*!<-8pt,1pt>{\hbox{5}},
"5";"5"+(0,32)**\dir{-}*\dir{>},
"5";"5"+(0,-32)**\dir{-}*\dir{>},
"6"*!<4.5pt,7pt>{\hbox{6}},
"6";"6"+(20,-20)**\dir{-}*\dir{>},
"6";"6"+(-20,20)**\dir{-}*\dir{>},
"7"*!<0pt,7pt>{\hbox{7}},
"7";"7"+(32,0)**\dir{-}*\dir{>},
"7";"7"+(-32,0)**\dir{-}*\dir{>},
"8"*!<-4.5pt,7pt>{\hbox{8}},
"8";"8"+(20,20)**\dir{-}*\dir{>},
"8";"8"+(-20,-20)**\dir{-}*\dir{>},
\endxy\hss}\vss}
\endtable

\medskip
The o-colourings for this vertex-orientation are (the colours assigned to each o-cycle
are shown in brackets):

\medskip
\centerline{\table{lllll}
Number     & 1               & 2                 & 3                 & 4\\
           & 1,5,8,3,7,1 (r) & 1,5,8,3,4,7,1 (r) & 1,5,6,8,4,7,1 (r) & 1,6,8,4,7,1 (r)\\
           & 4,5,6,8,4 (g)   & 4,5,6,8,4 (g)     & 3,4,5,8,3 (g)     & 3,4,5,8,3 (g)\\
           & 2,3,4,7,2 (b)   & 2,3,7,2 (b)       & 2,3,7,2 (b)       & 1,2,6,5,1 (g)\\
           & 1,2,6,1   (y)   & 1,2,6,1 (g)       & 1,2,6,1 (g)       & 2,3,7,2 (b)\\ 
	   \noalign{\vskip5pt}
	   \multicolumn{5}{c}{{ Table} 1}
	   \endtable}

\medskip
\leavevmode\hskip-\parindent\vbox{\hsize=3in We have shown the first o-cycle of length 7, namely $1,2,3,4,5,8,6,1$, and the other one
is obtained by reflecting this one across the axis of (vertex-orientation) symmetry through
vertices 2 and 8. It is evident that this o-cycle can't participate in an o-colouring of the graph, as
the cycle $7,4,8,3,7$ would have to be an o-cycle, and it fails to meet the requirement at vertex 7.}
\hskip .25in
\lower .01in
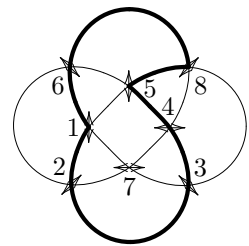
\hbox{
\xy /r.18pt/:,
(164,246.883)="1";
(138.363,195.579)="flex9"; 
(250.956,5.35538)="flex10"; 
(372.97,124.618)="3";
(360.934,194.119)="flex11"; 
(291.806,290.621)="flex12"; 
(247.555,333.211)="5";
(194.847,361.28)="flex13"; 
(5.40639,250.364)="flex14"; 
(126.818,125.817)="2";
(196.94,136.654)="flex15"; 
(291.991,201.795)="flex16"; 
(333.411,245.967)="4";
(361.597,299.41)="flex17"; 
(248.868,494.957)="flex18"; 
(123.927,373.58)="6";
(135.326,300.193)="flex19"; 
(205.009,202.868)="flex20"; 
(248.457,161.851)="7";
(301.274,134.724)="flex21"; 
(494.369,248.9)="flex22"; 
(373.625,373.095)="8";
(300.715,360.973)="flex23"; 
(203.676,290.969)="flex24"; 
"1";"flex9"**[thicker][thicker]\crv{ (152.845,231.234) & (144.725,213.705) },
"flex9";"2"*[o]=(0,0){\,}**[thicker][thicker]\crv{ (130.48,173.117) & (125.273,149.591) },
"2"*[o]=(0,0){\,};"flex10"**[thicker][thicker]\crv{ (131.173,58.8004) & (184.622,5) },
"flex10";"3"**[thicker][thicker]\crv{ (316.302,5.70547) & (368.725,58.6935) },
"3";"flex11"**[thicker][thicker]\crv{ (374.499,148.354) & (369.2,171.835) },
"flex11";"4"*[o]=(0,0){\,}**[thicker][thicker]\crv{ (354.094,212.558) & (345.224,230.236) },
"4"*[o]=(0,0){\,};"flex12"**[thicker][thicker]\crv{ (321.173,262.264) & (306.065,276.077) },
"flex12";"5"**[thicker][thicker]\crv{ (277.455,305.26) & (263.852,320.752) },
"5";"flex13"**\crv{ (231.65,345.371) & (213.611,354.322) },
"flex13";"6"*[o]=(0,0){\,}**\crv{ (172.101,369.714) & (148.148,375.237) },
"6"*[o]=(0,0){\,};"flex14"**\crv{ (57.878,369.063) & (5.24994,315.931) },
"flex14";"2"**\crv{ (5.56503,183.878) & (59.6302,130.642) },
"2";"flex15"**\crv{ (150.656,124.105) & (174.309,129.03) },
"flex15";"7"*[o]=(0,0){\,}**\crv{ (215.151,142.789) & (232.801,150.692) },
"7"*[o]=(0,0){\,};"flex16"**\crv{ (264.534,173.311) & (278.02,187.866) },
"flex16";"4"**\crv{ (306.304,216.064) & (321.257,229.803) },
"4";"flex17"**\crv{ (345.582,262.156) & (354.661,280.387) },
"flex17";"8"*[o]=(0,0){\,}**\crv{ (370.217,323.053) & (375.521,347.985) },
"8"*[o]=(0,0){\,};"flex18"**[thicker][thicker]\crv{ (368.535,440.525) & (315.44,495) },
"flex18";"6"**[thicker][thicker]\crv{ (182.537,494.914) & (129.66,440.751) },
"6";"flex19"**[thicker][thicker]\crv{ (121.795,348.608) & (126.7,323.71) },
"flex19";"1"*[o]=(0,0){\,}**[thicker][thicker]\crv{ (142.313,281.145) & (151.712,263.033) },
"1"*[o]=(0,0){\,};"flex20"**\crv{ (176.162,230.899) & (190.913,217.161) },
"flex20";"7"**\crv{ (219.015,188.667) & (232.479,173.813) },
"7";"flex21"**\crv{ (264.402,149.915) & (282.433,141.163) },
"flex21";"3"*[o]=(0,0){\,}**\crv{ (324.36,126.834) & (348.657,122.436) },
"3"*[o]=(0,0){\,};"flex22"**\crv{ (439.853,130.62) & (493.989,182.898) },
"flex22";"8"**\crv{ (494.75,315.17) & (440.712,368.319) },
"8";"flex23"**[thicker][thicker]\crv{ (348.767,374.865) & (324.139,369.427) },
"flex23";"5"*[o]=(0,0){\,}**[thicker][thicker]\crv{ (281.808,354.148) & (263.594,345.343) },
"5"*[o]=(0,0){\,};"flex24"**\crv{ (231.325,320.936) & (217.796,305.602) },
"flex24";"1"**\crv{ (189.928,276.722) & (175.5,263.015) },
"1"*!<6pt,0pt>{\hbox{1}},
"1";"1"+(0,32)**\dir{-}*\dir{>},
"1";"1"+(-,-32)**\dir{-}*\dir{>},
"2"*!<4.5pt,-7pt>{\hbox{2}},
"2";"2"+(20,20)**\dir{-}*\dir{>},
"2";"2"+(-20,-20)**\dir{-}*\dir{>},
"3"*!<-4.5pt,-7pt>{\hbox{3}},
"3";"3"+(20,-20)**\dir{-}*\dir{>},
"3";"3"+(-20,20)**\dir{-}*\dir{>},
"4"*!<.5pt,-8pt>{\hbox{4}},
"4";"4"+(32,0)**\dir{-}*\dir{>},
"4";"4"+(-32,0)**\dir{-}*\dir{>},
"5"*!<-8pt,1pt>{\hbox{5}},
"5";"5"+(0,32)**\dir{-}*\dir{>},
"5";"5"+(0,-32)**\dir{-}*\dir{>},
"6"*!<4.5pt,7pt>{\hbox{6}},
"6";"6"+(20,-20)**\dir{-}*\dir{>},
"6";"6"+(-20,20)**\dir{-}*\dir{>},
"7"*!<0pt,7pt>{\hbox{7}},
"7";"7"+(32,0)**\dir{-}*\dir{>},
"7";"7"+(-32,0)**\dir{-}*\dir{>},
"8"*!<-4.5pt,7pt>{\hbox{8}},
"8";"8"+(20,20)**\dir{-}*\dir{>},
"8";"8"+(-20,-20)**\dir{-}*\dir{>},
\endxy}

\medspace
\noindent{\bf Example 4.3.}
We complete our discussion with a case study of the basic polyhedral graph $6^*$, the 4-regular
simple graph that is obtained as an alternating six crossing projection of the borromean rings.
This graph has an automorphism group of size 48, and the natural action of the
automorphism group on the set of vertex orientations of $6^*$ has seven orbits. We offer a representative
of each orbit below, and for each, we present the complete collection of o-cycles, as well as every
way of decomposing the edge set into edge-disjoint o-cycles (what we refer to as o-colourings).
For each, we label the different o-colourings with indices based at 0, and then for each o-cycle, we indicate
its length and, by listing the indices, the different o-colourings in which the o-cycle participates.

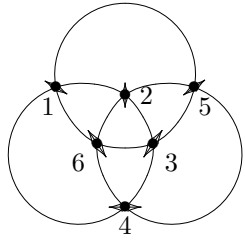
\begin{tabular}{cc}
 \begin{tabular}{c}
\hbox{\xy /r.20pt/:,
(119.113,311.02)="1";
(191.876,317.16)="flex7"; 
(283.884,255.041)="flex8"; 
(302.887,204.917)="3";
(291.639,144.326)="flex9"; 
(46.7613,118.174)="flex10"; 
(249.972,470.153)="flex11"; 
(380.834,311.02)="5";
(349.774,244.943)="flex12"; 
(249.969,196.31)="flex13"; 
(197.053,204.892)="6";
(150.169,244.915)="flex14"; 
(249.984,296.571)="2";
(308.08,317.151)="flex15"; 
(453.206,118.124)="flex16"; 
(249.947,84.3883)="4";
(208.268,144.307)="flex17"; 
(216.073,255.021)="flex18"; 
"1";"flex7"**\crv{ (142.378,319.265) & (167.397,320.203) },
"flex7";"2"*[o]=(5,5){\,}**\crv{ (212.611,314.583) & (233.207,309.056) },
"2"*[o]=(5,5){\,};"flex8"**\crv{ (264.485,285.78) & (274.858,270.676) },
"flex8";"3"**\crv{ (292.909,239.406) & (300.8,222.871) },
"3";"flex9"**\crv{ (305.301,184.151) & (299.777,163.562) },
"flex9";"4"*[o]=(5,5){\,}**\crv{ (282.03,121.611) & (268.714,100.414) },
"4"*[o]=(5,5){\,};"flex10"**\crv{ (186.091,29.8623) & (88.5278,45.8136) },
"flex10";"1"**\crv{ (5,190.526) & (39.9736,282.974) },
"1"*[o]=(5,5){\,};"flex11"**\crv{ (103.842,393.609) & (166.419,470.147) },
"flex11";"5"**\crv{ (333.535,470.16) & (396.127,393.611) },
"5";"flex12"**\crv{ (376.341,286.753) & (364.647,264.619) },
"flex12";"3"*[o]=(5,5){\,}**\crv{ (337.174,228.275) & (322.087,213.205) },
"3"*[o]=(5,5){\,};"flex13"**\crv{ (286.29,197.752) & (268.024,196.315) },
"flex13";"6"**\crv{ (231.916,196.306) & (213.65,197.733) },
"6";"flex14"**\crv{ (177.851,213.175) & (162.766,228.246) },
"flex14";"1"*[o]=(5,5){\,}**\crv{ (135.293,264.598) & (123.601,286.744) },
"2";"flex15"**\crv{ (266.758,309.053) & (287.35,314.575) },
"flex15";"5"*[o]=(5,5){\,}**\crv{ (332.556,320.193) & (357.571,319.259) },
"5"*[o]=(5,5){\,};"flex16"**\crv{ (460.009,282.979) & (495,190.491) },
"flex16";"4"**\crv{ (411.413,45.7587) & (313.821,29.8402) },
"4";"flex17"**\crv{ (231.186,100.41) & (217.868,121.596) },
"flex17";"6"*[o]=(5,5){\,}**\crv{ (200.137,163.542) & (194.63,184.131) },
"6"*[o]=(5,5){\,};"flex18"**\crv{ (199.149,222.848) & (207.046,239.384) },
"flex18";"2"**\crv{ (225.103,270.663) & (235.476,285.776) },
"1"*{\bullet}*!<3pt,7pt>{\hbox{1}},
"2"*{\bullet}*!<-8pt,0pt>{\hbox{2}},
"3"*{\bullet}*!<-7pt,7pt>{\hbox{3}},
"4"*{\bullet}*!<0pt,7pt>{\hbox{4}},
"5"*{\bullet}*!<-4pt,7pt>{\hbox{5}},
"6"*{\bullet}*!<7pt,7pt>{\hbox{6}},
"flex7";"flex14"**\dir{}?(.5)="1:2-6",
"1";"1:2-6"**\dir{}?(.4)="x","1";"x"**\dir{-}*\dir{>},
"flex10";"flex11"**\dir{}?(.5)="1:4-5",
"1";"1:4-5"**\dir{}?(.7)="y",
"1";"1"+"1"-"y"**\dir{-}*\dir{>},
"flex7";"flex15"**\dir{}?(.5)="2:1-5",
"2";"2:1-5"**\dir{}?(1.1)="y",
"2";"y"**\dir{-}*\dir{>},
"2";"2"+"2"-"y"**\dir{-}*\dir{>},
"3";"3"+(-17,-25)**\dir{-}*\dir{>},
"3";"3"-(-17,-25)**\dir{-}*\dir{>},
"4";"4"+(-30,0)**\dir{-}*\dir{>},
"4";"4"+(30,0)**\dir{-}*\dir{>},
"5";"5"+(20,14)**\dir{-}*\dir{>},
"5";"5"-(20,14)**\dir{-}*\dir{>},
"6";"6"+(17,-25)**\dir{-}*\dir{>},
"6";"6"-(17,-25)**\dir{-}*\dir{>},
\endxy} \\
 \begin{tabular}{ll}
  0 & (1,2,3,4,1),(1,5,3,6,1),(2,5,4,6,2)\\
  1 & (1,2,6,3,5,4,1),(1,5,2,3,4,6,1)\\
  2 & (1,2,6,4,3,5,1),(1,4,5,2,3,6,1)\\
 \end{tabular}
\end{tabular}
&
{
\begin{tabular}{lll}
\multicolumn{3}{c}{List of o-cycles (11)}\\
 (1,2,6,3,5,4,1) & 6 & 0\\
 (1,4,5,2,3,6,1) & 6 & 1\\
 (1,5,2,3,4,6,1) & 6 & 2\\
 (1,2,6,4,3,5,1) & 6 &  \\
 (1,4,5,3,6,1)   & 5 & 0\\
 (1,5,2,3,6,1)   & 5 & 2\\
 (1,2,6,3,5,1)   & 5 &\\
 (1,5,3,4,6,1)   & 5 &\\
 (2,5,4,6,2)     & 4 & 0\\
 (1,2,3,4,1)     & 4 & 1\\
 (1,5,3,6,1)     & 4 & \\
\end{tabular}} 
\end{tabular}

\medskip
\centerline{\begin{tabular}{cc}
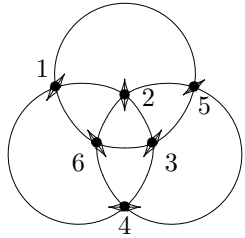
 \begin{tabular}{c}
\hbox{\xy /r.20pt/:,
(119.113,311.02)="1";
(191.876,317.16)="flex7"; 
(283.884,255.041)="flex8"; 
(302.887,204.917)="3";
(291.639,144.326)="flex9"; 
(46.7613,118.174)="flex10"; 
(249.972,470.153)="flex11"; 
(380.834,311.02)="5";
(349.774,244.943)="flex12"; 
(249.969,196.31)="flex13"; 
(197.053,204.892)="6";
(150.169,244.915)="flex14"; 
(249.984,296.571)="2";
(308.08,317.151)="flex15"; 
(453.206,118.124)="flex16"; 
(249.947,84.3883)="4";
(208.268,144.307)="flex17"; 
(216.073,255.021)="flex18"; 
"1";"flex7"**\crv{ (142.378,319.265) & (167.397,320.203) },
"flex7";"2"*[o]=(5,5){\,}**\crv{ (212.611,314.583) & (233.207,309.056) },
"2"*[o]=(5,5){\,};"flex8"**\crv{ (264.485,285.78) & (274.858,270.676) },
"flex8";"3"**\crv{ (292.909,239.406) & (300.8,222.871) },
"3";"flex9"**\crv{ (305.301,184.151) & (299.777,163.562) },
"flex9";"4"*[o]=(5,5){\,}**\crv{ (282.03,121.611) & (268.714,100.414) },
"4"*[o]=(5,5){\,};"flex10"**\crv{ (186.091,29.8623) & (88.5278,45.8136) },
"flex10";"1"**\crv{ (5,190.526) & (39.9736,282.974) },
"1"*[o]=(5,5){\,};"flex11"**\crv{ (103.842,393.609) & (166.419,470.147) },
"flex11";"5"**\crv{ (333.535,470.16) & (396.127,393.611) },
"5";"flex12"**\crv{ (376.341,286.753) & (364.647,264.619) },
"flex12";"3"*[o]=(5,5){\,}**\crv{ (337.174,228.275) & (322.087,213.205) },
"3"*[o]=(5,5){\,};"flex13"**\crv{ (286.29,197.752) & (268.024,196.315) },
"flex13";"6"**\crv{ (231.916,196.306) & (213.65,197.733) },
"6";"flex14"**\crv{ (177.851,213.175) & (162.766,228.246) },
"flex14";"1"*[o]=(5,5){\,}**\crv{ (135.293,264.598) & (123.601,286.744) },
"2";"flex15"**\crv{ (266.758,309.053) & (287.35,314.575) },
"flex15";"5"*[o]=(5,5){\,}**\crv{ (332.556,320.193) & (357.571,319.259) },
"5"*[o]=(5,5){\,};"flex16"**\crv{ (460.009,282.979) & (495,190.491) },
"flex16";"4"**\crv{ (411.413,45.7587) & (313.821,29.8402) },
"4";"flex17"**\crv{ (231.186,100.41) & (217.868,121.596) },
"flex17";"6"*[o]=(5,5){\,}**\crv{ (200.137,163.542) & (194.63,184.131) },
"6"*[o]=(5,5){\,};"flex18"**\crv{ (199.149,222.848) & (207.046,239.384) },
"flex18";"2"**\crv{ (225.103,270.663) & (235.476,285.776) },
"1"*{\bullet}*!<5pt,-7pt>{\hbox{1}},
"2"*{\bullet}*!<-9pt,1pt>{\hbox{2}},
"3"*{\bullet}*!<-7pt,7pt>{\hbox{3}},
"4"*{\bullet}*!<0pt,7pt>{\hbox{4}},
"5"*{\bullet}*!<-4pt,7pt>{\hbox{5}},
"6"*{\bullet}*!<7pt,7pt>{\hbox{6}},
"1";"1"+(17,25)**\dir{-}*\dir{>},
"1";"1"-(17,25)**\dir{-}*\dir{>},
"2";"2"+(0,30)**\dir{-}*\dir{>},
"2";"2"-(0,30)**\dir{-}*\dir{>},
"3";"3"+(-17,-25)**\dir{-}*\dir{>},
"3";"3"-(-17,-25)**\dir{-}*\dir{>},
"4";"4"+(-30,0)**\dir{-}*\dir{>},
"4";"4"+(30,0)**\dir{-}*\dir{>},
"5";"5"+(20,14)**\dir{-}*\dir{>},
"5";"5"-(20,14)**\dir{-}*\dir{>},
"6";"6"+(17,-25)**\dir{-}*\dir{>},
"6";"6"-(17,-25)**\dir{-}*\dir{>},
\endxy}\\
\begin{tabular}{ll}
 0 & (1,2,3,4,1),(1,5,3,6,1),(2,5,4,6,2)\\
 1 & (1,2,3,6,1),(1,4,3,5,1),(2,5,4,6,2)\\
 2 & (1,2,6,3,5,4,1),(1,5,2,3,4,6,1)\\
\end{tabular}
\end{tabular}
&
\begin{tabular}{lll}
\multicolumn{3}{c}{List of o-cycles (11)}\\
(1,2,6,3,5,4,1) & 6 & 0\\
(1,5,2,3,4,6,1) & 6 &\\
(1,4,3,2,5,1)   & 5 & 0\\
(1,2,3,4,6,1)   & 5 & 1\\
(1,5,2,3,6,1)   & 5 & 2\\
(1,5,3,4,6,1)   & 5 &\\
(1,4,3,5,1)     & 4 & 1\\
(1,2,3,4,1)     & 4 & 2\\
(1,2,3,6,1)     & 4 &\\
(1,5,3,6,1)     & 4 &\\
(2,5,4,6,2)     & 4 & \setbox0=\hbox{0}\hbox to \wd0{0,1\hss}\\
\end{tabular}
\end{tabular}}

\medskip
\centerline{\table{cc}
\table{c}
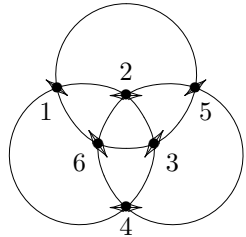
\hbox{\xy /r.20pt/:,
(119.113,311.02)="1";
(191.876,317.16)="flex7"; 
(283.884,255.041)="flex8"; 
(302.887,204.917)="3";
(291.639,144.326)="flex9"; 
(46.7613,118.174)="flex10"; 
(249.972,470.153)="flex11"; 
(380.834,311.02)="5";
(349.774,244.943)="flex12"; 
(249.969,196.31)="flex13"; 
(197.053,204.892)="6";
(150.169,244.915)="flex14"; 
(249.984,296.571)="2";
(308.08,317.151)="flex15"; 
(453.206,118.124)="flex16"; 
(249.947,84.3883)="4";
(208.268,144.307)="flex17"; 
(216.073,255.021)="flex18"; 
"1";"flex7"**\crv{ (142.378,319.265) & (167.397,320.203) },
"flex7";"2"*[o]=(5,5){\,}**\crv{ (212.611,314.583) & (233.207,309.056) },
"2"*[o]=(5,5){\,};"flex8"**\crv{ (264.485,285.78) & (274.858,270.676) },
"flex8";"3"**\crv{ (292.909,239.406) & (300.8,222.871) },
"3";"flex9"**\crv{ (305.301,184.151) & (299.777,163.562) },
"flex9";"4"*[o]=(5,5){\,}**\crv{ (282.03,121.611) & (268.714,100.414) },
"4"*[o]=(5,5){\,};"flex10"**\crv{ (186.091,29.8623) & (88.5278,45.8136) },
"flex10";"1"**\crv{ (5,190.526) & (39.9736,282.974) },
"1"*[o]=(5,5){\,};"flex11"**\crv{ (103.842,393.609) & (166.419,470.147) },
"flex11";"5"**\crv{ (333.535,470.16) & (396.127,393.611) },
"5";"flex12"**\crv{ (376.341,286.753) & (364.647,264.619) },
"flex12";"3"*[o]=(5,5){\,}**\crv{ (337.174,228.275) & (322.087,213.205) },
"3"*[o]=(5,5){\,};"flex13"**\crv{ (286.29,197.752) & (268.024,196.315) },
"flex13";"6"**\crv{ (231.916,196.306) & (213.65,197.733) },
"6";"flex14"**\crv{ (177.851,213.175) & (162.766,228.246) },
"flex14";"1"*[o]=(5,5){\,}**\crv{ (135.293,264.598) & (123.601,286.744) },
"2";"flex15"**\crv{ (266.758,309.053) & (287.35,314.575) },
"flex15";"5"*[o]=(5,5){\,}**\crv{ (332.556,320.193) & (357.571,319.259) },
"5"*[o]=(5,5){\,};"flex16"**\crv{ (460.009,282.979) & (495,190.491) },
"flex16";"4"**\crv{ (411.413,45.7587) & (313.821,29.8402) },
"4";"flex17"**\crv{ (231.186,100.41) & (217.868,121.596) },
"flex17";"6"*[o]=(5,5){\,}**\crv{ (200.137,163.542) & (194.63,184.131) },
"6"*[o]=(5,5){\,};"flex18"**\crv{ (199.149,222.848) & (207.046,239.384) },
"flex18";"2"**\crv{ (225.103,270.663) & (235.476,285.776) },
"1"*{\bullet}*!<4pt,9pt>{\hbox{1}},
"2"*{\bullet}*!<0pt,-9pt>{\hbox{2}},
"3"*{\bullet}*!<-7pt,7pt>{\hbox{3}},
"4"*{\bullet}*!<0pt,7pt>{\hbox{4}},
"5"*{\bullet}*!<-4pt,9pt>{\hbox{5}},
"6"*{\bullet}*!<7pt,7pt>{\hbox{6}},
%
"1";"1"+(20,-15)**\dir{-}*\dir{>},
"1";"1"-(20,-15)**\dir{-}*\dir{>},
"2";"2"+(30,0)**\dir{-}*\dir{>},
"2";"2"-(30,0)**\dir{-}*\dir{>},
"3";"3"+(-17,-25)**\dir{-}*\dir{>},
"3";"3"-(-17,-25)**\dir{-}*\dir{>},
"4";"4"+(-30,0)**\dir{-}*\dir{>},
"4";"4"+(30,0)**\dir{-}*\dir{>},
"5";"5"+(20,15)**\dir{-}*\dir{>},
"5";"5"-(20,15)**\dir{-}*\dir{>},
"6";"6"+(17,-25)**\dir{-}*\dir{>},
"6";"6"-(17,-25)**\dir{-}*\dir{>},
\endxy}\\
\table{ll}
0 & (1,2,3,4,1),(1,5,3,6,1),(2,5,4,6,2)\\
1 & (1,2,5,1),(1,4,5,3,6,1),(2,3,4,6,2)\\
2 & (1,2,5,4,1),(1,5,3,6,1),(2,3,4,6,2)\\
3 & (1,2,5,4,1),(1,5,3,4,6,1),(2,3,6,2)\\
\endtable
\endtable
&
\table{lll}
\multicolumn{3}{c}{List of o-cycles (9)}\\
(1,4,5,3,6,1) & 5 & 0\\
(1,5,3,4,6,1) & 5 & 1\\
(1,2,3,4,1)   & 4 & \setbox0=\hbox{0}\hbox to \wd0{2,3\hss}\\
(1,2,5,4,1)   & 4 & \setbox0=\hbox{0}\hbox to \wd0{0,2\hss}\\
(1,5,3,6,1)   & 4 & 1\\
(2,5,4,6,2)   & 4 & 3\\
(2,3,4,6,2)   & 4 & 0\\
(1,2,5,1)     & 3 & \setbox0=\hbox{0}\hbox to \wd0{1,2\hss}\\
(2,3,6,2)     & 3 & 3\\
\endtable
\endtable}

\medskip
\centering\table{cc}
\table{c}
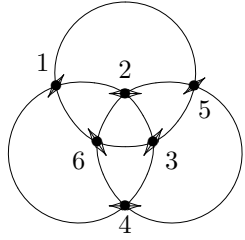
\hbox{\xy /r.20pt/:,
(119.113,311.02)="1";
(191.876,317.16)="flex7"; 
(283.884,255.041)="flex8"; 
(302.887,204.917)="3";
(291.639,144.326)="flex9"; 
(46.7613,118.174)="flex10"; 
(249.972,470.153)="flex11"; 
(380.834,311.02)="5";
(349.774,244.943)="flex12"; 
(249.969,196.31)="flex13"; 
(197.053,204.892)="6";
(150.169,244.915)="flex14"; 
(249.984,296.571)="2";
(308.08,317.151)="flex15"; 
(453.206,118.124)="flex16"; 
(249.947,84.3883)="4";
(208.268,144.307)="flex17"; 
(216.073,255.021)="flex18"; 
"1";"flex7"**\crv{ (142.378,319.265) & (167.397,320.203) },
"flex7";"2"*[o]=(5,5){\,}**\crv{ (212.611,314.583) & (233.207,309.056) },
"2"*[o]=(5,5){\,};"flex8"**\crv{ (264.485,285.78) & (274.858,270.676) },
"flex8";"3"**\crv{ (292.909,239.406) & (300.8,222.871) },
"3";"flex9"**\crv{ (305.301,184.151) & (299.777,163.562) },
"flex9";"4"*[o]=(5,5){\,}**\crv{ (282.03,121.611) & (268.714,100.414) },
"4"*[o]=(5,5){\,};"flex10"**\crv{ (186.091,29.8623) & (88.5278,45.8136) },
"flex10";"1"**\crv{ (5,190.526) & (39.9736,282.974) },
"1"*[o]=(5,5){\,};"flex11"**\crv{ (103.842,393.609) & (166.419,470.147) },
"flex11";"5"**\crv{ (333.535,470.16) & (396.127,393.611) },
"5";"flex12"**\crv{ (376.341,286.753) & (364.647,264.619) },
"flex12";"3"*[o]=(5,5){\,}**\crv{ (337.174,228.275) & (322.087,213.205) },
"3"*[o]=(5,5){\,};"flex13"**\crv{ (286.29,197.752) & (268.024,196.315) },
"flex13";"6"**\crv{ (231.916,196.306) & (213.65,197.733) },
"6";"flex14"**\crv{ (177.851,213.175) & (162.766,228.246) },
"flex14";"1"*[o]=(5,5){\,}**\crv{ (135.293,264.598) & (123.601,286.744) },
"2";"flex15"**\crv{ (266.758,309.053) & (287.35,314.575) },
"flex15";"5"*[o]=(5,5){\,}**\crv{ (332.556,320.193) & (357.571,319.259) },
"5"*[o]=(5,5){\,};"flex16"**\crv{ (460.009,282.979) & (495,190.491) },
"flex16";"4"**\crv{ (411.413,45.7587) & (313.821,29.8402) },
"4";"flex17"**\crv{ (231.186,100.41) & (217.868,121.596) },
"flex17";"6"*[o]=(5,5){\,}**\crv{ (200.137,163.542) & (194.63,184.131) },
"6"*[o]=(5,5){\,};"flex18"**\crv{ (199.149,222.848) & (207.046,239.384) },
"flex18";"2"**\crv{ (225.103,270.663) & (235.476,285.776) },
"1"*{\bullet}*!<5pt,-9pt>{\hbox{1}},
"2"*{\bullet}*!<0pt,-9pt>{\hbox{2}},
"3"*{\bullet}*!<-7pt,7pt>{\hbox{3}},
"4"*{\bullet}*!<0pt,7pt>{\hbox{4}},
"5"*{\bullet}*!<-4pt,9pt>{\hbox{5}},
"6"*{\bullet}*!<7pt,7pt>{\hbox{6}},
%
"1";"1"+(15,21)**\dir{-}*\dir{>},
"1";"1"-(15,21)**\dir{-}*\dir{>},
"2";"2"+(30,0)**\dir{-}*\dir{>},
"2";"2"-(30,0)**\dir{-}*\dir{>},
"3";"3"+(-17,-25)**\dir{-}*\dir{>},
"3";"3"-(-17,-25)**\dir{-}*\dir{>},
"4";"4"+(-30,0)**\dir{-}*\dir{>},
"4";"4"+(30,0)**\dir{-}*\dir{>},
"5";"5"+(20,15)**\dir{-}*\dir{>},
"5";"5"-(20,15)**\dir{-}*\dir{>},
"6";"6"+(17,-25)**\dir{-}*\dir{>},
"6";"6"-(17,-25)**\dir{-}*\dir{>},
\endxy}\\
\table{ll}
 0 & (1,2,3,4,1),(1,5,3,6,1),(2,5,4,6,2)\\
 1 & (1,2,3,6,1),(1,4,3,5,1),(2,5,4,6,2)\\
 2 & (1,2,5,4,1),(1,5,3,6,1),(2,3,4,6,2)\\
 3 & (1,2,5,4,1),(1,5,3,4,6,1),(2,3,6,2)\\
 4 & (1,2,5,4,6,1),(1,4,3,5,1),(2,3,6,2)\\
\endtable
\endtable
&
\table{lll}
\multicolumn{3}{c}{List of o-cycles (11)}\\
(1,2,3,4,6,1) & 5 & 0\\
(1,5,3,4,6,1) & 5 & 1\\
(1,2,5,4,6,1) & 5 & \\
(2,5,4,6,2)   & 4 & \setbox0=\hbox{0}\hbox to \wd0{0,1\hss}\\
(1,2,5,4,1)   & 4 & \setbox0=\hbox{0}\hbox to \wd0{0,2\hss}\\
(1,4,3,5,1)   & 4 & \setbox0=\hbox{0}\hbox to \wd0{1,4\hss}\\
(1,2,3,4,1)   & 4 & \setbox0=\hbox{0}\hbox to \wd0{2,3\hss}\\
(1,2,3,6,1)   & 4 & 4\\
(1,5,3,6,1)   & 4 & 3\\
(2,3,4,6,2)   & 4 & 2\\
(2,3,6,2)     & 3 & \setbox0=\hbox{0}\hbox to \wd0{3,4\hss}\\
\endtable
\endtable

\medskip
\centerline{\table{cc}
\table{c}
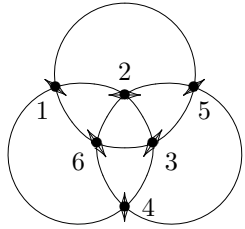
\hbox{\xy /r.20pt/:,
(119.113,311.02)="1";
(191.876,317.16)="flex7"; 
(283.884,255.041)="flex8"; 
(302.887,204.917)="3";
(291.639,144.326)="flex9"; 
(46.7613,118.174)="flex10"; 
(249.972,470.153)="flex11"; 
(380.834,311.02)="5";
(349.774,244.943)="flex12"; 
(249.969,196.31)="flex13"; 
(197.053,204.892)="6";
(150.169,244.915)="flex14"; 
(249.984,296.571)="2";
(308.08,317.151)="flex15"; 
(453.206,118.124)="flex16"; 
(249.947,84.3883)="4";
(208.268,144.307)="flex17"; 
(216.073,255.021)="flex18"; 
"1";"flex7"**\crv{ (142.378,319.265) & (167.397,320.203) },
"flex7";"2"*[o]=(5,5){\,}**\crv{ (212.611,314.583) & (233.207,309.056) },
"2"*[o]=(5,5){\,};"flex8"**\crv{ (264.485,285.78) & (274.858,270.676) },
"flex8";"3"**\crv{ (292.909,239.406) & (300.8,222.871) },
"3";"flex9"**\crv{ (305.301,184.151) & (299.777,163.562) },
"flex9";"4"*[o]=(5,5){\,}**\crv{ (282.03,121.611) & (268.714,100.414) },
"4"*[o]=(5,5){\,};"flex10"**\crv{ (186.091,29.8623) & (88.5278,45.8136) },
"flex10";"1"**\crv{ (5,190.526) & (39.9736,282.974) },
"1"*[o]=(5,5){\,};"flex11"**\crv{ (103.842,393.609) & (166.419,470.147) },
"flex11";"5"**\crv{ (333.535,470.16) & (396.127,393.611) },
"5";"flex12"**\crv{ (376.341,286.753) & (364.647,264.619) },
"flex12";"3"*[o]=(5,5){\,}**\crv{ (337.174,228.275) & (322.087,213.205) },
"3"*[o]=(5,5){\,};"flex13"**\crv{ (286.29,197.752) & (268.024,196.315) },
"flex13";"6"**\crv{ (231.916,196.306) & (213.65,197.733) },
"6";"flex14"**\crv{ (177.851,213.175) & (162.766,228.246) },
"flex14";"1"*[o]=(5,5){\,}**\crv{ (135.293,264.598) & (123.601,286.744) },
"2";"flex15"**\crv{ (266.758,309.053) & (287.35,314.575) },
"flex15";"5"*[o]=(5,5){\,}**\crv{ (332.556,320.193) & (357.571,319.259) },
"5"*[o]=(5,5){\,};"flex16"**\crv{ (460.009,282.979) & (495,190.491) },
"flex16";"4"**\crv{ (411.413,45.7587) & (313.821,29.8402) },
"4";"flex17"**\crv{ (231.186,100.41) & (217.868,121.596) },
"flex17";"6"*[o]=(5,5){\,}**\crv{ (200.137,163.542) & (194.63,184.131) },
"6"*[o]=(5,5){\,};"flex18"**\crv{ (199.149,222.848) & (207.046,239.384) },
"flex18";"2"**\crv{ (225.103,270.663) & (235.476,285.776) },
"1"*{\bullet}*!<5pt,9pt>{\hbox{1}},
"2"*{\bullet}*!<0pt,-9pt>{\hbox{2}},
"3"*{\bullet}*!<-7pt,7pt>{\hbox{3}},
"4"*{\bullet}*!<-9pt,0pt>{\hbox{4}},
"5"*{\bullet}*!<-4pt,9pt>{\hbox{5}},
"6"*{\bullet}*!<7pt,7pt>{\hbox{6}},
%
"1";"1"+(20,-15)**\dir{-}*\dir{>},
"1";"1"-(20,-15)**\dir{-}*\dir{>},
"2";"2"+(30,0)**\dir{-}*\dir{>},
"2";"2"-(30,0)**\dir{-}*\dir{>},
"3";"3"+(-17,-25)**\dir{-}*\dir{>},
"3";"3"-(-17,-25)**\dir{-}*\dir{>},
"4";"4"+(0,27)**\dir{-}*\dir{>},
"4";"4"-(0,27)**\dir{-}*\dir{>},
"5";"5"+(20,15)**\dir{-}*\dir{>},
"5";"5"-(20,15)**\dir{-}*\dir{>},
"6";"6"+(17,-25)**\dir{-}*\dir{>},
"6";"6"-(17,-25)**\dir{-}*\dir{>},
\endxy}\\
\table{ll}
 0 & (1,2,3,4,1),(1,5,3,6,1),(2,5,4,6,2)\\
 1 & (1,2,5,1),(1,4,6,1),(2,3,6,2),(3,4,5,3)\\
\endtable
\endtable
&
\table{lll}
\multicolumn{3}{c}{List of o-cycles (7)}\\
(1,2,3,4,1) & 4 & 0\\
(2,5,4,6,2) & 4 & 0\\
(1,5,3,6,1) & 4 & 1\\
(1,4,6,1)   & 3 & 0\\
(1,2,5,1)   & 3 & 1\\
(2,3,6,2)   & 3 & 1\\
(3,4,5,3)   & 3 & 1\\
\endtable
\endtable}

\medskip
\centerline{\table{cc}
\table{c}
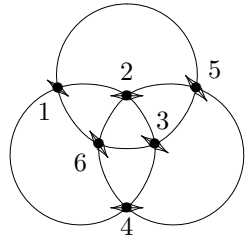
\hbox{\xy /r.20pt/:,
(119.113,311.02)="1";
(191.876,317.16)="flex7"; 
(283.884,255.041)="flex8"; 
(302.887,204.917)="3";
(291.639,144.326)="flex9"; 
(46.7613,118.174)="flex10"; 
(249.972,470.153)="flex11"; 
(380.834,311.02)="5";
(349.774,244.943)="flex12"; 
(249.969,196.31)="flex13"; 
(197.053,204.892)="6";
(150.169,244.915)="flex14"; 
(249.984,296.571)="2";
(308.08,317.151)="flex15"; 
(453.206,118.124)="flex16"; 
(249.947,84.3883)="4";
(208.268,144.307)="flex17"; 
(216.073,255.021)="flex18"; 
"1";"flex7"**\crv{ (142.378,319.265) & (167.397,320.203) },
"flex7";"2"*[o]=(5,5){\,}**\crv{ (212.611,314.583) & (233.207,309.056) },
"2"*[o]=(5,5){\,};"flex8"**\crv{ (264.485,285.78) & (274.858,270.676) },
"flex8";"3"**\crv{ (292.909,239.406) & (300.8,222.871) },
"3";"flex9"**\crv{ (305.301,184.151) & (299.777,163.562) },
"flex9";"4"*[o]=(5,5){\,}**\crv{ (282.03,121.611) & (268.714,100.414) },
"4"*[o]=(5,5){\,};"flex10"**\crv{ (186.091,29.8623) & (88.5278,45.8136) },
"flex10";"1"**\crv{ (5,190.526) & (39.9736,282.974) },
"1"*[o]=(5,5){\,};"flex11"**\crv{ (103.842,393.609) & (166.419,470.147) },
"flex11";"5"**\crv{ (333.535,470.16) & (396.127,393.611) },
"5";"flex12"**\crv{ (376.341,286.753) & (364.647,264.619) },
"flex12";"3"*[o]=(5,5){\,}**\crv{ (337.174,228.275) & (322.087,213.205) },
"3"*[o]=(5,5){\,};"flex13"**\crv{ (286.29,197.752) & (268.024,196.315) },
"flex13";"6"**\crv{ (231.916,196.306) & (213.65,197.733) },
"6";"flex14"**\crv{ (177.851,213.175) & (162.766,228.246) },
"flex14";"1"*[o]=(5,5){\,}**\crv{ (135.293,264.598) & (123.601,286.744) },
"2";"flex15"**\crv{ (266.758,309.053) & (287.35,314.575) },
"flex15";"5"*[o]=(5,5){\,}**\crv{ (332.556,320.193) & (357.571,319.259) },
"5"*[o]=(5,5){\,};"flex16"**\crv{ (460.009,282.979) & (495,190.491) },
"flex16";"4"**\crv{ (411.413,45.7587) & (313.821,29.8402) },
"4";"flex17"**\crv{ (231.186,100.41) & (217.868,121.596) },
"flex17";"6"*[o]=(5,5){\,}**\crv{ (200.137,163.542) & (194.63,184.131) },
"6"*[o]=(5,5){\,};"flex18"**\crv{ (199.149,222.848) & (207.046,239.384) },
"flex18";"2"**\crv{ (225.103,270.663) & (235.476,285.776) },
"1"*{\bullet}*!<5pt,9pt>{\hbox{1}},
"2"*{\bullet}*!<0pt,-9pt>{\hbox{2}},
"3"*{\bullet}*!<-3pt,-9pt>{\hbox{3}},
"4"*{\bullet}*!<0pt,7pt>{\hbox{4}},
"5"*{\bullet}*!<-7pt,-7pt>{\hbox{5}},
"6"*{\bullet}*!<7pt,7pt>{\hbox{6}},
%
"1";"1"+(20,-15)**\dir{-}*\dir{>},
"1";"1"-(20,-15)**\dir{-}*\dir{>},
"2";"2"+(30,0)**\dir{-}*\dir{>},
"2";"2"-(30,0)**\dir{-}*\dir{>},
"3";"3"+(-25,16)**\dir{-}*\dir{>},
"3";"3"-(-25,16)**\dir{-}*\dir{>},
"4";"4"+(-30,0)**\dir{-}*\dir{>},
"4";"4"+(30,0)**\dir{-}*\dir{>},
"5";"5"+(20,-24)**\dir{-}*\dir{>},
"5";"5"-(20,-24)**\dir{-}*\dir{>},
"6";"6"+(17,-25)**\dir{-}*\dir{>},
"6";"6"-(17,-25)**\dir{-}*\dir{>},
\endxy}\\
\table{ll}
0 & (1,2,3,4,1),(1,5,3,6,1),(2,5,4,6,2)\\
1 & (1,2,3,4,1),(1,5,4,6,1),(2,5,3,6,2)\\
2 & (1,2,3,5,1),(1,4,3,6,1),(2,5,4,6,2)\\
3 & (1,2,5,4,1),(1,5,3,6,1),(2,3,4,6,2)\\
\endtable
\endtable
&
\table{lll}
\multicolumn{3}{c}{List of o-cycles (9)}\\
(1,2,3,4,1) & length 4 & \setbox0=\hbox{0}\hbox to \wd0{0,1\hss}\\
(2,5,4,6,2) & length 4 & \setbox0=\hbox{0}\hbox to \wd0{0,2\hss}\\
(1,5,3,6,1) & length 4 & \setbox0=\hbox{0}\hbox to \wd0{0,3\hss}\\
(1,5,4,6,1) & length 4 & 1\\
(2,5,3,6,2) & length 4 & 1\\
(1,2,3,5,1) & length 4 & 2\\
(1,4,3,6,1) & length 4 & 2\\
(1,2,5,4,1) & length 4 & 3\\
(2,3,4,6,2) & length 4 & 3\\
\endtable
\endtable}

\medskip
\centerline{\table{cc}
\table{c}
\hbox{\xy /r.20pt/:,
(119.113,311.02)="1";
(191.876,317.16)="flex7"; 
(283.884,255.041)="flex8"; 
(302.887,204.917)="3";
(291.639,144.326)="flex9"; 
(46.7613,118.174)="flex10"; 
(249.972,470.153)="flex11"; 
(380.834,311.02)="5";
(349.774,244.943)="flex12"; 
(249.969,196.31)="flex13"; 
(197.053,204.892)="6";
(150.169,244.915)="flex14"; 
(249.984,296.571)="2";
(308.08,317.151)="flex15"; 
(453.206,118.124)="flex16"; 
(249.947,84.3883)="4";
(208.268,144.307)="flex17"; 
(216.073,255.021)="flex18"; 
"1";"flex7"**\crv{ (142.378,319.265) & (167.397,320.203) },
"flex7";"2"*[o]=(5,5){\,}**\crv{ (212.611,314.583) & (233.207,309.056) },
"2"*[o]=(5,5){\,};"flex8"**\crv{ (264.485,285.78) & (274.858,270.676) },
"flex8";"3"**\crv{ (292.909,239.406) & (300.8,222.871) },
"3";"flex9"**\crv{ (305.301,184.151) & (299.777,163.562) },
"flex9";"4"*[o]=(5,5){\,}**\crv{ (282.03,121.611) & (268.714,100.414) },
"4"*[o]=(5,5){\,};"flex10"**\crv{ (186.091,29.8623) & (88.5278,45.8136) },
"flex10";"1"**\crv{ (5,190.526) & (39.9736,282.974) },
"1"*[o]=(5,5){\,};"flex11"**\crv{ (103.842,393.609) & (166.419,470.147) },
"flex11";"5"**\crv{ (333.535,470.16) & (396.127,393.611) },
"5";"flex12"**\crv{ (376.341,286.753) & (364.647,264.619) },
"flex12";"3"*[o]=(5,5){\,}**\crv{ (337.174,228.275) & (322.087,213.205) },
"3"*[o]=(5,5){\,};"flex13"**\crv{ (286.29,197.752) & (268.024,196.315) },
"flex13";"6"**\crv{ (231.916,196.306) & (213.65,197.733) },
"6";"flex14"**\crv{ (177.851,213.175) & (162.766,228.246) },
"flex14";"1"*[o]=(5,5){\,}**\crv{ (135.293,264.598) & (123.601,286.744) },
"2";"flex15"**\crv{ (266.758,309.053) & (287.35,314.575) },
"flex15";"5"*[o]=(5,5){\,}**\crv{ (332.556,320.193) & (357.571,319.259) },
"5"*[o]=(5,5){\,};"flex16"**\crv{ (460.009,282.979) & (495,190.491) },
"flex16";"4"**\crv{ (411.413,45.7587) & (313.821,29.8402) },
"4";"flex17"**\crv{ (231.186,100.41) & (217.868,121.596) },
"flex17";"6"*[o]=(5,5){\,}**\crv{ (200.137,163.542) & (194.63,184.131) },
"6"*[o]=(5,5){\,};"flex18"**\crv{ (199.149,222.848) & (207.046,239.384) },
"flex18";"2"**\crv{ (225.103,270.663) & (235.476,285.776) },
"1"*{\bullet}*!<5pt,-7pt>{\hbox{1}},
"2"*{\bullet}*!<0pt,-9pt>{\hbox{2}},
"3"*{\bullet}*!<-7pt,7pt>{\hbox{3}},
"4"*{\bullet}*!<0pt,9pt>{\hbox{4}},
"5"*{\bullet}*!<-7pt,-7pt>{\hbox{5}},
"6"*{\bullet}*!<7pt,7pt>{\hbox{6}},
"1";"1"+(13,21)**\dir{-}*\dir{>},
"1";"1"-(13,21)**\dir{-}*\dir{>},
"2";"2"+(30,0)**\dir{-}*\dir{>},
"2";"2"-(30,0)**\dir{-}*\dir{>},
"3";"3"+(-17,-25)**\dir{-}*\dir{>},
"3";"3"-(-17,-25)**\dir{-}*\dir{>},
"4";"4"+(27,0)**\dir{-}*\dir{>},
"4";"4"-(27,0)**\dir{-}*\dir{>},
"5";"5"+(20,-24)**\dir{-}*\dir{>},
"5";"5"-(20,-24)**\dir{-}*\dir{>},
"6";"6"+(17,-25)**\dir{-}*\dir{>},
"6";"6"-(17,-25)**\dir{-}*\dir{>},
\endxy}\\
\table{ll}
0 & (1,2,3,4,1),(1,5,3,6,1),(2,5,4,6,2)\\
1 & (1,2,3,4,1),(1,5,4,6,1),(2,5,3,6,2)\\
2 & (1,2,3,4,6,1),(1,4,5,1),(2,5,3,6,2)\\
3 & (1,2,3,6,1),(1,4,3,5,1),(2,5,4,6,2)\\
4 & (1,2,3,6,1),(1,4,5,1),(2,5,3,4,6,2)\\
5 & (1,2,5,3,4,1),(1,5,4,6,1),(2,3,6,2)\\
6 & (1,2,5,3,4,6,1),(1,4,5,1),(2,3,6,2)\\
7 & (1,2,5,3,6,1),(1,4,5,1),(2,3,4,6,2)\\
8 & (1,2,5,4,1),(1,5,3,6,1),(2,3,4,6,2)\\
9 & (1,2,5,4,1),(1,5,3,4,6,1),(2,3,6,2)\\
10 & (1,2,5,4,6,1),(1,4,3,5,1),(2,3,6,2)\\
\endtable
\endtable
&
\table{lll}
\multicolumn{3}{c}{List of o-cycles (18)}\\
(1,2,5,3,4,6,1) & 6 & \setbox0=\hbox{0}\hbox to \wd0{0,1\hss}\\
(1,2,3,4,6,1)   & 5 & 2\\
(1,2,5,3,4,1)   & 5 & \setbox0=\hbox{0}\hbox to \wd0{3,4\hss}\\
(1,2,5,3,6,1)   & 5 & 5\\
(1,2,5,4,6,1)   & 5 & 6\\
(2,5,3,4,6,2)   & 5 & 7\\
(1,5,3,4,6,1)   & 5 & \setbox0=\hbox{0}\hbox to \wd0{8,9\hss}\\
(1,2,3,4,1)     & 4 & \setbox0=\hbox{0}\hbox to \wd0{0,3\hss}\\
(2,5,4,6,2)     & 4 & \setbox0=\hbox{0}\hbox to \wd0{0,8\hss}\\
(1,2,3,6,1)     & 4 & \setbox0=\hbox{0}\hbox to \wd0{1,2\hss}\\
(1,4,3,5,1)     & 4 & \setbox0=\hbox{0}\hbox to \wd0{1,5\hss}\\
(1,5,3,6,1)     & 4 & \setbox0=\hbox{0}\hbox to \wd0{2,4,6,7\hss}\\
(1,5,4,6,1)     & 4 & \setbox0=\hbox{0}\hbox to \wd0{3,10\hss}\\
(1,2,5,4,1)     & 4 & 4\\
(2,5,3,6,2)     & 4 & 9\\
(2,3,4,6,2)     & 4 & 10\\
(2,3,6,2)       & 3 & \setbox0=\hbox{0}\hbox to \wd0{5,6,9,10\hss}\\
(1,4,5,1)       & 3 & \setbox0=\hbox{0}\hbox to \wd0{7,8\hss}\\
\endtable
\endtable}

\end{document}